\documentclass{article}%
\usepackage[style=alphabetic,maxnames=99,maxalphanames=4,maxbibnames=10]{biblatex}
\usepackage[final]{microtype}%
\usepackage{amsthm,mathtools}%
\usepackage{xcolor}%
\usepackage{float}%
\usepackage[export]{adjustbox}%
\usepackage{tikz, tikz-cd, mathtools, amssymb, stmaryrd}\usetikzlibrary{calc}
\usetikzlibrary{decorations.pathmorphing}

\tikzset{curve/.style={settings={#1},to path={(\tikztostart)
    .. controls ($(\tikztostart)!\pv{pos}!(\tikztotarget)!\pv{height}!270:(\tikztotarget)$)
    and ($(\tikztostart)!1-\pv{pos}!(\tikztotarget)!\pv{height}!270:(\tikztotarget)$)
    .. (\tikztotarget)\tikztonodes}},
    settings/.code={\tikzset{quiver/.cd,#1}
        \def\pv##1{\pgfkeysvalueof{/tikz/quiver/##1}}},
    quiver/.cd,pos/.initial=0.35,height/.initial=0}

\tikzset{tail reversed/.code={\pgfsetarrowsstart{tikzcd to}}}
\tikzset{2tail/.code={\pgfsetarrowsstart{Implies[reversed]}}}
\tikzset{2tail reversed/.code={\pgfsetarrowsstart{Implies}}}
\tikzset{no body/.style={/tikz/dash pattern=on 0 off 1mm}}
    \usepackage{mathpartir}
\newtheorem{theorem}{Theorem}[section]%
\newtheorem{lemma}[theorem]{Lemma}%
\newtheorem{observation}[theorem]{Observation}%
\newtheorem{corollary}[theorem]{Corollary}%
\theoremstyle{definition}%
\newtheorem{definition}[theorem]{Definition}%
\newtheorem{informaldefinition}[theorem]{Informal definition}%
\newtheorem{construction}[theorem]{Construction}%
\newtheorem{notation}[theorem]{Notation}%
\newtheorem{explication}[theorem]{Explication}%
\newtheorem{proposition}[theorem]{Proposition}%
\newtheorem{attitude}[theorem]{Attitude}%
\newtheorem{example}[theorem]{Example}%
\newtheorem{remark}[theorem]{Remark}%
\newtheorem{convention}[theorem]{Convention}%
\newtheorem{acknowledgements}{Acknowledgements}%
\usepackage[colorlinks=true,linkcolor={blue!30!black}]{hyperref}%
\usepackage{newpxmath,newpxtext}%
\usepackage{cleveref}%
\usepackage[mode=buildmissing]{standalone}%
\usepackage{hyperref}
\usepackage{url}

    \title{Towards a double operadic theory of systems}\author{
        Sophie Libkind
      \and{}
        David Jaz Myers
      }
\begin{document}
\maketitle
    \begin{abstract}
        We present a unified framework for categorical systems theory which packages a collection of open systems, their interactions, and their maps into a symmetric monoidal loose right \textbf{module of systems} over a symmetric monoidal \textbf{double category of interfaces and interactions}. As examples, we give detailed descriptions of (1) the module of open Petri nets over  undirected wiring diagrams and (2) the module of deterministic Moore machines over lenses. We define several  pseudo-functorial constructions of modules of systems in the form of \textbf{doctrines of systems theories}. In particular, we introduce doctrines for port-plugging systems, variable sharing systems, and generalized Moore machines, each of which  generalizes existing work in categorical systems theory. Finally, we observe how diagrammatic interaction patterns are free processes in particular doctrines. 
    \end{abstract}

    \tableofcontents
    \section{Introduction}\label{djm-009H}\subsection{Categorical systems theories}\label{ssl-002T}\par{}Categorical systems theory is the branch of applied category theory that uses methods of categorical algebra to aid in the \emph{modular} design and \emph{compositional} analysis of complex systems.\par{}No one can quite agree on what a ``system'' is; and they shouldn't have to. Rather, different people working on different problems devise different particular notions of system --- \textbf{systems theories} --- to best address the needs of the problem at hand. ``Systems'' might be systems of equations such as ODEs or PDEs, (non-)deterministic automata of varying forms, Markov processes or Markov decision processes, Hamiltonians on symplectic manifolds, Lagrangian relations, or a variety of diagrammatic languages used to describe systems such as circuit diagrams, stock-flow diagrams, Petri nets, Tonti diagrams, labelled transition systems, weighted graphs, and many others.\par{}	In general, the term ``system'' only means that a whole has been \emph{composed} of many parts. We use the term \emph{categorical systems theory} to refer to a large body of work that uses categorical algebra to organize the patterns of composition by which complex systems may be formed out of simpler components, and giving methods for analysing the behavior of composite systems in terms of their components' behavior and the pattern by which they were composed. The basic ideas of categorical systems theory are:\subsubsection{Basic ideas of categorical systems theory}\label{djm-00FO}\begin{enumerate}\item{}Any system interacts with its environment through its \textbf{interface}, which can be described separately from the system itself.
	\item{}All interactions of the system with its environment take place through its interface, so that from the point of view of the environment, all we need to know about a system is what is going on at the interface.
	\item{}Systems interact with other systems through their respective interfaces through ongoing \textbf{processes}. So, to understand complex systems in terms of their component subsystems, we need to understand the processes through which interfaces can interact. We refer to a particularly simple (for example, state-free) interaction as a \textbf{composition pattern}.
	\item{}Given (a) a process describing how system interfaces may interact and (b) component systems with those interfaces, we should derive a \textbf{composite} system that is the totality of the component subsystems interacting according to that process. This application of processes to component systems allows for the \textbf{modular} design of systems.
	\item{}For a given systems theory, we can derive certain behaviors, facts, and features of a composite system from the behaviors, facts, and features of its component subsystems and the process that connects them. Such a method is known as a \textbf{compositionality theorem}. It is not always possible to characterize the behaviors of composite systems in terms of their component subsystems, but we often can compose the behaviors of composite systems to derive \emph{some} of the behaviors of their composites, and then study the difference.\end{enumerate}\par{}	In this paper, we put forward a general framework for categorical systems theory in the form of \emph{symmetric monoidal loose right modules} (of systems) over \emph{symmetric monoidal double categories} (of interfaces and interactions). This framework builds on the ``double indexed'' approach \cite{jaz-2021-double} and its expansion in the manuscript \cite{jaz-2021-book}, as well as the operadic theory of systems \cite{libkind-2022-operadic}. We aim to kick off a program of organizing categorical systems theory in this \emph{double operadic} point of view by:
	\begin{enumerate}\item{}			Constructing, pseudo-functorially in basic data, a number of symmetric monoidal loose right modules of systems of various kinds, including 
			\begin{enumerate}\item{}					Systems that compose via lens composition, such as systems of ODEs, partially observable Markov decision processes, and Moore machines of various kinds.
				\item{}					Systems that compose by gluing together parts (pushout), such as Petri nets, stock-flow diagrams, causal loop diagrams, and other diagrammatic presentations of systems.
				\item{}					Systems that compose by sharing variables, as in Jan Willems' \emph{behavioral} approach to control theory \cite{willems-2007-behavioral}, or Lagrangian relations.\end{enumerate}
		\item{}			Giving a method for restricting the interactions between systems to those generated by particular data, allowing us to extend the above pseudo-functorial constructions to produce algebras for (double) operads of various sorts of \emph{wiring diagrams}.\end{enumerate}
			Pseudo-functoriality of the above constructions can be used to give a number of black-boxing compositionality theorems along the lines of Fong and Sarazola's \emph{recipes} \cite{fong-2018-recipe}.\par{}	In the remainder of this introduction, we will briefly review existing approaches to categorical systems theory and where categorial systems theory would benefit from our ``double operadic'' approach  (in \hyperref[ssl-002U]{Section \ref{ssl-002U}}); we will then sketch, informally, the double operadic view on systems theory (in \hyperref[djm-009L]{Section \ref{djm-009L}}); finally, we will summarize the contributions of this paper (in \hyperref[ssl-002W]{Section \ref{ssl-002W}}) and sketch some future directions for the double operadic theory of systems (in \hyperref[ssl-002X]{Section \ref{ssl-002X}}).\subsection{Operadic systems theories and process theories}\label{ssl-002U}\par{}	Category theorists make the informal ideas of \hyperref[djm-00FO]{§ \ref{djm-00FO}} formal by associating to each systems theory (ODEs, Markov decision processes, causal loop diagrams, etc.) a categorical structure, and expressing \emph{compositionality theorems} as \emph{morphisms} between these structures. Speaking generally, categorical systems theory has relied on two key gadgets: \textbf{algebras for operads}  and \textbf{symmetric monoidal categories}, perhaps with extra structure. 
	\begin{itemize}\item{}			In the \emph{operadic} point of view (see, e.g. \cite{foley-2021-operads}, \cite{libkind-2022-operadic}, \cite{spivak-2013-operad}, \cite{baez-2021-keystones}), we use two mathematical structures to define a systems theory:
			\begin{enumerate}\item{}					A (symmetric, multi-object) \textbf{operad} \(W\) whose objects are system interfaces and whose morphisms \(w : I_1, \ldots , I_n \to  J\) are ways that \(n\) systems (each with interface \(I_j\) for \(0 \leq  j \leq  n\)) may be composed into a single system with interface \(J\) (often known as \emph{wiring diagrams} \cite{yau-2018-operads});				
				\item{}					An \textbf{algebra} \(S \colon  W \to  \mathsf {Set}\) which associates to every interface \(I \in  W\) the set \(S(I)\) of systems with this interface, and associates to every composition pattern \(w : I_1, \ldots , I_n \to  J\) the operation \(S(w) \colon   S(I_1) \times  \cdots  \times  S(I_n) \to  S(J)\) of composition by that pattern.\end{enumerate} 
				A compositionality theorem then takes the form of a morphism between the operads of composition patterns and, relative to this, a morphism of the algebras of systems. Examples of this pattern include \cite{schultz-2019-dynamical}, \cite{spivak-2015-steady}, \cite{dangelo-2025-dependent}, \cite{libkind-2022-algebraic}, and \cite{breiner-2020-modeling}.
		\item{}			In the \emph{symmetric monoidal} point of view (see, e.g. \cite{baez-2010-physics}, \cite{baez-2014-categories}, \cite{baez-2017-props}, \cite{fritz-2020-synthetic}, \cite{coecke-2016-mathematical}), we use a symmetric monoidal category \(S\) (perhaps with extra structure) to define a systems theory. 
			\begin{enumerate}\item{}				It is the morphisms of \(S\) which are the systems in this theory; a system \(s : I \to  O\) in this point of view always has an input interface \(I\) and an output interface \(O\).	
				\item{}Systems compose in series using composition of morphisms in \(S\), and they compose in parallel using the symmetric monoidal product of \(S\).
				\item{}					Since systems in this point of view always have an input and output, we will follow the seminal \cite{katis-1997-bicategories} and call such symmetric monoidal categories \textbf{process theories} (to separate them from the operadic \emph{systems theories}).\end{enumerate}
			A compositionality theorem takes the form of a \emph{lax symmetric monoidal functor} between symmetric monoidal categories. Examples of this pattern include \cite{fong-2018-recipe}, \cite{baez-2016-compositional}, \cite{baez-2015-compositional}, \cite{baez-2017-props}.\end{itemize}\par{}Both operad algebras and symmetric monoidal categories miss out on a crucial way that systems relate to each other: \textbf{maps between systems}.  
	 In systems theory, maps between systems include:
	\begin{enumerate}\item{}\emph{Refinements}, \emph{coarse grainings}, and \emph{simulations}.
		\item{}\emph{Subsystem inclusions} such as \emph{attractors}.
		\item{}			Behaviors of systems such as \emph{trajectories}, \emph{traces}, as well as behaviors satisfying special properties such as \emph{periodic orbits} and \emph{steady states}.
		\item{}			Lyapunov and control-barrier functions which witness certain properties of the behavior of systems.			\end{enumerate}
	Our theories of systems must account for maps between systems in order to study system behavior, express notions of refinement, coarse-graining, and bisimulation which are needed for model surrogacy, and witness safety and liveness properties of systems.\par{}	Various approachs to categorical systems theory have successfully accounted for maps between systems. Indeed, much of categorical systems theory has focused primarily on systems and their maps, and not on the composition of the systems themselves. The highly successful application of \emph{coalgebras} (\cite{rutten-2000-universal}, \cite{kurz-2001-coalgebras}, \cite{cirstea-2013-from}) to systems theory exemplifies how the study of maps between systems can capture the structure of the systems themselves.\par{}Baez et. al. \cite{baez-2017-network} extended the operadic systems theories to account for maps between systems by using operad algebras valued in \emph{categories} of systems, rather than sets of systems. However, this approach only expresses maps between systems with identical interfaces; many interesting maps, such as trajetories, steady states, and Lyapunov functions must act on the interface as well.\par{}	The seminal work of Katis, Sabadini, and Walters \cite{katis-1997-bicategories}, introduced \emph{symmetric monoidal bicategories} as \emph{process theories} (see also \cite{gadducci-1999-bicategorical}, \cite{boccali-2023-bicategories}, \cite{capucci-2022-towards}, \cite{di-2020-spangraph}). This line of work expressed maps between processes as well as sequential and parallel composition of processes.  However, as with the operadic systems theories, this approach does not express maps between systems with different interaces.\par{}In his thesis \cite{courser-2020-open}, Courser further extends the process theory lineage by considering \emph{symmetric monoidal double categories} of systems (see also \cite{baez-2017-coarsegraining}, \cite{baez-2022-structured}, \cite{loregian-2025-monads}, \cite{master-2021-composing}, \cite{baez-2020-open}). The upgrade to double categories introduces maps between processes with different interfaces. However, the process theory approach obscures the composition patterns that interconnect systems, and presupposes that system interfaces have a clean division into input and output. See \hyperref[djm-00FR]{Remark \ref{djm-00FR}} for a more detailed comparison between our approach and Courser's.\par{}In this paper, we put forward a new approach categorical systems theory using \emph{algebras for double operads}, or, more precisely, \emph{symmetric monoidal loose right modules} of symmetric monoidal double categories. We will call our approach the \textbf{double operadic theory of systems} (or ``DOTS''). Double operadic systems theory combines the best of the operadic systems theory and process theory approaches. As in operadic systems theories, a double operadic systems theory treats composition patterns as concrete objects that can act on systems.  As in Courser's processes theories, a double operadic systems theory includes maps between systems with different interfaces.\subsection{Double operadic systems theory}\label{djm-009L}\par{}The \hyperref[djm-00FO]{basic ideas of categorical systems theory} leave a lot of questions unanswered. What is, or what could be a system? What is an interface? What is a process of interactions? between systems? What is a behavior of a system, and how can we study it categorically? There is no single answer to this suite of questions. Rather, we may package answers to this suite of questions into a \textbf{systems theory}.
\begin{informaldefinition}[{Systems theory}]\label{djm-009K}
\par{}A \textbf{systems theory} is a particular way to answer the following questions:
\begin{enumerate}\item{}What does it mean to be a \textbf{system}? Does it have a notion of states, or of behaviors? Or is it a diagram describing the way some primitive parts are organized?    
\item{}What should the \textbf{interface} of a system be?  
\item{}How can interfaces be connected through \textbf{processes of interaction} or \textbf{composition patterns}?  
\item{}How are systems \textbf{composed} through processes of interaction between their interfaces?  
\item{}How do systems \emph{refine}, \emph{simulate}, or \emph{coarse-grain} each other? How may systems be included as \emph{sub-systems}? In other words, what is a \textbf{map} between systems, and what is a corresponding notion of map for interfaces? 
\item{}When can maps between systems be composed along the same interactions as the systems?\end{enumerate}\end{informaldefinition}
\par{}The questions in \hyperref[djm-009K]{Informal definition \ref{djm-009K}} suggests the following ontology for systems theory, which we may package into a very small (strict) double category: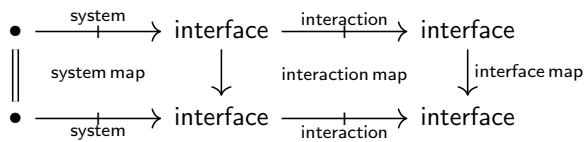
\begin{figure}[H]\begin{center}
\begin {tikzcd}[column sep = huge]
	\bullet  & {\mathsf {interface}} & {\mathsf {interface}} \\
	\bullet  & {\mathsf {interface}} & {\mathsf {interface}}
	\arrow [""{name=0, anchor=center, inner sep=0}, "{\mathsf {system}}", "\shortmid "{marking}, from=1-1, to=1-2]
	\arrow [Rightarrow, no head, from=1-1, to=2-1]
	\arrow [""{name=1, anchor=center, inner sep=0}, "{\mathsf {interaction}}", "\shortmid "{marking}, from=1-2, to=1-3]
	\arrow [from=1-2, to=2-2]
	\arrow ["{\mathsf {interface\, map}}", from=1-3, to=2-3]
	\arrow [""{name=2, anchor=center, inner sep=0}, "{\mathsf {system}}"', "\shortmid "{marking}, from=2-1, to=2-2]
	\arrow [""{name=3, anchor=center, inner sep=0}, "{\mathsf {interaction}}"', "\shortmid "{marking}, from=2-2, to=2-3]
	\arrow ["{\mathsf {system\, map}}"{description}, draw=none, from=0, to=2]
	\arrow ["{\mathsf {interaction\, map}}"{description}, draw=none, from=1, to=3]
\end {tikzcd}
\end{center}\caption{Ontology of double categorical systems theory}\label{djm-009N}\end{figure}\par{}We read the elements of this double category as the kinds of entities determined by theory of systems: \(\mathsf {system}\)s, \(\mathsf {interfaces}\), \(\mathsf {interactions}\), and so on. The extra unnamed object \(\bullet \) is there just to give \(\mathsf {system}\) a domain. If we suppose that all displayed elements except for \(\mathsf {system}\) are identities in this double category, we get a very compact way of expressing the following axioms for systems theories:
\begin{enumerate}\item{}Every \textbf{system} admits an \textbf{interface}: we have \(\mathsf {system} : \bullet  \mathrel {\mkern 3mu\vcenter {\hbox {$\shortmid $}}\mkern -10mu{\to }} \mathsf {interface}\).
	\item{}An \textbf{interaction} has an "inner" interface and an "outer" interface: we have \(\mathsf {interaction} : \mathsf {interface} \mathrel {\mkern 3mu\vcenter {\hbox {$\shortmid $}}\mkern -10mu{\to }} \mathsf {interface}\).
	\item{}We may compose a system with a interaction to get a new system:
		
\begin{center}
			\begin {tikzcd}
		\bullet  & {\mathsf {interface}} & {\mathsf {interface}}
		\arrow ["{\mathsf {system}}"', "\shortmid "{marking}, from=1-1, to=1-2]
		\arrow ["{\mathsf {system}}", "\shortmid "{marking}, bend left, from=1-1, to=1-3]
		\arrow ["{\mathsf {interaction}}"', "\shortmid "{marking}, from=1-2, to=1-3]
	\end {tikzcd}
		\end{center}
	\item{}There is a notion of \textbf{map} for interfaces (\(\mathsf {interface\,map} : \mathsf {interface} \to  \mathsf {interface}\)), and a notion of map for systems (\(\mathsf {system\, map}\)); every system map comes with a corresponding map from the interface of its domain to the interface of its codomain. We may compose interface maps to get new interface maps (since \(\mathsf {interface\, map} \circ  \mathsf {interface\, map} = \mathsf {interface\, map}\)) and similarly we may compose system maps to get new system maps.	
	
\begin{center}
		\begin {tikzcd}
			\bullet  & {\mathsf {interface}} && \bullet  & {\mathsf {interface}} \\
			\bullet  & {\mathsf {interface}} & {=} \\
			\bullet  & {\mathsf {interface}} && \bullet  & {\mathsf {interface}}
			\arrow [""{name=0, anchor=center, inner sep=0}, "{\mathsf {system}}", "\shortmid "{marking}, from=1-1, to=1-2]
			\arrow [equals, from=1-1, to=2-1]
			\arrow ["{\mathsf {interface\, map}}", from=1-2, to=2-2]
			\arrow [""{name=1, anchor=center, inner sep=0}, "{\mathsf {system}}", "\shortmid "{marking}, from=1-4, to=1-5]
			\arrow [equals, from=1-4, to=3-4]
			\arrow ["{\mathsf {interface\, map}}", from=1-5, to=3-5]
			\arrow [""{name=2, anchor=center, inner sep=0}, "\shortmid "{marking}, from=2-1, to=2-2]
			\arrow [equals, from=2-1, to=3-1]
			\arrow ["{\mathsf {interface\, map}}", from=2-2, to=3-2]
			\arrow [""{name=3, anchor=center, inner sep=0}, "{\mathsf {system}}"', "\shortmid "{marking}, from=3-1, to=3-2]
			\arrow [""{name=4, anchor=center, inner sep=0}, "{\mathsf {system}}"', "\shortmid "{marking}, from=3-4, to=3-5]
			\arrow ["{\mathsf {system\,map}}"{description}, shorten <=4pt, shorten >=4pt, Rightarrow, from=0, to=2]
			\arrow ["{\mathsf {system\, map}}"{description}, shorten <=9pt, shorten >=9pt, Rightarrow, from=1, to=4]
			\arrow ["{\mathsf {system\,map}}"{description}, shorten <=4pt, shorten >=4pt, Rightarrow, from=2, to=3]
		\end {tikzcd}
	\end{center}
	\item{}There is a notion of \textbf{map of interactions}, which we might also think of as a \textbf{compatibility} between interface maps and composition patterns. Given a map of interactions, we may compose a system map along those interaction maps to get a new system map between the composed systems:
	
\begin{center}
		\begin {tikzcd}
	\bullet  & {\mathsf {interface}} & {\mathsf {interface}} \\
	\bullet  & {\mathsf {interface}} & {\mathsf {interface}} \\
	& {=} \\
	\bullet  && {\mathsf {interface}} \\
	\bullet  && {\mathsf {interface}}
	\arrow [""{name=0, anchor=center, inner sep=0}, "{\mathsf {system}}", "\shortmid "{marking}, from=1-1, to=1-2]
	\arrow [equals, from=1-1, to=2-1]
	\arrow [""{name=1, anchor=center, inner sep=0}, "{\mathsf {interaction}}", "\shortmid "{marking}, from=1-2, to=1-3]
	\arrow [from=1-2, to=2-2]
	\arrow ["{\mathsf {interface\, map}}", from=1-3, to=2-3]
	\arrow [""{name=2, anchor=center, inner sep=0}, "{\mathsf {system}}"', "\shortmid "{marking}, from=2-1, to=2-2]
	\arrow [""{name=3, anchor=center, inner sep=0}, "{\mathsf {interaction}}"', "\shortmid "{marking}, from=2-2, to=2-3]
	\arrow [""{name=4, anchor=center, inner sep=0}, "\shortmid "{marking}, from=4-1, to=4-3]
	\arrow [equals, from=4-1, to=5-1]
	\arrow ["{\mathsf {interface\, map}}", from=4-3, to=5-3]
	\arrow [""{name=5, anchor=center, inner sep=0}, "\shortmid "{marking}, from=5-1, to=5-3]
	\arrow ["{\mathsf {system\,map}}"{description}, shorten <=4pt, shorten >=4pt, Rightarrow, from=0, to=2]
	\arrow ["{\mathsf {ineraction\,map}}"{description}, shorten <=4pt, shorten >=4pt, Rightarrow, from=1, to=3]
	\arrow ["{\mathsf {system\,map}}"{description}, shorten <=4pt, shorten >=4pt, Rightarrow, from=4, to=5]
\end {tikzcd}
	\end{center}\end{enumerate}\par{}In general, we may think of the double category \hyperref[djm-009N]{Figure \ref{djm-009N}} as our ontology for systems theory. A particular systems theory is the necessary assumptions to form a collection of entities of the above kinds, able to compose in the above ways. In other words, a systems theory will give rise to a (generally not strict) double category \(\mathbb {S}\) equipped with a double functor to the \hyperref[djm-009N]{above simple double category}, \emph{labelling} each of its elements with the appropriate kinds. We call such a labeled double category a \emph{module of systems}, and we may picture a module of systems as follows: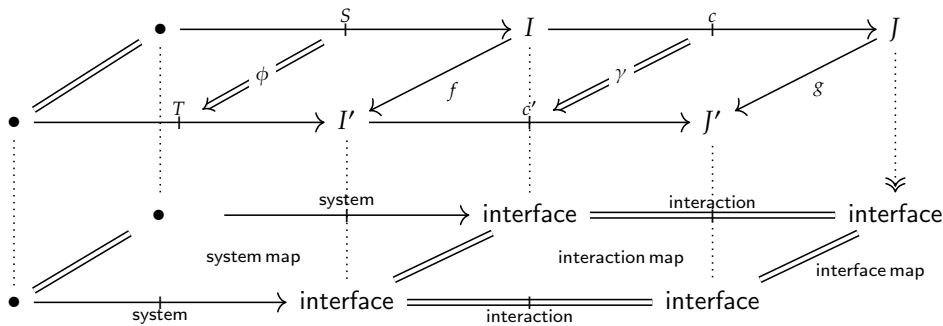
\begin{figure}[H]\begin{center}
\begin {tikzcd}
	& \bullet  && I && J \\
	\bullet  && {I'} && {J'} \\
	& {\phantom {int}\bullet \phantom {ace}} && {\mathsf {interface}} && {\mathsf {interface}} \\
	\bullet  && {\mathsf {interface}} && {\mathsf {interface}}
	\arrow [""{name=0, anchor=center, inner sep=0}, "S", "\shortmid "{marking}, from=1-2, to=1-4]
	\arrow [Rightarrow, no head, from=1-2, to=2-1]
	\arrow [dotted, no head, from=1-2, to=3-2]
	\arrow [""{name=1, anchor=center, inner sep=0}, "c", "\shortmid "{marking}, from=1-4, to=1-6]
	\arrow ["f", from=1-4, to=2-3]
	\arrow [dotted, no head, from=1-4, to=3-4]
	\arrow ["g", from=1-6, to=2-5]
	\arrow [dotted, two heads, from=1-6, to=3-6]
	\arrow [""{name=2, anchor=center, inner sep=0}, "T", "\shortmid "{marking}, from=2-1, to=2-3]
	\arrow [dotted, no head, from=2-1, to=4-1]
	\arrow [""{name=3, anchor=center, inner sep=0}, "{c'}", "\shortmid "{marking}, from=2-3, to=2-5]
	\arrow [dotted, no head, from=2-3, to=4-3]
	\arrow [dotted, no head, from=2-5, to=4-5]
	\arrow [""{name=4, anchor=center, inner sep=0}, "{\mathsf {system}}", "\shortmid "{marking}, from=3-2, to=3-4]
	\arrow [Rightarrow, no head, from=3-2, to=4-1]
	\arrow [""{name=5, anchor=center, inner sep=0}, "{\mathsf {interaction}}", "\shortmid "{marking}, Rightarrow, no head, from=3-4, to=3-6]
	\arrow [Rightarrow, no head, from=3-4, to=4-3]
	\arrow ["{\mathsf {interface\, map}}", Rightarrow, no head, from=3-6, to=4-5]
	\arrow [""{name=6, anchor=center, inner sep=0}, "{\mathsf {system}}"', "\shortmid "{marking}, from=4-1, to=4-3]
	\arrow [""{name=7, anchor=center, inner sep=0}, "{\mathsf {interaction}}"', "\shortmid "{marking}, Rightarrow, no head, from=4-3, to=4-5]
	\arrow ["\phi "{description}, shorten <=9pt, shorten >=9pt, Rightarrow, from=0, to=2]
	\arrow ["\gamma "{description}, shorten <=9pt, shorten >=9pt, Rightarrow, from=1, to=3]
	\arrow ["{\mathsf {system\, map}}"{description}, draw=none, from=4, to=6]
	\arrow ["{\mathsf {interaction\, map}}"{description}, draw=none, from=5, to=7]
\end {tikzcd}
\end{center}\caption{Systems theory as labelled double category}\label{djm-009O}\end{figure}\par{}Not quite pictured above is the most basic kind of composition of systems: the \textbf{parallel product}, where two systems are composed but \emph{do not interact}. We write the paralell product as \(\parallel \) for both systems and interfaces, so that in the above picture we have \(\mathsf {system} \parallel  \mathsf {system} = \mathsf {system}\) (that is, the parallel product of two systems is a system) and \(\mathsf {interface} \parallel  \mathsf {interface} = \mathsf {interface}\) (the parallell product of two interfaces is an interface), and so on for the other elements of a module of systems. This leads us to our formal definition in \hyperref[ssl-0047]{Section \ref{ssl-0047}} of a \textbf{module of systems} as a \emph{symmetric monoidal object} of an appropriate 2-category of \emph{loose (right) modules of double categories}: a symmetric monoidal loose right module (of systems and maps) over a symmetric monoidal double category (of interfaces, interactions, and maps between them).\par{}We will see this point of view on systems theories worked out with full examples in \hyperref[ssl-0013]{Section \ref{ssl-0013}}.
\begin{remark}[{Symmetric monoidal loose right modules and algebras of double operads}]\label{djm-00H9}
\par{}	One might wonder why we have decided to use the name "double operadic systems theory" when our main objects of study are not algebras for double operads but instead symmetric monoidal loose right modules on symmetric monoidal double categories.\par{}	In the 1-categorical case, symmetric monoidal categories \(\mathsf {C}\) may be identified with operads (here meaning symmetric multicategories) which have all tensors. In this case, algebras for \(\mathsf {C}\) (which may be defined as operad morphisms \(A : \mathsf {C} \to  \mathsf {Set}\)) correspond to lax symmetric monoidal functors \(A : \mathsf {C} \to  \mathsf {Set}\). Such functors are \emph{symmetric monoidal right modules} (also known as profunctors) over the symmetric monoidal category \(\mathsf {C}\).\par{}	To facilitate the constructions we perform in this paper, we make use of the same trick. We will use symmetric monoidal \emph{loose} right modules over a symmetric monoidal double category in place of algebras for double operads. We put off a formal comparison between these two notions to later work. However, it is true that in the double categorical setting, these two notions diverge more strongly than in the 1-categorical setting. Nevertheless, the constructions we present here suffice to give algebras for double operads.\par{}	For the expert reader, here is a specific conjecture concerning their relationship: the 2-category of symmetric monoidal loose right modules over isofibrant double categories and \emph{pseudo symmetric monoidal pseudo double functors} whose unitors and laxators are \hyperref[djm-00EX]{\emph{companion commuter transformations}} is equivalent to the 2-category of algebras for isofibrant double operads with all tensors and pseudo-morphisms whose universal comparison map between tensors is companion to an isomorphism.\end{remark}
\subsection{Contributions of this paper}\label{ssl-002W}\par{}	In this paper, we put forward a general setting for the  modular design and compositional analysis of complex systems that accounts for both the composition of systems and of the maps between them. Below we present several key aspects of our approach to categorical systems theory and their relation to past work in categorical systems theory. 

	\begin{enumerate}\item{}We define a \hyperref[ssl-0047]{\textbf{module of systems}} as \textbf{symmetric monoidal loose right module} (of systems and their maps) over a symmetric monoidal double category (of interfaces, interactions, and their maps). That is, we see loose morphisms in symmetric monoidal double categories as \emph{interaction processes} by which systems can be composed. Seeing operads as process theories which act on systems by composing them is the approach taken in the recent \cite{selby-2025-generalised}; we add to this approach the \emph{maps} of systems and their interfaces.
		\item{}			We  give a number of (pseudo-)functorial constructions of systems theories, incorporating existing work and providing new examples, including:
			\begin{enumerate}\item{}					For any  (symmetric monoidal double category), we associate its systems theory of \emph{initial processes} \(S \colon  \ast  \mathrel {\mkern 3mu\vcenter {\hbox {$\shortmid $}}\mkern -10mu{\to }} I\). See \hyperref[ssl-003Z]{Section \ref{ssl-003Z}} for this construction.
				\item{}					For any category \(C\) with pushouts (such as a category of Petri nets \cite{baez-2020-open}, stock and flow diagrams \cite{baez-2023-compositional}, or other \emph{presentations} of systems), together with some marked objects \(P \subseteq  C\) called \emph{ports}, we associate the systems theory of objects of \(C\) composing by gluing together their ports as specified by \emph{undirected wiring diagrams} \cite{fong-2018-hypergraph}. See \hyperref[ssl-0041]{Section \ref{ssl-0041}} for this construction, with undirected wiring diagrams covered in \hyperref[ssl-005D]{Example \ref{ssl-005D}}.
				\item{}					For any category \(C\) with pullbacks, we will associate a \emph{behavioral} systems theory, where systems are "variable sets" of system variables, and composition is by \emph{sharing} exposed variables. This categorifies Willem's \emph{behavior approach} to systems theory \cite{willems-2007-behavioral} (as put forward by Schultz, Spivak, and Vasilakopolous \cite{schultz-2019-dynamical}); it also includes compositional Hamiltonian and Lagrangian mechanics (as in \cite{baez-2021-open} and in Section 1.3.2 of \cite{schreiber-20xx-differential}). See \hyperref[ssl-0041]{Section \ref{ssl-0041}} for this construction.
				\item{}					For any tangent category \cite{cockett-2019-differential}, we will associate a systems theory of ODEs composing by \emph{directed wiring diagrams} \cite{vagner-2014-algebras}. Indeed, we extend this to any "first order differential structure (FODS)" as in \cite{capucci-2024-fibrational}. See \hyperref[djm-00HF]{Section \ref{djm-00HF}}.
				\item{}					For any lax monoidal endofunctor \(F\) on a cartesian category, we will associate a systems theory of \emph{open \(F\)-coalgebras} or \(F\)-Moore machines composing by lenses, including directed wiring diagrams. This includes partially observable Markov decision processes (for \(F\) a probability monad), as well as ordinary (possibly non-deterministic) Moore machines. See \hyperref[djm-00EE]{Section \ref{djm-00EE}}.\end{enumerate}
			These constructions are \emph{pseudo-functorial} in their underlying data; we therefore see that systems theories organize into \hyperref[djm-00A8]{\textbf{doctrines}} according to the data needed to specify them. This pseudo-functoriality gives us a number of \emph{recipes} for ``black-boxing'' functors --- in this case, pseudo-symmetric monoidal pseudo-functors between symmetric monoidal loose right modules. This extends the approach of Fong and Sarazola \cite{fong-2018-recipe} in constructing morphisms of hypergraph categories from morphisms between decoration data for corelations.  Our theorems above give recipes not only for giving morphisms between hypergraph (double) categories (symmetric monoidal modules over cospan double categories), but also between modules for lens double categories of, e.g. ODEs or POMDPs.
		\item{}			We will observe that many of the diagrammatic languages used in categorical systems theories arise as \emph{free processes} in the various doctrines of systems theories. In particular, undirected wiring diagrams are free processes in both span- and cospan-based systems theories, while directed wiring diagrams are free processes in lens-based systems theories. See \hyperref[djm-00G5]{Section \ref{djm-00G5}}. This let's us recover hypergraph categories of decorated or structured cospans and of relational systems which compose by sharing variables, as well as the algebras for wiring diagrams of directed operads considered by Spivak et. al. \cite{schultz-2019-dynamical}. All constructions are pseudo-functorial, giving \emph{recipes} for ``black boxing functors'' between a variety of different sorts of systems.\end{enumerate}\subsection{Future work}\label{ssl-002X}\par{}This paper is the beginning of a series on the double operadic theory of systems, based on the \emph{formal category theory} of symmetric monoidal loose right modules on symmetric monoidal double categories. In future work we will:
\begin{enumerate}\item{}		Examine the \emph{Yoneda theory} of systems theories (extending the work on representable morphisms of systems theories in Chapter 5 of \cite{jaz-2021-book}). We will show that the simplest form of Willems' behavioral approach to systems theory, as categorified in the sheaf approach of Schultz, Spivak, and Vasilakopolous \cite{schultz-2019-dynamical}, is a \emph{discrete opfibration classifier} in the 2-category of systems theories. In particular, features of systems theories \emph{representable} by maps (such as \emph{trajectories}, \emph{steady states}, etc. but also \emph{control-barrier functions} and \emph{cocycles}) give (sometimes lax) morphisms of systems theories into Willems' style behavioral systems theories (as demonstrated in the manuscript \cite{jaz-2021-book}). This gives a robust class of compositionality theorems. We will explore how time variation in system behaviors arises out of a choice of a category of clock-systems which represent time-varying behavior, and connect this with the sheaf theoretic approach of \cite{schultz-2019-dynamical}.
	\item{}		Explore \emph{assume-guarantee} reasoning and compositional system validation using \emph{slice} systems theories. Since maps between systems not only represent behaviors of those systems but also can express the satisfaction of certain properties (for example, Lyapunov functions \cite{ames-2025-categorical} which witness stability), we can form \emph{certified systems theories} as \textbf{slices} (comma objects) in the 2-category of systems theories. Compositionality within these certified systems theories gives a form of \emph{assume-guarantee} reasoning, or compositional model checking.\end{enumerate}
\begin{acknowledgements}[{}]\label{djm-00JI}
\par{}	We would like to thank Mitchell Riley for his careful reading and comments during the drafting of this paper. This project was funded by the Advanced Research + Invention Agency
(ARIA).\end{acknowledgements}
\section{Preliminaries}\label{ssl-003R}\par{}	In this section, we will collect a number of preliminary notions which we will use throughout the remainder of the paper. In particular, we will review the some necessary 2-category theory; recall the definition of the 2-category of double categories which we will use in this paper; and recall the theory of \hyperref[djm-00AF]{\emph{adequate triples}} which can be used to give tight control over the \(\mathbb {S}\mathsf {pan}\)-construction. Finally, we will review Spivak's notion of \emph{lens} (\cite{spivak-2019-generalized}) via the \(\mathbb {S}\mathsf {pan}\)-construction for adequate triples.\subsection{2-Categories and \(\mathcal {F}\)-sketches}\label{djm-00J7}\par{}	In this section, we review some necessary background from the theory of 2-categories and 2-categorical algebra.\par{}	By a 2-category we will mean a category enriched in the cartesian monoidal (1-)category of categories. That is, by 2-category we will always mean something strict. For a comprehensive account, we refer the reader to \cite{johnson-2020-dimensional}; for a more informal introduction, see \cite{lack-2009-categories}. We note that as a category of categories enriched over a cartesian monoidal category with all limits, the category of 2-categories admits all pullbacks.\par{}	We will make extensive use of constructions involving 2-categories with (strict, 2-categorical) products. For this reason, we recall the following definition of a \emph{cartesian} 2-category.
\begin{definition}[{Cartesian 2-category}]\label{djm-00AN}
\par{}A 2-category \(\mathcal {K}\) is \textbf{cartesian} when it has all finite 2-categorical products --- products as categories enriched in \(\mathsf {Cat}\). That is, for each pair of objects \(A\) and \(B\) in \(\mathcal {K}\), there is an object \(A \times  B\) together with 2-natural isomorphisms of categories:
\[\mathcal {K}(Z, A \times  B) \cong  \mathcal {K}(Z, A) \times  \mathcal {K}(Z, B),\]
and there is a terminal object \(1\) for which
\[\mathcal {K}(Z, 1) \cong  1.\]\par{}A 2-functor which preserves finite products is called a \textbf{cartesian} 2-functor.\end{definition}
\par{}	It is straightforward to show that cartesian 2-categories, and cartesian 2-functors form a category, and that this 2-category admits all pullbacks (constructed in \(2\mathcal {C}\mathsf {at}\)).\par{}	The reason we are interested in cartesian 2-categories is they are the appropriate setting for describing 2-algebraic structure --- such as symmetric monoidal or cartesian structure --- on objects. To fix an account of 2-algebra in this paper, we will use the \(\mathcal {F}\)-sketch formalism of the excellent \cite{arkor-2024-enhanced}. We will review this briefly here.\par{}	An \(\mathcal {F}\)-category (originally defined in \cite{lack-2012-enhanced}) is a 2-category with a class of its 1-cells marked \emph{inert}. An \(\mathcal {F}\)-functor is a 2-functor which sends inert 1-cells to inert 1-cells.
\begin{remark}[{Terminology for \(\mathcal {F}\)-categories}]\label{djm-00J8}
\par{}	In both \cite{lack-2012-enhanced} and \cite{arkor-2024-enhanced}, the marked 1-cells of an \(\mathcal {F}\)-category are called \emph{tight} and the general 1-cells are called \emph{loose}. In this paper, we will use the terms ``tight'' and ``loose'' to refer to the two sorts of morphisms in a double category; for this reason, we will name the marked 1-cells of an \(\mathcal {F}\)-sketch \emph{inert} in this paper. As for why we chose the word ``inert'', see below.\par{}	In Remark 5.11 of \cite{arkor-2024-enhanced}, the authors describe how in many cases the inert morphisms of an \(\mathcal {F}\)-category are chosen as display maps, and they give an example of the theory of a category with hom-sets displayed over sets of objects. The simplex category \(\Delta \) may be equipped with the structure of a limit sketch whose models are categories (via the \emph{Segal conditions}); the maps which display the hom-sets (\(X(\Delta [1])\)) over the object-sets (\(X(\Delta [0])\)) for such a model \(X\) are ``inert'' maps in the active-inert factorization system on \(\Delta \). This is why we choose the name ``inert''.\end{remark}
\par{}	We will generally consider our 2-categories as \emph{chordate} \(\mathcal {F}\)-categories, meaning that all of their 1-cells are marked as inert. Any cartesian 2-category, when considered as a chordate \(\mathcal {F}\)-category, supports the definition of (symmetric) monoidal objects --- pseudo-monoids --- as described in Definition 4.13 of \cite{arkor-2024-enhanced}. To describe the 2-category \(\mathcal {S}\mathsf {M}(\mathcal {K})\) of symmetric pseudo-monoids in a cartesian 2-category \(\mathcal {K}\), we will appeal to the theory of \(\mathcal {F}\)-sketches.\par{}	An \(\mathcal {F}\)-sketch (Definition 5.8 of \cite{arkor-2024-enhanced}) is an \(\mathcal {F}\)-category equipped with some 2-categorical cones marked as formal limit cones. A model of an \(\mathcal {F}\)-sketch \(\mathcal {A}\) in another \(\mathcal {K}\) is an \(\mathcal {F}\)-functor \(M : \mathcal {A} \to  \mathcal {K}\) which sends these marked cones to marked cones; usually, we take \(\mathcal {K}\) to be a suitably complete 2-category and mark the limit cones, so that a model becomes an \(\mathcal {F}\)-functor \(M : \mathcal {A} \to  \mathcal {K}\) which sends the marked formal limit cones to actual limit cones.\par{}	See Definition 5.10 of \cite{arkor-2024-enhanced} for a definition of the \(\mathcal {F}\)-category \(\mathcal {M}\mathsf {od}_w(\mathcal {A}; \mathcal {K})\) for the \(\mathcal {F}\)-category of models of an \(\mathcal {F}\)-sketch \(\mathcal {A}\) with \(w\)-morphisms for \(w\) a \emph{weakness}: strict, pseudo, lax, or colax. With this definition in mind, we may make the following definitions which we will use in this paper.
\begin{definition}[{\(\mathcal {S}\mathsf {M}\) and \(\mathcal {C}\mathsf {art}\)}]\label{djm-00J9}
\par{}	Let \(\mathcal {K}\) be a \hyperref[djm-00AN]{\emph{cartesian 2-category}}. We then define
	\[
		\mathcal {S}\mathsf {M}(\mathcal {K}) := \mathcal {M}\mathsf {od}_p(\mathcal {A}_{\mathcal {S}\mathsf {M}}; \mathcal {K})
	\]
	and 
	\[
		\mathcal {C}\mathsf {art}(\mathcal {K}) := \mathcal {M}\mathsf {od}_p(\mathcal {A}_{\mathcal {C}\mathsf {art}}; \mathcal {K})
	\]
	to be the 2-categories of models and pseudo-morphisms where \(\mathcal {A}_{\mathcal {S}\mathsf {M}}\) is the \(\mathcal {F}\)-sketch of symmetric pseudo-monoids (extending Example 5.12 of \cite{arkor-2024-enhanced} by a symmtery isomorphism in an evident way), and \(\mathcal {A}_{\mathcal {C}\mathsf {art}}\) is the \(\mathcal {F}\)-sketch of cartesian objects (pseudo-monoids where multiplication is adjoint to diagonal and unit is adjoint to the terminal map).\end{definition}
\par{}	Noting the covariant representablility of symmetric monoidal objects by the symmetric pseudo-monoid \(\mathcal {F}\)-sketch in this way makes it clear that not only is the assignment \(\mathcal {K} \mapsto  \mathcal {S}\mathsf {M}(\mathcal {K})\) functorial with respect to cartesian 2-functors, but also that this assignment preserves pullbacks of 2-categories:
	
\begin{center}
		\begin {tikzcd}
	{\mathcal {K}_1} & {\mathcal {K}_3} && {\mathcal {S}\mathsf {M}(\mathcal {K}_1)} & {\mathcal {S}\mathsf {M}(\mathcal {K}_1)} \\
	{\mathcal {K}_2} & {\mathcal {K}_4} && {\mathcal {S}\mathsf {M}(\mathcal {K}_1)} & {\mathcal {S}\mathsf {M}(\mathcal {K}_1)}
	\arrow [from=1-1, to=1-2]
	\arrow [from=1-1, to=2-1]
	\arrow [""{name=0, anchor=center, inner sep=0}, from=1-2, to=2-2]
	\arrow [from=1-4, to=1-5]
	\arrow [""{name=1, anchor=center, inner sep=0}, from=1-4, to=2-4]
	\arrow [from=1-5, to=2-5]
	\arrow [""{name=2, anchor=center, inner sep=0}, from=2-1, to=2-2]
	\arrow [""{name=2p, anchor=center, inner sep=0}, phantom, from=2-1, to=2-2, start anchor=center, end anchor=center]
	\arrow [""{name=3, anchor=center, inner sep=0}, from=2-4, to=2-5]
	\arrow ["\lrcorner "{anchor=center, pos=0.125}, draw=none, from=1-1, to=2p]
	\arrow [shorten <=20pt, shorten >=20pt, maps to, from=0, to=1]
	\arrow ["\lrcorner "{anchor=center, pos=0.125}, draw=none, from=1-4, to=3]
\end {tikzcd}
	\end{center}\par{}	Finally, we recall the wonderful \emph{symmetry of internalization} result from Theorem 7.5 of \cite{arkor-2024-enhanced}, which states that if \(\mathcal {A}\) and \(\mathcal {B}\) are \(\mathcal {F}\)-sketches which whose marked limits are entirely \emph{inert}, then 
	\[
		\mathcal {M}\mathsf {od}_w(\mathcal {A}; \mathcal {M}\mathsf {od}_{\overline {w}}(\mathcal {B}; \mathcal {K})) \cong 
		\mathcal {M}\mathsf {od}_{\overline {w}}(\mathcal {B}; \mathcal {M}\mathsf {od}_{w}(\mathcal {A}; \mathcal {K}))
	\]
	for any weakness \(w\), where \(\overline {w}\) is the dual weakness (so, \(\overline {\mathsf {lax}} = \mathsf {colax}\) and vice versa, but \(\overline {\mathsf {pseudo}} = \mathsf {pseudo}\)).\par{}	We will use the the theory of \(\mathcal {F}\)-sketches fluently to compute what it means to be a symmetric monoidal object of a variety of 2-categories under consideration in this paper.\subsection{Double categorical preliminaries}\label{djm-00JA}\par{}	In this section, we review the 2-category \(\mathcal {D}\mathsf {bl}\) of double categories, (pseudo-)double functors, and tight transformations. For a comprehensive review of double category theory, we review the reader to \cite{grandis-2019-higher}; all the necessary concepts we will use are also written explicitly, in conventions closer to our own, in the excellent \cite{lambert-2024-cartesian}.
\begin{definition}[{The 2-category of double categories}]\label{djm-00IR}
\par{}We take \(\mathcal {D}\mathsf {bl}\) to be the 2-category of (pseudo-)double categories, pseudo-double functors, and tight transformations. Specifically, we define
\[\mathcal {D}\mathsf {bl} := \mathcal {M}\mathsf {od}_p(\mathcal {A}_{\mathsf {PseudoCat}}; \mathcal {C}\mathsf {at})\]
to be the 2-category of models of the \hyperref[djm-00J7]{\(\mathcal {F}\)-sketch} of pseudo-categories (Example 5.13 of \cite{arkor-2024-enhanced}) valued in the 2-category of categories.\par{}See Defintion 5.2.1 of \cite{leinster-2004-higher} for a full definition of a double category (there called a ``weak'' double category) with unbiased, \(n\)-ary composition, and see Definition 3.5.1 of \cite{grandis-2019-higher} for a definition of lax double functor (pseudo-double functors are lax double functors whose unitor and laxitors are isomorphisms) and Definition 3.5.4 of \cite{grandis-2019-higher} for a definition of tight transformation (there called a ``horizontal transformation'').\end{definition}

\begin{convention}[{Double categorical terminology}]\label{djm-00JC}
\par{}	A (pseudo-)double category \(\mathbb {D}\) has an underlying span
	
\begin{center}
	\begin {tikzcd}
	& {\mathsf {Loose}({\mathbb {D}})} \\
	{ \mathsf {Tight}({\mathbb {D}})} && { \mathsf {Tight}({\mathbb {D}})}
	\arrow ["{\mathsf {dom}}"', from=1-2, to=2-1]
	\arrow ["{\mathsf {codom}}", from=1-2, to=2-3]
\end {tikzcd}
	\end{center}

	where \(\mathsf {Tight}({\mathbb {D}})\) is its category of \emph{objects} and \emph{tight} morphisms, and \(\mathsf {Loose}({\mathbb {D}})\) is its category of \emph{loose} morphisms and \emph{squares}. That is, we refer to those morphisms of \(\mathbb {D}\) which compose strictly as ``tight'' morphisms, and those that compose up to coherent isomorphism as ``loose'' morphisms. We refer to composition of morphisms in \(\mathsf {Tight}({\mathbb {D}})\) and \(\mathsf {Loose}({\mathbb {D}})\) as \emph{tight composition}, and composition of \(\mathbb {D}\) qua pseudo-category in \(\mathcal {C}\mathsf {at}\) as \emph{loose composition}.\par{}	We will tend to draw our tight morphisms going down the page and loose morphisms going across the page.We avoid using the terminology ‘vertical’ and ‘horizontal’ to avoid confusion with other conventions in the double categorical literature.\end{convention}
\par{}	We recall that the 2-category of double categories is cartesian, and that products of double categories may be constructed in \(\mathcal {C}\mathsf {at}\).
\begin{lemma}[{\(\mathcal {D}\mathsf {bl}\) is cartesian}]\label{djm-00JB}
\par{}	The \hyperref[djm-00IR]{2-category \(\mathcal {D}\mathsf {bl}\) of double categories} is \hyperref[djm-00AN]{cartesian}.\end{lemma}
\par{}	As such, we may consider \emph{symmetric monoidal double categories}, which are \hyperref[djm-00J9]{symmetric monoidal objects} in the 2-category \(\mathcal {D}\mathsf {bl}\). Similarly, \emph{cartesian double categories} (as considered, for example, in \cite{lambert-2024-cartesian}) are \hyperref[djm-00J9]{cartesian objects of \(\mathcal {D}\mathsf {bl}\)}.\par{}	Finally, we recall two technical definitions which will appear in the upcoming \hyperref[djm-00IQ]{Section \ref{djm-00IQ}} and be important constraints in the rest of the paper.
\begin{lemma}[{Tight isomorphisms are commuter cells}]\label{djm-00EY}
\par{}	Let \(f : x \Rightarrow  y\) be a tight isomorphism between loose maps:
	
\begin{center}
		\begin {tikzcd}
	{x_0} & {x_1} \\
	{y_0} & {y_1}
	\arrow [""{name=0, anchor=center, inner sep=0}, "x", "\shortmid "{marking}, from=1-1, to=1-2]
	\arrow ["{f_0}"', from=1-1, to=2-1]
	\arrow ["{f_1}", from=1-2, to=2-2]
	\arrow [""{name=1, anchor=center, inner sep=0}, "y"', "\shortmid "{marking}, from=2-1, to=2-2]
	\arrow ["f"{description}, draw=none, from=0, to=1]
\end {tikzcd}
	\end{center}

	Suppose that \(f_0\) and \(f_1\) are companions; then \(f\) is a companion \hyperref[djm-00EX]{commuter cell}. Dually, if \(f_0\) and \(f_1\) are conjoints, then \(f\) is a conjoint commuter cell.\begin{proof}\par{}	The inverse of the transposed cell \(f^{<}\) is the transposed cell associated to the tight inverse of \(f\).\end{proof}\end{lemma}

\begin{definition}[{Commuter cells}]\label{djm-00EX}
\par{}	A cell 
	
\begin{center}
		\begin {tikzcd}
	{x_0} & {x_1} \\
	{y_0} & {y_1}
	\arrow [""{name=0, anchor=center, inner sep=0}, "x", "\shortmid "{marking}, from=1-1, to=1-2]
	\arrow ["{f_0}"', from=1-1, to=2-1]
	\arrow ["{f_1}", from=1-2, to=2-2]
	\arrow [""{name=1, anchor=center, inner sep=0}, "y"', "\shortmid "{marking}, from=2-1, to=2-2]
	\arrow ["f"{description}, draw=none, from=0, to=1]
\end {tikzcd}
	\end{center}

	in a double category is said to be a \textbf{companion commuter} (Definition 8.1 of \cite{paré-2023-retrocells}) if \(f_0\) and \(f_1\) have companions \(f_0^{>}\) and \(f_1^{>}\) respectively and the associated globular transpose cell \(f^{>}\) of \(f\) is a tight isomorphism:
	
\begin{center}
		\begin {tikzcd}
	{x_0} & {x_0} & {x_1} & {y_1} \\
	{x_0} & {y_0} & {y_1} & {y_1}
	\arrow ["\shortmid "{marking}, equals, from=1-1, to=1-2]
	\arrow [equals, from=1-1, to=2-1]
	\arrow [""{name=0, anchor=center, inner sep=0}, "x", "\shortmid "{marking}, from=1-2, to=1-3]
	\arrow ["{f_0}"', from=1-2, to=2-2]
	\arrow ["{f_1^{>}}", "\shortmid "{marking}, from=1-3, to=1-4]
	\arrow ["{f_1}", from=1-3, to=2-3]
	\arrow [equals, from=1-4, to=2-4]
	\arrow ["{f_0^{>}}"', "\shortmid "{marking}, from=2-1, to=2-2]
	\arrow [""{name=1, anchor=center, inner sep=0}, "y"', "\shortmid "{marking}, from=2-2, to=2-3]
	\arrow ["\shortmid "{marking}, equals, from=2-3, to=2-4]
	\arrow ["f"{description}, draw=none, from=0, to=1]
\end {tikzcd}
	\end{center}

	Dually, \(f\) is a \textbf{conjoint commuter} if \(f_0\) and \(f_1\) have conjoints, and the associated globular transpose of \(f\) is a tight isomorphism.\end{definition}
\subsection{Adequate triples, spans, and lenses}\label{djm-00JD}\par{}	In this section, we will review the notion of \hyperref[djm-0087]{\emph{adequate triple}} (due to \cite{haugseng-2020-twovariable}) and the resulting \hyperref[djm-0088]{\emph{span} construction} and its special case, the \hyperref[djm-00B6]{\emph{lens construction}} of Spivak \cite{spivak-2019-generalized}.\subsubsection{Adequate triples and spans}\label{djm-00AF}\par{}In this section, we'll give a general span construction making use of a (1-categorical version of) \cite{haugseng-2020-twovariable}\hyperref[djm-0087]{\emph{adequate triples}}. An adequate triple is just the data needed on a category in order to perform the span construction on it.
\begin{definition}[{Adequate triple}]\label{djm-0087}
\par{}An \textbf{adequate triple} \((\mathsf {C}, (L, R))\) (Definition 1.2 of \cite{haugseng-2020-twovariable}) consists of a category \(\mathsf {C}\) together with two classes of maps \(L\) and \(R\) in \(\mathsf {C}\) satisfying the following axioms:
\begin{enumerate}\item{}Both \(L\) and \(R\) contain all identities and are closed under composition.
    \item{}In any pullback square as below with \(\ell  \in  L\) and \(r \in  R\),

\begin{center}
\begin {tikzcd}
	& p \\
	x && y \\
	& z
	\arrow ["{\ell ' \in  L}"', dashed, from=1-2, to=2-1]
	\arrow ["{r' \in  R}", dashed, from=1-2, to=2-3]
	\arrow ["\lrcorner "{anchor=center, pos=0.125, rotate=-45}, draw=none, from=1-2, to=3-2]
	\arrow ["{r \in  R}"', from=2-1, to=3-2]
	\arrow ["{\ell  \in  L}", from=2-3, to=3-2]
\end {tikzcd}
\end{center}

the arrow \(r'\) is in \(R\) and \(\ell '\) is in \(L\). Furthermore, for any cospan as drawn in solid above, there is such a pullback. We will call these pullback squares "\(L\)-\(R\) pullbacks".\end{enumerate}
The 2-category \(\mathcal {A}\mathsf {dTr}\) consists of adequate triples, functors preserving the classes \(L\) and \(R\) and preserving \(L\)-\(R\) pullbacks, and arbitrary natural transformations.\end{definition}
\par{}The main reason for defining the notion of \hyperref[djm-0087]{adequate triples} is that they capture (almost) the minimally necessary data and properties of a category needed to perform the span construction.
\begin{construction}[{Span construction of an adequate triple}]\label{djm-0088}
\par{}	We construct a 2-functor \(\mathbb {S}\mathsf {pan} : \mathcal {A}\mathsf {dTr} \to  \mathcal {D}\mathsf {bl}_u\) sending an adequate triple \((\mathsf {C}, (L, R))\) to the double category \(\mathbb {S}\mathsf {pan}(\mathsf {C}, (L, R))\) whose tight category is \(\mathsf {C}\) and whose loose morphisms are spans
	
\begin{center}
		\begin {tikzcd}
			& y \\
			x && z
			\arrow ["{\ell  \in  L}"', from=1-2, to=2-1]
			\arrow ["{r \in  R}", from=1-2, to=2-3]
		\end {tikzcd}
	\end{center}

	whose left leg is in \(L\) and whose right leg is in \(R\), with composition given by pullback (noting that an adequate triple asks for precisely the sorts of pullbacks we need to compose these spans). Squares are the usual maps of spans.\par{}A map \(F : (\mathsf {C}_1, (L_1, R_1)) \to  (\mathsf {C}_2, (L_2, R_2))\) of adequate triples gets sent to the double functor given by applying \(F\) to all elements.\begin{proof}\par{}That the span construction gives a 2-functor into double categories appears as Proposition 3.26 of \cite{dawson-2010-span}. While this performs the usual construction beginning with a category with all pullbacks, we only actually ever take the \(L\)-\(R\) pullbacks guarenteed by the axioms of an \hyperref[djm-0087]{adequate triple}, and double functoriality only requires preserving these \(L\)-\(R\) pullbacks. Since \(F\) preserves identities, and since identity spans are those with both legs identities, \(\mathbb {S}\mathsf {pan}(F)\) is a strictly unitary double functor.\end{proof}\end{construction}

\begin{lemma}[{Companions and conjoints in \(\mathbb {S}\mathsf {pan}\) of an adequate triple}]\label{djm-00F4}
\par{}	Let \((\mathsf {C}, (L, R))\) be an \hyperref[djm-0087]{adequate triple}. Then \(f : c_0 \to  c_1\) is a companion (resp. conjoint) in \hyperref[djm-0088]{the span construction \(\mathbb {S}\mathsf {pan}(\mathsf {C}, (L, R))\)} if and only if it is in \(R\) (resp. \(L\)).\begin{proof}\par{}	In a span double category, \(f : c_0 \to  c_1\) has a companion given by 
	
\begin{center}
		\begin {tikzcd}
	& {c_0} \\
	{c_0} && {c_0}
	\arrow [equals, from=1-2, to=2-1]
	\arrow ["f", from=1-2, to=2-3]
\end {tikzcd}
	\end{center}

	This only works in \(\mathbb {S}\mathsf {pan}(\mathsf {C}, (L, R))\) if \(f \in  R\), by construction. Dually for conjoints.\end{proof}\end{lemma}
\par{}The \hyperref[djm-0088]{span construction} gives us a 2-functorial way to produce double categories. We can easily extend the \hyperref[djm-0088]{span construction} to produce \emph{symmetric monoidal} double categories by observing that it preserves cartesian products, and therefore preserves symmetric monoidal objects.
\begin{lemma}[{Adequate triples form a cartesian 2-category}]\label{djm-00AK}
\par{}The 2-category \(\mathcal {A}\mathsf {dTr}\) of \hyperref[djm-0087]{adequate triples} has cartesian products given by 
\[(\mathsf {C}_1, (L_1, R_1)) \times  (\mathsf {C}_2, (L_2, R_2)) = (\mathsf {C}_1 \times  \mathsf {C}_2, (L_1 \times  L_2, R_1 \times  R_2)).\]\begin{proof}\par{}This is straightforward to verify. The projections preserve (and jointly reflect) both classes by definition, and a square in \(\mathsf {C}_1 \times  \mathsf {C}_2\) is a pullback if and only if both its components are pullbacks.\end{proof}\end{lemma}
\par{}Next, we explicate symmetric monoidal adequate triples and show that they are preserved by the span construction.
\begin{lemma}[{Symmetric monoidal adequate triples}]\label{djm-00AQ}
\par{}A \textbf{symmetric monoidal adequate triple} (an object of \(\mathcal {S}\mathsf {M}{(\mathcal {A}\mathsf {dTr})}\)) is equivalently a symmetric monoidal category \((\mathsf {C}, \otimes , 1)\) equipped with the structure of an adequate triple \((\mathsf {C}, (L, R))\), such that:
\begin{enumerate}\item{}The monoidal product preserves both classes: if \(f\) and \(g\) are both in \(L\) or both in \(R\), then so is \(f \otimes  g\).
	\item{}		The monoidal product preserves \(L\)-\(R\) pullbacks: if the square below are pullbacks
		
\begin{center}
			\begin {tikzcd}
				{C_1} & {C_3} & {C_1'} & {C_3'} \\
				{C_2} & {C_4} & {C_2'} & {C_4'}
				\arrow ["{\in  R}", from=1-1, to=1-2]
				\arrow ["{\in  L}"', from=1-1, to=2-1]
				\arrow ["\lrcorner "{anchor=center, pos=0.125}, draw=none, from=1-1, to=2-2]
				\arrow ["{\in  L}", from=1-2, to=2-2]
				\arrow ["{\in  R}", from=1-3, to=1-4]
				\arrow ["{\in  L}"', from=1-3, to=2-3]
				\arrow ["\lrcorner "{anchor=center, pos=0.125}, draw=none, from=1-3, to=2-4]
				\arrow ["{\in  L}", from=1-4, to=2-4]
				\arrow ["{\in  R}"', from=2-1, to=2-2]
				\arrow ["{\in  R}"', from=2-3, to=2-4]
			\end {tikzcd}
		\end{center}

		then the diagram below is a pullback.
		
\begin{center}
			\begin {tikzcd}
				{C_1 \otimes  C_1'} & {C_3 \otimes  C_3'} \\
				{C_2 \otimes  C_2'} & {C_4 \otimes  C_4'}
				\arrow ["{\in  R}", from=1-1, to=1-2]
				\arrow ["{\in  L}"', from=1-1, to=2-1]
				\arrow ["\lrcorner "{anchor=center, pos=0.125}, draw=none, from=1-1, to=2-2]
				\arrow ["{\in  L}", from=1-2, to=2-2]
				\arrow ["{\in  R}"', from=2-1, to=2-2]
			\end {tikzcd}
		\end{center}\end{enumerate}\end{lemma}

\begin{lemma}[{The span construction is a cartesian 2-functor}]\label{djm-00AP}
\par{}The \hyperref[djm-0088]{span construction} \(\mathbb {S}\mathsf {pan} : \mathcal {A}\mathsf {dTr} \to  \mathcal {D}\mathsf {bl}_u\) is a \hyperref[djm-00AN]{cartesian 2-functor}.\begin{proof}\par{}This is straightforward by the construction of products \hyperref[djm-00AK]{in \(\mathcal {A}\mathsf {dTr}\)} and . By \hyperref[djm-00AK]{Lemma \ref{djm-00AK}}, the product in \(\mathcal {A}\mathsf {dTr}\) given by the product of the underlying categories, and a span in the product is a pair of spans.\end{proof}\end{lemma}
\subsubsection{Lex and rex categories}\label{ssl-004Y}\par{}An important set of examples of \hyperref[djm-0087]{adequate triples} come from lex and rex categories, which we define below.\paragraph{Lex categories}\label{djm-00JF}\par{}First, we can recover the ordinary span construction from the \hyperref[djm-0088]{span construction for adequate triples} by taking a category with finite limits (called a \hyperref[djm-001Z]{lex category}) and turning it into an \hyperref[djm-0087]{adequate triple} where both classes contain all maps.
\begin{definition}[{The 2-category of categories with finite limits}]\label{djm-001Z}
\par{}Define \(\mathcal {L}\mathsf {ex}\) to be the 2-category whose objects are locally small categories with finite limits (also called \textbf{lex categories}), whose morphisms are functors which preserve finite limits (also called \textbf{lex functors}), and whose 2-cells are natural transformations.\end{definition}

\begin{lemma}[{The 2-category \(\mathcal {L}\mathsf {ex}\) of lex categories is cartesian}]\label{djm-00AT}
\par{}The 2-category \hyperref[djm-001Z]{\(\mathcal {L}\mathsf {ex}\) of lex categories} is \hyperref[djm-00AN]{cartesian}, and its finite products may be computed on underlying categories.\begin{proof}\par{}The terminal category is lex, and the terminal morphisms into it all preserve finite limits. Limits in products are computed component-wise, so if two categories have finite limits then so does their product, and the projections preserve these limits.\end{proof}\end{lemma}

\begin{construction}[{Adequate triple of a lex category}]\label{djm-008C}
\par{}There is a faithful, cartesian 2-functor \(\mathcal {L}\mathsf {ex} \to  \mathcal {A}\mathsf {dTr}\) from \hyperref[djm-001Z]{the 2-category of lex categories} to \hyperref[djm-0087]{the 2-category of adequate triples} which sends a lex category \(\mathsf {C}\) to the triple \((\mathsf {C}, (\mathsf {all}, \mathsf {all}))\) where both classes consist of all maps. Any lex functor trivially preserves both classes. It also preserves all pullbacks and hence all \(L\)-\(R\) pullbacks.\end{construction}

\begin{observation}[{Span of adequate triple from lex category is ordinary span construction}]\label{djm-00AS}
\par{}The \hyperref[djm-0088]{span construction} \(\mathbb {S}\mathsf {pan}(\mathsf {C}, (\mathsf {all}, \mathsf {all}))\) of the \hyperref[djm-0088]{adequate triple} \hyperref[djm-008C]{constructed from a lex category \(\mathsf {C}\)} is the ordinary span construction \(\mathbb {S}\mathsf {pan}(\mathsf {C})\) as described, for example, in Section 3 of \cite{dawson-2010-span}.\end{observation}
\par{}In addition the the usual but rather austere spans of sets, we also find spans in richer categories that can act as systems in their own right \cite{baez-2021-open}. For example, Lagrangian correspondences.
\begin{example}[{
	Lagrangian spans
}]\label{djm-00FU}
\par{}	In 1.3.2 of \cite{schreiber-20xx-differential}, Schreiber describes a simple construction of a double category faithfully including symplectic manifolds and Lagrangian correspondences. We summarize this construction here, as an example of the variable sharing doctrine.\par{}	Let \(H\) be a topos of \emph{smooth sets} such as sheaves on the site of smooth manifolds and open covers (see, \cite{giotopoulos-2023-field} for a comprehensive introduction to this topos), or Dubuc's topos \cite{dubuc-1979-sur}. For any of these sites, there is a notion of differential form definable on the objects of the site --- for smooth manifold \(M\), the set \(\Lambda ^n(M)\) of alternating \(n\)-forms on \(M\) --- stable under pullback and defined locally; that is, there is a sheaf \(\Lambda ^n\) which \emph{represents} differential \(n\)-forms in the sense that a form \(\omega  \in  \Lambda ^n(M)\) is equivalently given by a map \(\omega  : M \to  \Lambda ^n\). Since forms are naturally an abelian group, \(\Lambda ^n\) is an abelian group object of \(H\).\par{}	Since the differential commutes with pullback of forms, there will be a map \(d : \Lambda ^2 \to  \Lambda ^3\); its kernel \(\Lambda ^2_{\mathsf {cl}}\) classifies closed 2-forms. A map \(\omega  : M \to  \Lambda ^2_{\mathsf {cl}}\) is therefore a closed alternating 2-form on \(M\), also known as a \emph{pre-symplectic structure}. Since precomposition by \(f : N \to  M\) sends \(\omega  : M \to  \Lambda ^2_{\mathsf {cl}}\) to its pullback \(f^{\ast }\omega  : N \to  \Lambda ^2_{\mathsf {cl}}\), we see that the slice topos \(H \downarrow  \Lambda ^{2}_{\mathsf {cl}}\) consists of the pre-symplectic manifolds and the \emph{symplectomorphisms} between them: smooth maps \(f : N \to  M\) such that \(f^{\ast } \omega _M = \omega _N\).\par{}	As a topos, \(H \downarrow  \Lambda ^2_{\mathsf {cl}}\) has finite limits, and therefore lives in \(\mathcal {L}\mathsf {ex}\). A span in \(H \downarrow  \Lambda ^2_{\mathsf {cl}}\) is therefore a commuting square
	
\begin{center}
		\begin {tikzcd}
	& S \\
	M && N \\
	& {\Lambda ^2_{\mathsf {cl}}}
	\arrow ["m"', from=1-2, to=2-1]
	\arrow ["n", from=1-2, to=2-3]
	\arrow ["{\omega _M}"', from=2-1, to=3-2]
	\arrow ["{\omega _N}", from=2-3, to=3-2]
\end {tikzcd}
	\end{center}

	so that \(m^{\ast }\omega _M = n^{\ast }\omega _N\), or \(m^{\ast }\omega _M - n^{\ast }\omega _N = 0\). If of the appropriate dimension (\(\mathsf {dim}(S) = \frac {1}{2}(\mathsf {dim}(M) + \mathsf {dim}(N))\)), then such a span is a \emph{Lagrangian correspondence} (see, for example \cite{weinstein-2009-symplectic}).\end{example}
\paragraph{Rex categories}\label{djm-00JG}\par{}As a straightforward dual of the span construction for a category with finite limits, we can find the cospan construction of a category with finite colimits as a special case of the \hyperref[djm-0088]{span construction for adequate triples}.
\begin{definition}[{2-Category of rex categories}]\label{djm-00AY}
\par{}We say that a category is \textbf{rex} when it has all finite colimits. The 2-category \(\mathcal {R}\mathsf {ex}\) has the rex categories as objects, finite colimit preserving functors as maps, and (general) natural transformations as 2-cells.\end{definition}
\par{}Since a finite colimit in \(\mathsf {C}\) is equivalently a finite limit in \(\mathsf {C}^{\mathsf {op}}\), and since \((-)^{\mathsf {op}} : \mathcal {C}\mathsf {at} \cong  \mathcal {C}\mathsf {at}\) is an involution of the 2-category \(\mathcal {C}\mathsf {at}\) of categories, we have the following duality.\par{}Note that in the following, for a 2-category \(\mathcal {K}\), we define \({\mathcal {K}}^{\mathsf {co}}\) to be the \textbf{2-dimensional dual} of \(\mathcal {K}\): its objects and 1-morphisms are the objects and 1-morphisms of \(\mathcal {K}\) and its 2-cells are dualized.
\begin{observation}[{Duality between lex and rex categories}]\label{djm-00B1}
\par{}The 2-functor \((-)^{\mathsf {op}} : {\mathcal {R}\mathsf {ex}}^{\mathsf {co}} \to  \mathcal {L}\mathsf {ex}\) sending a rex category \(C\) to its dual \(C^{\mathsf {op}}\) is an isomorphism of 2-categories with inverse also given by taking the opposite.\end{observation}
\par{}This is what lets us explicitly transport results for lex categories. For clarity, we record them in full nonetheless.
\begin{lemma}[{The 2-category \(\mathcal {R}\mathsf {ex}\) of rex categories is cartesian}]\label{djm-00AZ}
\par{}The \hyperref[djm-00AY]{2-category \(\mathcal {R}\mathsf {ex}\) of rex categories} is \hyperref[djm-00AN]{cartesian}, and products are constructed in \(\mathcal {C}\mathsf {at}\).\end{lemma}

\begin{construction}[{Adequate triple of a rex category}]\label{djm-008D}
\par{}There is a faithful, cartesian 2-functor \({\mathcal {R}\mathsf {ex}}^{\mathsf {co}} \to  \mathcal {A}\mathsf {dTr}\) from the 2-dimensional dual of the 2-category of finite cocomplete categories and finitely cocontinuous functors to \hyperref[djm-0087]{the 2-category of adequate triples} given by sending \(\mathsf {C}\) to the triple \((\mathsf {C}^{\mathsf {op}}, (\mathsf {all}, \mathsf {all}))\) where both classes consist of all maps.\end{construction}

\begin{observation}[{Cospan construction is a span construction}]\label{djm-00B0}
\par{}The \hyperref[djm-0088]{span construction} \(\mathbb {S}\mathsf {pan}(\mathsf {C}^{\mathsf {op}}, (\mathsf {all}, \mathsf {all}))\) of the \hyperref[djm-008D]{adequate triple associated to a rex category \(\mathsf {C}\)} is the tight dual \({\mathbb {C}\mathsf {ospan}(\mathsf {C})}^{\mathsf {op}}\) of the cospan construction of \(\mathsf {C}\).\end{observation}
\paragraph{Trivial adequate triples}\label{djm-00JH}\par{}It is worth noting that every category \(\mathsf {C}\) may be considered as a trivial adequate triple \((\mathsf {C}, (\mathsf {id}, \mathsf {id}))\).
\begin{observation}[{Span construction of trivial adequate triple is tight inclusion}]\label{djm-00EM}
\par{}	The \hyperref[djm-0088]{span construction} \(\mathbb {S}\mathsf {pan}(\mathsf {C}, (\mathsf {id}, \mathsf {id}))\) of a  is isomorphic to the loosely discrete double category \(\mathbb {T}\mathsf {}(\mathsf {C})\) on \(\mathsf {C}\).
	
\begin{center}
		\begin {tikzcd}
	& {\mathcal {A}\mathsf {dTr}} \\
	{\mathcal {C}\mathsf {at}} && {\mathcal {D}\mathsf {bl}}
	\arrow ["{\mathbb {S}\mathsf {pan}}", from=1-2, to=2-3]
	\arrow ["{\mathsf {triv}}", from=2-1, to=1-2]
	\arrow ["{\mathbb {T}\mathsf {}(-)}"', from=2-1, to=2-3]
\end {tikzcd}
	\end{center}\end{observation}
\subsubsection{Fibrations and functor lenses}\label{djm-00B5}\par{}In this section, we'll see how the \hyperref[djm-0088]{span construction for adequate triples} can be used to construct the double categories of \emph{lenses} (in the sense of \cite{spivak-2019-generalized}) and \emph{charts} (in the sense of Definition 3.3.0.13 of \cite{jaz-2021-book}). These double categories appear as Definition 4.1 of \cite{jaz-2021-double} (and Definition 3.5.0.6 of \cite{jaz-2021-book});here, we'll make use of Proposition 4.2 of \cite{jaz-2021-double} (proven in Theorem 3.9 of \cite{myers-2020-cartesian}.) to see them as double categories of spans.
\begin{definition}[{Cartesian fibration}]\label{djm-008F}
\par{}The 2-category \(\mathcal {F}\mathsf {ib}\) of cartesian fibrations (also known as \emph{Grothendieck fibrations}) consists of the cartesian fibrations \(\pi  : E \to  B\), the cartesian functors which are strictly commuting squares

\begin{center}
\begin {tikzcd}
	E & {E'} \\
	B & {B'}
	\arrow ["{\overline {f}}", from=1-1, to=1-2]
	\arrow ["\pi "', from=1-1, to=2-1]
	\arrow ["{\pi '}", from=1-2, to=2-2]
	\arrow ["f"', from=2-1, to=2-2]
\end {tikzcd}
\end{center}

where \(\overline {f}\) sends cartesian morphisms to cartesian morphisms. A 2-cell is a pair of maps \((\alpha  : f \Rightarrow  g, \overline {\alpha } : \overline {f} \to  \overline {g})\) for which \(\pi '\overline {\alpha } = \alpha  \pi \).\end{definition}

\begin{example}[{Display map fibrations}]\label{djm-00DZ}
\par{}A useful example of \hyperref[djm-0087]{cartesian fibrations} are \emph{display map} fibrations.
\begin{definition}[{Display map category}]\label{djm-00E0}
\par{}	A \textbf{display map category} \((\mathsf {C}, \mathsf {D})\) is a category \(\mathsf {C}\) equipped with a class of maps \(\mathsf {D} \subseteq  \mathsf {Arr}(\mathsf {C})\) which is closed under pullback along arbitrary maps in \(\mathsf {C}\). That is, for any solid diagram below
	
\begin{center}
		\begin {tikzcd}
	\bullet  & \bullet  \\
	\bullet  & \bullet 
	\arrow ["{f'}", dashed, from=1-1, to=1-2]
	\arrow ["{d'}"', dashed, from=1-1, to=2-1]
	\arrow ["\lrcorner "{anchor=center, pos=0.125}, draw=none, from=1-1, to=2-2]
	\arrow ["d", from=1-2, to=2-2]
	\arrow ["f"', from=2-1, to=2-2]
\end {tikzcd}
	\end{center}

	with \(d \in  \mathsf {D}\), the dashed pullback exists and \(d'\) is also in \(\mathsf {D}\).\par{}	The 2-category \(\mathcal {D}\mathsf {isp}\) of display map categories consists of the display map categories, functors which preserve display maps and their pullbacks, and arbitrary natural transformations.\end{definition}

\begin{lemma}[{2-Category of display map categories is cartesian}]\label{djm-00E2}
\par{}	The \hyperref[djm-00E0]{2-category \(\mathcal {D}\mathsf {isp}\)} of display map categories is \hyperref[djm-00AN]{cartesian}, with 
	\[(C_1, D_1) \times  (C_2, D_2) := (C_1 \times  C_2, D_1 \times  D_2) \]\end{lemma}
\par{}Any cartesian category can be equipped with a display map structure given by its product projections.
\begin{theorem}[{Display map category associated to a cartesian category}]\label{djm-00E3}
\par{}	Let \(C\) be a cartesian category, having finite products. Then the class \(\mathsf {proj}\) of left projections \(\pi  : X \times  Y \to  X\) in \(C\) equips \(C\) with the structure of a \hyperref[djm-00E0]{display map category}.\par{}	Moreover, this construction gives a cartesian 2-functor \((-, \mathsf {proj}) : \mathcal {C}\mathsf {art} \to  \mathcal {D}\mathsf {isp}\).\begin{proof}\par{}	In a cartesian category, any square of the following form is a pullback:
	
\begin{center}
\begin {tikzcd}
	{Z \times  Y} & {X \times  Y} \\
	Z & X
	\arrow ["{f \times  \mathrm {id}_Y}", from=1-1, to=1-2]
	\arrow ["{\pi _Z}"', from=1-1, to=2-1]
	\arrow ["\lrcorner "{anchor=center, pos=0.125}, draw=none, from=1-1, to=2-2]
	\arrow ["{\pi _X}", from=1-2, to=2-2]
	\arrow ["f"', from=2-1, to=2-2]
\end {tikzcd}		
	\end{center}

	This can be straightforwardly verified using the universal property of the product, and it shows that \((C, \mathsf {proj})\) is a display map category.\par{}	Any cartesian functor preserves projections from cartesian products, so induces a display map functor. Finally, the products in \(\mathcal {C}\mathsf {art}\) and \(\mathcal {D}\mathsf {isp}\) are constructed in \(\mathcal {C}\mathsf {at}\), and so this inclusion preserves them.\end{proof}\end{theorem}

\begin{theorem}[{Display map fibration}]\label{djm-00E1}
\par{}	For any \hyperref[djm-00E0]{display map category} \((C, D)\), the codomain projection \(\partial _1 : D \to  C\) is a cartesian fibration. Moreover, this gives a cartesian 2-functor \(\mathcal {D}\mathsf {isp} \to  \mathcal {F}\mathsf {ib}\).\begin{proof}\par{}	Cartesian lifts are given by pullback, which exist by assumption. Any display map functor gives a cartesian functor since by assumption it preserves these pullbacks.\par{}	Finally, it is straightforwardly cartesian, since cartesian products in each of these 2-categories are constructed down in \(\mathcal {C}\mathsf {at}\).\end{proof}\end{theorem}

\begin{definition}[{Simple fibrations}]\label{djm-00E4}
\par{}	The cartesian 2-functor 
	\[\mathsf {Simp} : \mathcal {C}\mathsf {art} \to  \mathcal {D}\mathsf {isp} \to  \mathcal {F}\mathsf {ib} \]
	composed of \hyperref[djm-00E3]{Theorem \ref{djm-00E3}} and \hyperref[djm-00E1]{Theorem \ref{djm-00E1}} sends a cartesian category \(\mathsf {C}\) to its \textbf{simple fibration} \(\mathsf {proj}_\mathsf {C} \xrightarrow {\partial _1} \mathsf {C}\) of (left) product projections.\end{definition}
\par{}	The simple fibrations appear in systems theories as the Grothendieck constructions of the indexed categories \(c \mapsto  \mathsf {cokl}(c \times  -) : C^{\mathsf {op}} \to  \mathcal {C}\mathsf {at}\) (Definition 2.6.2.4 of \cite{jaz-2021-book}). We prove the equivalence between these two constructions here for reference.
\begin{lemma}[{Simple fibration is Grothendieck construction of \(\mathsf {cokl}(c \times  -)\)}]\label{djm-00E5}
\par{}	The simple fibration associated to \(C\) is equivalently the Grothendieck construction of the indexed category \(c \mapsto  \mathsf {cokl}(c \times  -) : C^{\mathsf {op}} \to  \mathcal {C}\mathsf {at}\).\begin{proof}\par{}	Under the equivalence between fibrations and indexed categories given by the Grothendieck construction, a display map fibration \(\partial _1 : D \to  C\) corresponds to the indexed category \(C \downarrow _D (-) : C^{\mathsf {op}} \to  \mathcal {C}\mathsf {at}\) given by sending \(c \in  C\) to the full subcategory of the slice \(C \downarrow  c\) spanned by the display maps \(D\). It therefore suffices to show that for any \(\mathsf {cokl}(c \times  -) \simeq  C \downarrow _{\mathsf {proj}} c\), natural in \(c \in  C\).\par{}	The objects of \(\mathsf {cokl}(c \times  -)\) are those of \(C\), while the objects of \(C \downarrow _{\mathsf {proj}} c\) are the product projections \(c \times  x \to  c\), which are in bijection with the objects \(x \in  C\). Furthermore, there is a bijection between cokleisli arrows \(f : c \times  x \to  y\) and maps \((\mathrm {id}_c, f) : c \times  x \to  c \times  y\) between product projections. It is straightforward to verify that this bijection is an isomorphism.\end{proof}\end{lemma}
\end{example}

\begin{construction}[{Adequate triple of a cartesian fibration}]\label{djm-008E}
\par{}There is a cartesian 2-functor \(\mathcal {F}\mathsf {ib} \to  \mathcal {A}\mathsf {dTr}\) sending a \hyperref[djm-008F]{cartesian fibration} \(\pi  : E \to  B\) to the adequate triple \((E, \mathsf {vert}, \mathsf {cart})\) whose left class is the class of vertical maps (such that \(\pi (f)\) is an isomorphism) and whose right class is the class of cartesian maps.\begin{proof}\par{}First, let's show that \((E, \mathsf {vert}, \mathsf {cart})\) is an \hyperref[djm-0087]{adequate triple}. This means showing that in any pullback square like so:

\begin{center}
\begin {tikzcd}
	& p \\
	x && y \\
	& z
	\arrow ["{\ell ' \in  \mathsf {vert}}"', dashed, from=1-2, to=2-1]
	\arrow ["{r' \in  \mathsf {cart}}", dashed, from=1-2, to=2-3]
	\arrow ["\lrcorner "{anchor=center, pos=0.125, rotate=-45}, draw=none, from=1-2, to=3-2]
	\arrow ["{r \in  \mathsf {cart}}"', from=2-1, to=3-2]
	\arrow ["{\ell  \in  \mathsf {vert}}", from=2-3, to=3-2]
\end {tikzcd}
\end{center}

we have that \(\ell '\) is vertical and \(r'\) is cartesian; and that any solid cospan above can be completed to such a pullback. This follows from Proposition 2.4 and Lemma 2.2 of \cite{myers-2020-cartesian}.\par{}Next, we need to check 2-functoriality. We note that the assignment \((\pi  : E \to  B) \mapsto  E\) gives a 2-functor \(\mathcal {F}\mathsf {ib} \to  \mathcal {C}\mathsf {at}\) by the \hyperref[djm-008F]{definition} of \(\mathcal {F}\mathsf {ib}\); it only remains to show that this 2-functor lands in \(\mathcal {A}\mathsf {dTr}\), which is to say that for any cartesian functor 

\begin{center}
\begin {tikzcd}
	E & {E'} \\
	B & {B'}
	\arrow ["{\overline {f}}", from=1-1, to=1-2]
	\arrow ["\pi "', from=1-1, to=2-1]
	\arrow ["{\pi '}", from=1-2, to=2-2]
	\arrow ["f"', from=2-1, to=2-2]
\end {tikzcd}
\end{center}

the functor \(\overline {f} : E \to  E'\) preserves both verticals, cartesians, and \(\mathsf {vert}\)-\(\mathsf {cart}\) pullbacks. Now, \(\overline {f}\) preserves cartesian maps by assumption; it preserves verticals because the above square commutes, so that if \(a : x \to  y\) in \(E\) is vertical (\(\pi  a\) is iso), then \(\pi ' \overline {f} a = f \pi  a\) is also iso. Finally, Lemma 2.2 of \cite{myers-2020-cartesian} says that any \(\mathsf {vert}\)-\(\mathsf {cart}\) square is a pullback, since \(\overline {f}\) preserves the classes it therefore preserves \(\mathsf {vert}\)-\(\mathsf {cart}\) pullbacks.\par{}Finally, we need to show that this 2-functor is cartesian. This follows straightforwardly from the fact that the cartesian product of two fibrations is \(\pi  \times  \pi ' : E \times  E' \to  B \times  B'\).\end{proof}\end{construction}
\par{}For the purposes of this paper, we'll make the following definition of the double category of Spivak lenses.
\begin{definition}[{Double category of Spivak lenses}]\label{djm-00B6}
\par{}The \textbf{double category \(\mathbb {L}\mathsf {ens}(\pi  : E \to  B)\) (often just \(\mathbb {L}\mathsf {ens}(E)\) for short) of Spivak lenses (\cite{spivak-2019-generalized}) for the fibration \(\pi  : E \to  B\)} is defined to be the \hyperref[djm-0088]{span construction} of the \hyperref[djm-008E]{adequate triple associated to \(\pi \)}:
\[\mathbb {L}\mathsf {ens}(E) := \mathbb {S}\mathsf {pan}(E, \mathsf {vert}, \mathsf {cart}).\]\end{definition}

\begin{notation}[{Spivak lenses}]\label{ssl-004S}
\par{}Let \(\pi  \colon  E \to  B\) be a cartesian fibration. 
  \begin{itemize}\item{}For \(I \in  E\), we often use \({I \choose  \pi (I)}\) to refer to the corresponding object in \(\mathbb {L}\mathsf {ens}(E)\). In other words, writing an object \({I \choose  O}\) in \(\mathbb {L}\mathsf {ens}(E)\) implies that \(I\) lives in the fiber over \(O\).
    \item{}A \textbf{lens} is a loose map in \(\mathbb {L}\mathsf {ens}(E)\) and corresponds to the pair of maps \((f, f^\#)\) which fit into the diagram below.
      
\begin{center}
        \begin {tikzcd}
          {I_1} & {f^*O_1} & {I_2} & E \\
          {O_1} & {O_1} & {O_2} & B
          \arrow [dashed, maps to, from=1-1, to=2-1]
          \arrow ["{f^\#}"', from=1-2, to=1-1]
          \arrow ["{\mathsf {lift}(f)}", from=1-2, to=1-3]
          \arrow [dashed, maps to, from=1-2, to=2-2]
          \arrow [dashed, maps to, from=1-3, to=2-3]
          \arrow ["\pi ", from=1-4, to=2-4]
          \arrow [equals, from=2-2, to=2-1]
          \arrow ["f"', from=2-2, to=2-3]
        \end {tikzcd}
      \end{center}

      We notate this pair \[{f^\# \choose  f} \colon  {I_1 \choose  O_1} \mathrel {\mkern 3mu\vcenter {\hbox {$\shortmid $}}\mkern -10mu{\to }} {I_2 \choose  O_2}.\]\end{itemize}\end{notation}

\begin{explication}[{Explication of \hyperref[djm-00B6]{Definition \ref{djm-00B6}}}]\label{djm-00B7}
\par{}\hyperref[djm-00B6]{Definition \ref{djm-00B6}} requires a bit of explanation. In Theorem 3.9 of \cite{myers-2020-cartesian}, it is shown that \(\mathbb {L}\mathsf {ens}(\pi  : E \to  B)\) (as defined in \hyperref[djm-00B6]{Definition \ref{djm-00B6}}) is equivalent to the \emph{Grothendieck double construction} of the indexed category \(E_{(-)} : B^{\mathsf {op}} \to  \mathcal {C}\mathsf {at}\) associated to \(\pi \) (Definition 3.8 of \cite{myers-2020-cartesian}; also Definition 4.1 of \cite{jaz-2021-double}). The Grothendieck double construction is a strict double category whose tight category is the Grothendieck construction of \(E_{(-)} : B^{\mathsf {op}} \to  \mathcal {C}\mathsf {at}\) and whose loose category is the Grothendieck construction of the pointwise-opposite \(E_{(-)}^{\mathsf {op}} : B^{\mathsf {op}} \to  \mathcal {C}\mathsf {at}\). This latter Grothendieck construction is Spivak's \emph{generalized lens} construction (Definition 3.3 of \cite{spivak-2019-generalized}).\par{}	If \(\pi  : E \to  B\) is presented as the Grothendieck construction of an indexed category \(E_{(-)} : B^{\mathsf {op}} \to  \mathcal {C}\mathsf {at}\), then a lens \({f^{\sharp } \choose  f} : {A^- \choose  A^+} \mathrel {\mkern 3mu\vcenter {\hbox {$\shortmid $}}\mkern -10mu{\to }} {B^- \choose  B^+}\) is a diagram of the following form:
	
\begin{center}
		\begin {tikzcd}
	{A^-} & {f^{\ast }B^-} & {B^-} \\
	{A^+} & {A^+} & {B^+}
	\arrow [maps to, from=1-1, to=2-1]
	\arrow ["{f^{\sharp }}"', from=1-2, to=1-1]
	\arrow ["{\mathsf {lift}(f)}", from=1-2, to=1-3]
	\arrow [maps to, from=1-2, to=2-2]
	\arrow ["\lrcorner "{anchor=center, pos=0.125}, draw=none, from=1-2, to=2-3]
	\arrow [maps to, from=1-3, to=2-3]
	\arrow [equals, from=2-1, to=2-2]
	\arrow ["f"', from=2-2, to=2-3]
\end {tikzcd}
	\end{center}\par{}	If \(\pi  : E \to  B\) is actually the \hyperref[djm-00E4]{simple fibration} of a cartesian category \(B\), then a lens the above diagram takes the following form:
	
\begin{center}
		\begin {tikzcd}
	{A^+ \times  A^-} & {A^+ \times  B^-} & {B^+ \times  B^-} \\
	{A^+} & {A^+} & {B^+}
	\arrow ["{\pi _{A^+}}"', from=1-1, to=2-1]
	\arrow ["{(\pi _{A^+}, f^{\sharp })}"', from=1-2, to=1-1]
	\arrow ["{A^+ \times  f}", from=1-2, to=1-3]
	\arrow ["{\pi _{A^+}}"', from=1-2, to=2-2]
	\arrow ["\lrcorner "{anchor=center, pos=0.125}, draw=none, from=1-2, to=2-3]
	\arrow ["{\pi _{B^+}}", from=1-3, to=2-3]
	\arrow [equals, from=2-1, to=2-2]
	\arrow ["f"', from=2-2, to=2-3]
\end {tikzcd}
	\end{center}

	which is determined by the pair of maps \(f : A^+ \to  B^+\) and \(f^{\sharp } : A^+ \times  B^- \to  A^-\) which characterize a usual ``polymorphic, lawless'' lens --- a morphism in the fiberwise opposite of the simple fibration.\end{explication}

\begin{example}[{Lens construction of simple fibration}]\label{djm-00E6}
\par{}	Let \(C\) be a cartesian category. Then the \hyperref[djm-00B6]{lens double category} of its \hyperref[djm-00E4]{simple fibration} has as its loose morphisms the lenses \({f^{\sharp } \choose  f} : {A^{\sharp } \choose  A} \leftrightarrows  {B^{\sharp } \choose  B}\) the usual cartesian lenses (see, e.g. Definition 1.3.1.1 of \cite{jaz-2021-book}) and as tight morphisms the \emph{charts} (Definition 3.3.0.1 of \cite{jaz-2021-book}).\end{example}
\section{Loose bimodules and loose right modules}\label{djm-00IQ}\par{}	In this paper, we'll make extensive use of the basic theory of \emph{loose bimodules} between double categories. Unfortunately, this basic theory has not yet been developed, at least as far as we could tell. We therefore set out to develop it ourselves, before splitting it off into a forthcoming companion paper \cite{forthcoming-2025-loosebimodules} for reasons of audience and scope.\par{}	In this section, we will introduce what we need of the theory of loose bimodules, including a full definition of the 2-category of loose bimodules and a (very straightforward) construction of loose hom bimodules. We will then sketch the equivalence between the definition we use here (``double barrels'') and the perhaps more expected definition as pseudo-bimodules, as well as sketch our two main ways for constructing new bimodules out of old ones: specifically, \emph{restriction} and \emph{collapse}. The sketches given here are previews of the complete proofs, which will appear in our forthcoming companion paper \cite{forthcoming-2025-loosebimodules}.\subsection{Loose bimodules}\label{djm-00IY}\par{}	We begin with the definition of \emph{loose bimodule}. In order to give a slick definition of the 2-category of loose bimodules, we'll take a \emph{labelling approach} which categorifies Joyal's \emph{barrels}. A loose bimodule will be identified with its \emph{collage} \(\mathbb {M}\) together with a labelling double functor \(\ell  : \mathbb {M} \to  \mathbb {L}\mathsf {oose}\) into the \hyperref[djm-0047]{\emph{walking loose arrow}}.
\begin{definition}[{Walking loose arrow}]\label{djm-0047}
\par{}The \textbf{walking loose arrow} \(\mathbb {L}\mathsf {oose}\) is the free strict double category generated by a single loose arrow \(2_{\ell } \colon   0 \mathrel {\mkern 3mu\vcenter {\hbox {$\shortmid $}}\mkern -10mu{\to }} 1\).\end{definition}

\begin{explication}[{Walking loose arrow}]\label{ssl-001G}
\par{}The \hyperref[djm-0047]{walking loose arrow double category} \(\mathbb {L}\mathsf {oose}\) has two objects, three loose arrows (two identities along with the walking loose arrow \(2_{\ell }\)), two identity tight arrows, and three identity squares. These are all pictured below.

  \begin{center}\begin {tikzcd}
    0 & 0 & 1 & 1 \\
    0 & 0 & 1 & 1
    \arrow ["\shortmid "{marking}, equals, from=1-1, to=1-2]
    \arrow [equals, from=1-1, to=2-1]
    \arrow ["\shortmid "{marking}, from=1-2, to=1-3]
    \arrow [equals, from=1-2, to=2-2]
    \arrow ["\shortmid "{marking}, equals, from=1-3, to=1-4]
    \arrow [equals, from=1-3, to=2-3]
    \arrow [equals, from=1-4, to=2-4]
    \arrow ["\shortmid "{marking}, equals, from=2-1, to=2-2]
    \arrow ["\shortmid "{marking}, from=2-2, to=2-3]
    \arrow ["\shortmid "{marking}, equals, from=2-3, to=2-4]
  \end {tikzcd}\end{center}\end{explication}
\par{}	A double functor \(\ell  : \mathbb {M} \to  \mathbb {L}\mathsf {oose}\) therefore labels each object of \(\mathbb {M}\) by either \(0\) or \(1\) (marking it either an object of the \emph{source} of the associated loose bimdoule, or the \emph{target}), and a loose morphism \(m : x \mathrel {\mkern 3mu\vcenter {\hbox {$\shortmid $}}\mkern -10mu{\to }} y\) of \(\mathbb {M}\) can either be labelled by the loose identities \(\mathrm {id}_0 : 0 \mathrel {\mkern 3mu\vcenter {\hbox {$\shortmid $}}\mkern -10mu{\to }} 0\) and \(\mathrm {id}_1 : 1 \mathrel {\mkern 3mu\vcenter {\hbox {$\shortmid $}}\mkern -10mu{\to }} 1\), or the walking loose arrow \(2 : 0 \mathrel {\mkern 3mu\vcenter {\hbox {$\shortmid $}}\mkern -10mu{\to }} 1\) itself. In the former case, we interpret \(m\) as a loose morphism of the source of \(\mathbb {M}\); in second case, as a loose morphism of the target; and in third case, where \(\ell (m) = 2\), we consider \(m\) as a loose \emph{heteromorphism}, or an element of the bimodule.\par{}	This story is repeated again for all the tight morphisms and squares. Composition in \(\mathbb {M}\) gives a coherently unital and associative action of the loose hetermorphisms (and heterosquares --- those labelled by the tight identity of \(2 : 0 \mathrel {\mkern 3mu\vcenter {\hbox {$\shortmid $}}\mkern -10mu{\to }} 1\)) on the left and right by those loose morphisms and squares labelled by \(0\) and \(1\) respectively.\par{}	In total, we get all the structure of a loose bimodule out of a labelling double functor \(\ell  : \mathbb {M} \to  \mathbb {L}\mathsf {oose}\). This observation leads us to the following slick definition of the whole 2-category of looes bimodules.
\begin{definition}[{Loose bimodule}]\label{djm-00IS}
\par{}A \textbf{loose bimodule} is a (necessarily strict) double functor \(M : \mathbb {M} \to  \mathbb {L}\mathsf {oose}\) into the \hyperref[djm-0047]{walking loose arrow}. The 2-category \({\ell }\mathcal {B}\mathsf {imod}\) of loose modules is the (strict) slice 2-category \(\mathcal {D}\mathsf {bl} \downarrow  \mathbb {L}\mathsf {oose}\) over the \hyperref[djm-0047]{walking loose arrow}. Explicitly, \({\ell }\mathcal {B}\mathsf {imod}\) consists of:
\begin{enumerate}\item{}	Objects are loose bimodules \(M : \mathbb {M} \to  \mathbb {L}\mathsf {oose}\).
\item{}Morphisms are strictly commuting triangles:

\begin{center}
\begin {tikzcd}
	{\mathbb {M}} && {\mathbb {N}} \\
	& \mathbb {L}\mathsf {oose}
	\arrow ["{F}", from=1-1, to=1-3]
	\arrow ["{M}"', from=1-1, to=2-2]
	\arrow ["{N}", from=1-3, to=2-2]
\end {tikzcd}
\end{center}

where \(F\) is a pseudo double functor.
\item{}2-cells are tight transformations \(\alpha  \colon  F \to  G\) so that \(N \alpha  = M\), again strictly.\end{enumerate}\end{definition}
\par{}	This slick definition of loose bimodules has a number of benefits. It becomes trivial to construct loose hom bimodules, for example.
\begin{definition}[{Hom loose bimodule}]\label{djm-005F}
\par{}Let \(\mathbb {D}\) be a pseudo-double category. Its \textbf{loose hom bimodule} \(\mathsf {Hom}^{l}(\mathbb {D})\) is defined to be the product projection \(  \mathbb {D} \times  \mathbb {L}\mathsf {oose} \to  \mathbb {L}\mathsf {oose}\).\par{}This gives us a cartesian 2-functor
\[\mathsf {Hom}^{l} := \mathbb {L}\mathsf {oose} \times  (-) : \mathcal {D}\mathsf {bl} \to  {\ell }\mathcal {B}\mathsf {imod}\]\par{}It is often convenient to notate the loose home bimodule \(\mathsf {Hom}^{l}(\mathbb {D}) \colon  \mathbb {D} \mathrel {\mkern 3mu\vcenter {\hbox {$\shortmid $}}\mkern -10mu{\to }} \mathbb {D}\) simply by \[\mathbb {D}(-,-) \colon  \mathbb {D} \mathrel {\mkern 3mu\vcenter {\hbox {$\shortmid $}}\mkern -10mu{\to }} \mathbb {D}.\]\end{definition}
\par{}	In \(\mathbb {D}(-,-)\), we have two full copies of \(\mathbb {D}\) appearing as \(\mathbb {D} \times  \{0\}\) and \(\mathbb {D} \times  \{1\}\), respectively, acting on the loose morphisms \(m : x \mathrel {\mkern 3mu\vcenter {\hbox {$\shortmid $}}\mkern -10mu{\to }} y\) of \(\mathbb {D}\) appearing as \((m, 2) : (x, 0) \mathrel {\mkern 3mu\vcenter {\hbox {$\shortmid $}}\mkern -10mu{\to }} (y, 1)\) by loose composition (and similarly for squares).\par{}	We will need to know a bit more about \(\mathcal {D}\mathsf {bl}\); namely, we need to know that it is \hyperref[djm-00AN]{cartesian}. Since it is a slice 2-category, its products may be constructed by taking pullback; but pullbacks of \emph{pseudo-}double functors do not always exist. In this case, however, \(\mathbb {L}\mathsf {oose}\) is just so strict that nothing can really go wrong.
\begin{lemma}[{Pullbacks of maps into \(\mathbb {L}\mathsf {oose}\)}]\label{djm-00IT}
\par{}Any two double functors 

\begin{center}
	\begin {tikzcd}
	& {\mathbb {D}} \\
	{\mathbb {C}} & {\mathbb {L}\mathsf {oose}}
	\arrow ["F", from=1-2, to=2-2]
	\arrow ["G"', from=2-1, to=2-2]
\end {tikzcd}
\end{center}

admit a pullback in the 2-category \(\mathcal {D}\mathsf {bl}\) of double categories and pseudo-functors, and it is constructed by taking the pullback of all underlying categories (of tight morphisms and of squares, respectively).\begin{proof}\par{}	This is straightforward to verify directly, but let's see an abstract argument.\par{}			In \cite{lack-2005-limits}, Lack shows that comma objects out of a lax morphism and into a strict morphism always exist in 2-categories of algebras for 2-monads and lax morphisms and are constructed in the underlying 2-category. Similarly, iso-commas out of a pseudo-morphism and into a strict morphism always exist in 2-categories of algebras for 2-monads and pseudo-monads, and are again constructed in the underlying 2-category. Since double categories are 2-monadic over graphs of categories, and since all functors into \(\mathbb {L}\mathsf {oose}\) are strict (since there is simply no room to be pseudo), we conclude that the iso-comma of the above cospan exists. But since \(\mathbb {L}\mathsf {oose}\) has only identity tight morphisms and squares, that iso-comma must already be a strict pullback.\end{proof}\end{lemma}
\par{}	As a corollary, the double category of \hyperref[djm-00IS]{loose bimodules} is \hyperref[djm-00AN]{cartesian}. The product of loose bimodules is straightforward: we take pairs of all things with the same labels.
\begin{corollary}[{}]\label{djm-00IU}
\par{}	The 2-category \({\ell }\mathcal {B}\mathsf {imod}\) of \hyperref[djm-00IS]{loose bimodules} is \hyperref[djm-00AN]{cartesian}, with the cartesian product given by taking the \hyperref[djm-00IT]{pullback of maps into \(\mathbb {L}\mathsf {oose}\)}.\end{corollary}

\begin{definition}[{Symmetric monoidal loose bimodules}]\label{ssl-0019}
\par{}Define the 2-category of \textbf{symmetric monoidal loose bimodules} to be \(\mathcal {S}\mathsf {M}{({\ell }\mathcal {B}\mathsf {imod})}\), the 2-category of .\end{definition}
\par{}	With \hyperref[djm-00IT]{Lemma \ref{djm-00IT}}, we can also formally define the source and target of a loose bimodule.
\begin{definition}[{Source and target of loose bimodule}]\label{djm-004A}
\par{}Let \(M \colon  \mathbb {M} \to  \mathbb {L}\mathsf {oose}\) be a . The \textbf{source} \(M_0\) of \(M\) is its pullback along \(0 : \bullet  \to  \mathbb {L}\mathsf {oose}\). Its \textbf{target} \(M_1\) is its pullback along \(1 : \bullet  \to  \mathbb {L}\mathsf {oose}\).\par{}Together, taking the source and target of a loose bimodule gives a cartesian 2-functor
\[(()_0, ()_1) \colon  {\ell }\mathcal {B}\mathsf {imod} \to  \mathcal {D}\mathsf {bl} \times  \mathcal {D}\mathsf {bl}\]\end{definition}

\begin{notation}[{Loose bimodule}]\label{ssl-0016}
\par{}We often denote a  \(M \colon  \mathbb {M} \to  \mathbb {L}\mathsf {oose}\) as a proarrow \[\mathbb {M}  \colon  \mathbb {D}_0 \mathrel {\mkern 3mu\vcenter {\hbox {$\shortmid $}}\mkern -10mu{\to }} \mathbb {D}_1\] where:
  \begin{itemize}\item{}\(\mathbb {D}_0 \) is the \hyperref[djm-004A]{source}  \(M_0\).
    \item{}\(\mathbb {D}_1 \) is the \hyperref[djm-004A]{target}  \(M_1\).\end{itemize}\end{notation}
\par{}	Finally, we introduce the notion of a loose \emph{right} module, which is simply a loose bimodule whose \hyperref[djm-004A]{source} is terminal.
\begin{definition}[{2-Category of loose right modules}]\label{djm-004M}
\par{}Define the 2-category of \textbf{loose right modules} to be the following pullback of 2-categories:\begin{center} \begin {tikzcd}
	{\ell }\mathcal {M}\mathsf {od}_{\mathsf {r}} & {\ell }\mathcal {B}\mathsf {imod} \\
	\mathcal {D}\mathsf {bl} & {\mathcal {D}\mathsf {bl} \times  \mathcal {D}\mathsf {bl}}
	\arrow [from=1-1, to=1-2]
	 \arrow [from=1-1, to=2-1]
	\arrow ["{(()_0, ()_1)}", from=1-2, to=2-2]
	\arrow [""{name=0, anchor=center, inner sep=0}, "{(\bullet , \mathrm {id})}"', from=2-1, to=2-2]
	\arrow ["\lrcorner "{anchor=center, pos=0.125}, draw=none, from=1-1, to=0]
\end {tikzcd}
\end{center}\end{definition}

\begin{explication}[{Loose right module}]\label{ssl-003S}
\par{}A loose right module \(\mathbb {M} \colon  \bullet  \mathrel {\mkern 3mu\vcenter {\hbox {$\shortmid $}}\mkern -10mu{\to }} \mathbb {D}_1\) is a loose bimodule whose source is the  and thus has a trivial left action. The \hyperref[djm-004A]{target} \(\mathbb {D}_1\) acts on the \hyperref[ssl-0036]{carrier} \(\mathsf {Car}(\mathbb {M})\) on the right. We say \(\mathbb {M}\) is a loose right module \emph{over} \(\mathbb {D}_1\).\end{explication}

\begin{proposition}[{The 2-category \({\ell }\mathcal {M}\mathsf {od}_{\mathsf {r}}\) is cartesian}]\label{ssl-001D}
\par{}The 2-category \({\ell }\mathcal {M}\mathsf {od}_{\mathsf {r}}\) is cartesian, and products are constructed in \({\ell }\mathcal {B}\mathsf {imod}\) (as pullbacks over \(\mathbb {L}\mathsf {oose}\)).\begin{proof}\par{}The source and target maps are cartesian because they are constructed via pullbacks. Since \({\ell }\mathcal {M}\mathsf {od}_{\mathsf {r}}\) is a pullback of cartesian 2-categories along cartesian maps, it is cartesian.\end{proof}\end{proposition}
\subsection{Loose bimodules as pseudo-bimodules}\label{djm-00IZ}\par{}	Let's see how the definition of \hyperref[djm-00IS]{loose bimodule} given in \hyperref[djm-00IS]{Definition \ref{djm-00IS}} relates to perhaps more expected definition of a pseudo-bimodule acted on the left and right by pseudo-categories in \(\mathcal {C}\mathsf {at}\). We will only sketch the equivalence here; a full proof will appear in \cite{forthcoming-2025-loosebimodules}.\par{}	First, we note that to any \hyperref[djm-00IS]{loose bimodule}, we may associate a category which we call its \emph{carrier}.
\begin{definition}[{Carrier of a loose bimodule}]\label{ssl-0036}
\par{}Let \(\ell  : \mathbb {M} \to  \mathbb {L}\mathsf {oose}\) be a loose bimodule. Its carrier \(\mathsf {Car} (\mathbb {M})\) is the category whose: 
  \begin{itemize}\item{}Objects are loose morphisms of \(\mathbb {M}\) living over the walking loose arrow. We refer to such loose morphisms as \textbf{heteromorphisms}.
    \item{}Morphisms are squares of \(\mathbb {M}\) living over the identity square of the walking loose arrow.\end{itemize}\par{}Writing the loose bimodule \(\mathbb {M} \to  \mathbb {L}\mathsf {oose}\) as a \hyperref[ssl-0016]{proarrow between its source and target} \[\mathbb {M} \colon  \mathbb {D}_0 \mathrel {\mkern 3mu\vcenter {\hbox {$\shortmid $}}\mkern -10mu{\to }} \mathbb {D}_1\] we have the following explication.\par{}An object of \(\mathsf {Car}(\mathbb {M})\) consists of a triple \[(d_0, d_1, m \colon  d_0 \mathrel {\mkern 3mu\vcenter {\hbox {$\shortmid $}}\mkern -10mu{\to }} d_1)\] where \(d_0\) is an object of \(\mathbb {D}_0\), \(d_1\) is an object of \(\mathbb {D}_1\), and \(m\) is a loose morphism in \(\mathbb {M}\) from \(d_0\) to \(d_1\).\par{}A morphism \((d_0, d_1, m \colon  d_0 \mathrel {\mkern 3mu\vcenter {\hbox {$\shortmid $}}\mkern -10mu{\to }} d_1) \to  (d_0', d_1', m' \colon  d_0' \mathrel {\mkern 3mu\vcenter {\hbox {$\shortmid $}}\mkern -10mu{\to }} d_1')\) consists of a tight morphism \(f_0 \colon  d_0 \to  d_0'\) in \(\mathbb {D}_0\), a tight morphism \(f_1 \colon  d_1 \to  d_1'\) in the \(\mathbb {D}_0\),  and a square \( \alpha  \colon  m \Rightarrow  m'\) in \(\mathbb {M}\):

  \begin{center}
    \begin {tikzcd}
    {d_0} & {d_1} \\
    {d_0'} & {d_1'} \\
    {}
    \arrow [""{name=0, anchor=center, inner sep=0}, "m", "\shortmid "{marking}, from=1-1, to=1-2]
    \arrow ["{f_0}"', from=1-1, to=2-1]
    \arrow ["{f_1}", from=1-2, to=2-2]
    \arrow [""{name=1, anchor=center, inner sep=0}, "{m'}"', "\shortmid "{marking}, from=2-1, to=2-2]
    \arrow ["\alpha ", shorten <=4pt, shorten >=4pt, Rightarrow, from=0, to=1]
  \end {tikzcd}
  \end{center}\end{definition}
\par{}	The \hyperref[ssl-0036]{carrier} of a \hyperref[djm-00IS]{loose bimodule} is spanned over the category of tight arrows of its \hyperref[djm-004A]{source and target}.
\begin{proposition}[{Projection from the carrier of a loose bimodule to its source and target}]\label{ssl-003D}
\par{}Given a \hyperref[djm-00IS]{loose bimodule} \(\mathbb {M} \colon  \mathbb {D}_0 \mathrel {\mkern 3mu\vcenter {\hbox {$\shortmid $}}\mkern -10mu{\to }} \mathbb {D}_1\), there are projections from its \hyperref[ssl-0036]{carrier} to the tight category of its \hyperref[djm-004A]{source and target} of \(\mathbb {M}\). This gives a span:  
  
\begin{center}
\begin {tikzcd}
	& {\mathsf {Car} (\mathbb {M})} \\
	{ \mathsf {Tight}({\mathbb {D}_0})} && { \mathsf {Tight}({\mathbb {D}_1})}
	\arrow ["{\mathsf {dom}}"', from=1-2, to=2-1]
	\arrow ["{\mathsf {codom}}", from=1-2, to=2-3]
\end {tikzcd}    
  \end{center}\end{proposition}
\par{}	With regard to these projections, the carrier \(\mathsf {Car}(\mathbb {M})\) picks up left and right actions by \(\mathbb {D}_0\) and \(\mathbb {D}_1\) by composition in \(\mathbb {M}\):
	
\begin{center}
		\begin {tikzcd}[column sep = tiny]
	{ \mathsf {Tight}({\mathbb {D}_0})} && { \mathsf {Tight}({\mathbb {D}_0})} && { \mathsf {Tight}({\mathbb {D}_1})} && { \mathsf {Tight}({\mathbb {D}_1})} \\
	& {\mathsf {Loose}({\mathbb {D}_0})} && {\mathsf {Car} (\mathbb {M})} && {\mathsf {Loose}({\mathbb {D}_1})} \\
	&& \bullet  && \bullet  \\
	&&& \bullet  \\
	&&& {\mathsf {Car} (\mathbb {M})} \\
	{ \mathsf {Tight}({\mathbb {D}_0})} &&&&&& { \mathsf {Tight}({\mathbb {D}_1})}
	\arrow [equals, from=1-1, to=6-1]
	\arrow [equals, from=1-7, to=6-7]
	\arrow [from=2-2, to=1-1]
	\arrow [from=2-2, to=1-3]
	\arrow ["{\mathsf {dom}}"', from=2-4, to=1-3]
	\arrow ["{\mathsf {codom}}", from=2-4, to=1-5]
	\arrow [from=2-6, to=1-5]
	\arrow [from=2-6, to=1-7]
	\arrow ["\lrcorner "{anchor=center, pos=0.125, rotate=135}, draw=none, from=3-3, to=1-3]
	\arrow [from=3-3, to=2-2]
	\arrow [from=3-3, to=2-4]
	\arrow ["\lrcorner "{anchor=center, pos=0.125, rotate=135}, draw=none, from=3-5, to=1-5]
	\arrow [from=3-5, to=2-4]
	\arrow [from=3-5, to=2-6]
	\arrow ["\lrcorner "{anchor=center, pos=0.125, rotate=135}, draw=none, from=4-4, to=2-4]
	\arrow [from=4-4, to=3-3]
	\arrow [from=4-4, to=3-5]
	\arrow [dashed, from=4-4, to=5-4]
	\arrow [from=5-4, to=6-1]
	\arrow [from=5-4, to=6-7]
\end {tikzcd}
	\end{center}\par{}	We will often find it more convenient to think of \hyperref[djm-00IS]{loose bimodules} \(\mathbb {M}\) as being displayed categories \(\mathsf {Car}(\mathbb {M})\) acted on the left and right by double categories. This will make the upcoming constructions of \emph{restriction} and \emph{collapse} more easy to understand as well.\par{}	Let's now sketch the relationship between \hyperref[djm-00IS]{loose bimodules} presented as double barrels (\hyperref[djm-00IS]{Definition \ref{djm-00IS}}) and as pseudo-bimodules between pseudo-categories in \(\mathcal {C}\mathsf {at}\); this will be a sneak peak at \cite{forthcoming-2025-loosebimodules}, though we do not follow exactly this strategy there.\par{}	As we saw in \hyperref[djm-00IR]{Definition \ref{djm-00IR}}, we may identify the 2-category \(\mathcal {D}\mathsf {bl}\) of double categories with the 2-category \(\mathcal {M}\mathsf {od}_p(\mathcal {A}_\mathsf {PseudoCat}; \mathcal {C}\mathsf {at})\) of models of the \(\mathcal {F}\)-sketch of pseudo-categories valued in the 2-category \(\mathcal {C}\mathsf {at}\) of categories. We may consider the strict double category \(\mathbb {L}\mathsf {oose} : \mathcal {A}_{\mathsf {PseudoCat}} \to  \mathcal {C}\mathsf {at}\) as such a model; it's 2-category of elements \(\mathcal {E}\mathsf {l}(\mathbb {L}\mathsf {oose})\) may be endowed with the structure of an \(\mathcal {F}\)-sketch by taking its marked cones to be those that project down to marked cones in \(\mathcal {A}_{\mathsf {PseudoCat}}\).\par{}	An inspection of \(\mathcal {E}\mathsf {l}(\mathbb {L}\mathsf {oose})\) reveals that it gives a \(\mathcal {F}\)-sketch for pseudo-bimodules between pseudo-categories. By adapting closure properties of 2-fibrations to \(\mathcal {F}\)-categories, we can prove a categorification of the usual theorem relating slice categories of presheaf categories to presheaves on categories of elements:
	\[
		\mathcal {M}\mathsf {od}_p(\mathcal {A}_{\mathsf {PseudoCat}};\mathcal {C}\mathsf {at}) \downarrow  \mathbb {L}\mathsf {oose} \simeq  \mathcal {M}\mathsf {od}_p(\mathcal {E}\mathsf {l}(\mathbb {L}\mathsf {oose}); \mathcal {C}\mathsf {at}).
	\]
	This gives the equivalence between the definition of loose bimodule we gave in \hyperref[djm-00IS]{Definition \ref{djm-00IS}} and pseudo-bimodules between pseudo-categories.\subsection{Collapse and restriction of loose bimodules}\label{djm-00J0}\par{}	In this section, we'll sketch two constructions of new loose bimodules from old: \emph{collapse} and \emph{restriction}. All of the loose modules of systems we construct in this paper will be done by restricting \hyperref[djm-005F]{loose hom bimodules}, and then potentially collapsing them.\par{}	We'll begin with the collapse, since it is more straightforward. We then describe the restriction of loose bimodules, and in particular its pseudo-functoriality over a 2-category of \emph{niches}.\subsubsection{Collapse of loose bimodules}\label{djm-00J1}\par{}	Given a \hyperref[djm-00IS]{loose bimodule} \(\mathbb {M} : \mathbb {D}_0 \mathrel {\mkern 3mu\vcenter {\hbox {$\shortmid $}}\mkern -10mu{\to }} \mathbb {D}_1\) we may consider it as a pseudo-bimodule \hyperref[ssl-003D]{spanned between its source \(\mathbb {D}_0\) and its target \(\mathbb {D}_1\)}:
	
\begin{center}
		\begin {tikzcd}
	& {\mathsf {Car} (\mathbb {M})} \\
	{ \mathsf {Tight}({\mathbb {D}_0})} && { \mathsf {Tight}({\mathbb {D}_1})}
	\arrow ["{\mathsf {dom}}"', from=1-2, to=2-1]
	\arrow ["{\mathsf {codom}}", from=1-2, to=2-3]
\end {tikzcd}    
	\end{center}

	The \textbf{collapse} \(\mathfrak {c}\mathbb {M}\) of \(\mathbb {M}\) is the \hyperref[djm-004M]{loose right module} given by replacing \(\mathbb {D}_0\) by the terminal double category:
	
\begin{center}
		\begin {tikzcd}
	& {\mathsf {Car} (\mathbb {M})} \\
	\bullet  && { \mathsf {Tight}({\mathbb {D}_1})}
	\arrow ["{!}"', from=1-2, to=2-1]
	\arrow ["{\mathsf {codom}}", from=1-2, to=2-3]
\end {tikzcd}
	\end{center}

	In particular, the collapse \(\mathfrak {c}\mathbb {M}\) of a loose bimodule \(\mathbb {M}\) has exactly the same \hyperref[ssl-003D]{carrier}.\par{}	This construction is a bit awkward to express with double barrels but evidently 2-functorial when viewed as a construction on pseudo-bimodules. We therefore have the following theorem:
\begin{theorem}[{Collapse of loose bimodule}]\label{djm-00J2}
\par{}	The \textbf{collapse} of a \hyperref[djm-00IS]{loose bimodule} into a \hyperref[djm-004M]{loose right module} gives a cartesian 2-functor 
	\[\mathfrak {c} : {\ell }\mathcal {B}\mathsf {imod} \to  {\ell }\mathcal {M}\mathsf {od}_{\mathsf {r}}.\]\end{theorem}
\subsubsection{Restriction of loose bimodules}\label{djm-00J3}\par{}	In this section, we will sketch the \textbf{restriction} of loose bimodules along pseudo-double functors. The construction itself itself is rather straightforward. Suppose we are given a \hyperref[djm-00IS]{loose bimodule} \(\mathbb {M} : \mathbb {D}_0 \mathrel {\mkern 3mu\vcenter {\hbox {$\shortmid $}}\mkern -10mu{\to }} \mathbb {D}_1\) and two double functors \(F_0 : \mathbb {E}_0 \to  \mathbb {D}_0\) and \(F_1 : \mathbb {E}_1 \to  \mathbb {D}_1\). We may then form the carrier of the restricted loose bimodule \(\mathbb {M}(F_0, F_1)\) as the limit
	
\begin{center}
\begin {tikzcd}[cramped]
      {\mathsf {Tight}({\mathbb {E}_0})} && {\mathsf {Car}(\mathbb {M})} && {\mathsf {Tight}({\mathbb {E}_1})} \\
      & {\mathsf {Tight}({\mathbb {D}_0})} && {\mathsf {Tight}({\mathbb {D}_1})}
      \arrow ["{\mathsf {Tight}({F_0})}"', from=1-1, to=2-2]
      \arrow ["{\mathsf {dom}}", from=1-3, to=2-2]
      \arrow ["{\mathsf {codom}}"', from=1-3, to=2-4]
      \arrow ["{\mathsf {Tight}({F_1})}", from=1-5, to=2-4]
    \end {tikzcd}
	\end{center}

	Explicitly, we will have \(\mathsf {Car}(\mathbb {M}(F_0, F_1))\) defined as follows:
	\begin{itemize}\item{}Objects are triples \((e_0 : {\mathsf {ob}}(\mathbb {E}_0), e_1: {\mathsf {ob}}(\mathbb {E}_1), m \colon  F_0 e_0 \mathrel {\mkern 3mu\vcenter {\hbox {$\shortmid $}}\mkern -10mu{\to }} F_1 e_1).\)
    \item{}A morphism \((e_0, e_1, m) \to  (e_0', e_1', m')\) is a triple \((f_0:e_0\to  e_0', f_1:e_1\to  e_1', \alpha )\) where \(\alpha \) is a square in \(\mathbb {M}\) of the following form:

  \begin{center}
        \begin {tikzcd}
          {F_0e_0} & {F_1e_1} \\
          {F_0e_0'} & {F_1e_1} \\
          {}
          \arrow [""{name=0, anchor=center, inner sep=0}, "m", "\shortmid "{marking}, from=1-1, to=1-2]
          \arrow ["{F_0f_0}"', from=1-1, to=2-1]
          \arrow ["{F_1f_1}", from=1-2, to=2-2]
          \arrow [""{name=1, anchor=center, inner sep=0}, "{m'}"', "\shortmid "{marking}, from=2-1, to=2-2]
          \arrow ["\alpha ", shorten <=4pt, shorten >=4pt, Rightarrow, from=0, to=1]
        \end {tikzcd}
      \end{center}\end{itemize}
	The action of \(\mathbb {E}_0\) and \(\mathbb {E}_1\) is by first applying \(F_i\), and then acting in \(\mathbb {M}\).\par{}	What is more interesting is the pseudo-functoriality of the above construction. We will need a very particular sort of pseudo-functoriality in a domain 2-category which we will describe now.
\begin{definition}[{The 2-category of niches}]\label{djm-00BS}
\par{}	Consider the following pullback of 2-categories:

  \begin{center}
			\begin {tikzcd}
				{\mathcal {N}\mathsf {}} & {\ell }\mathcal {B}\mathsf {imod} \\
				{\mathsf {2}\mathcal {C}\mathsf {at}^c([1],\mathcal {D}\mathsf {bl})^2}& {\mathcal {D}\mathsf {bl}^2}
				\arrow [from=1-1, to=1-2]
				\arrow [from=1-1, to=2-1]
				\arrow [from=1-2, to=2-2]
				\arrow [""{name=0, anchor=center, inner sep=0}, "{{\mathsf {codom}}^2}"', from=2-1, to=2-2]
				\arrow [""{name=0p, anchor=center, inner sep=0}, phantom, from=2-1, to=2-2, start anchor=center, end anchor=center]
				\arrow ["\lrcorner "{anchor=center, pos=0.125}, draw=none, from=1-1, to=0p]
			\end {tikzcd}
		\end{center}

Here \(\mathsf {2}\mathcal {C}\mathsf {at}^c\) is the 3-category of 2-categories, 2-functors, colax natural transformations, 
and modifications, so that \(\mathsf {2}\mathcal {C}\mathsf {at}^c([1],\mathcal {D}\mathsf {bl})\) is the 2-category of pseudo double functors and 
colax-commuting squares. (This will be explicated below.)\par{}We define the 2-category \(\mathcal {N}\mathsf {iche}\) to be the wide and locally full sub-2-category of \(\mathcal {N}\mathsf {}\)
  on those 1-morphisms whose component in \(\mathsf {2}\mathcal {C}\mathsf {at}^c([1],\mathcal {D}\mathsf {bl})^2\) lies in 
  \(\mathsf {2}\mathcal {C}\mathsf {at}^{\mathsf {colax}}(\Delta [1], \mathcal {D}\mathsf {bl})_{\mathsf {conj}} \times  \mathsf {2}\mathcal {C}\mathsf {at}^{\mathsf {colax}}(\Delta [1], \mathcal {D}\mathsf {bl})_{\mathsf {comp}}\), where 
  \(\mathsf {2}\mathcal {C}\mathsf {at}^{\mathsf {colax}}(\Delta [1], \mathcal {D}\mathsf {bl})_{\mathsf {conj}}\)  (resp. \(\mathsf {2}\mathcal {C}\mathsf {at}^{\mathsf {colax}}(\Delta [1], \mathcal {D}\mathsf {bl})_{\mathsf {comp}}\)) is the .
  
  \par{}	We also define the 2-category \(\mathcal {N}\mathsf {iche}^{\cong }\) to be the sub-2-category of \(\mathcal {N}\mathsf {iche}\) consisting of all objects, but only the 1-cells for which the colaxators are isomorphisms (and all 2-cells).
\begin{explication}[{The 2-category \(\mathcal {N}\mathsf {iche}\)}]\label{ssl-003L}
\par{}Explicitly:
  \begin{enumerate}\item{}An object of \(\mathcal {N}\mathsf {iche}\) consists of a loose bimodule \(M : \mathbb {M} \to  \mathbb {L}\mathsf {oose}\) together with two double functors \(F_0 : \mathbb {E}_0 \to  M_0\) and \(F_1 : \mathbb {E}_1 \to  M_1\). We will call an object of \(\mathcal {N}\mathsf {iche}\) a \textbf{niche}, and can draw it in this way:
      
\begin{center}
        \begin {tikzcd}
          {\mathbb {E}_0} & {\mathbb {E}_1} \\
          {M_0} & {M_1}
          \arrow ["{F_0}"', from=1-1, to=2-1]
          \arrow ["{F_1}",from=1-2, to=2-2]
          \arrow ["M"', "{\shortmid }"{marking}, from=2-1, to=2-2]
        \end {tikzcd}
      \end{center}
    \item{}A 1-cell of the above pullback consists of a  map \(a : M_0 \to  M_1\) of modules and double functors \(f_i \colon  \mathbb {E}_{0i} \to  \mathbb {E}_{1i}\), together with two colax-commuting squares of double functors of the following form:
      
\begin{center}
        \begin {tikzcd}
          {\mathbb {E}_{0i}} & {\mathbb {E}_{1i}} \\
          {M_{0i}} & {M_{1i}}
          \arrow ["{f_i}", from=1-1, to=1-2]
          \arrow ["{F_{0i}}"', from=1-1, to=2-1]
          \arrow ["{F_{1i}}", from=1-2, to=2-2]
          \arrow ["{\overline {f}_i}"{description}, Rightarrow, from=2-1, to=1-2]
          \arrow ["M"', "\shortmid "{marking}, from=2-1, to=2-2]
        \end {tikzcd}
      \end{center}

      where \(\overline {f}_0\) is a \hyperref[djm-00EX]{companion commuter transformation} and \(\overline {f}_1\) is a \hyperref[djm-00EX]{conjoint commuter transformation}.
    \item{}A 2-cell of the above pullback consists of tight transformations \(\alpha  : a_0 \to  a_1\), and \(\phi _{i} : f_{0i} \to  f_{1i}\) so that the following equations hold:
      
\begin{center}
        \begin {tikzcd}
          {\mathbb {E}_{0i}} & {\mathbb {E}_{1i}} && {\mathbb {E}_{00}} & {\mathbb {E}_{1i}} \\
          {\mathbb {E}_{0i}} & {\mathbb {E}_{1i}} & {=} & {M_{0i}} & {M_{1i}} \\
          {M_{0i}} & {M_{1i}} && {M_{0i}} & {M_{1i}}
          \arrow ["{f_{1i}}", from=1-1, to=1-2]
          \arrow [Rightarrow, no head, from=1-1, to=2-1]
          \arrow [Rightarrow, no head, from=1-2, to=2-2]
          \arrow ["{f_{1i}}", from=1-4, to=1-5]
          \arrow ["{F_{0i}}"', from=1-4, to=2-4]
          \arrow ["{F_{1i}}", from=1-5, to=2-5]
          \arrow ["{\phi _i}", shorten <=5pt, shorten >=5pt, Rightarrow, from=2-1, to=1-2]
          \arrow ["{f_{0i}}", from=2-1, to=2-2]
          \arrow ["{F_{0i}}"', from=2-1, to=3-1]
          \arrow ["{F_{1i}}", from=2-2, to=3-2]
          \arrow ["{\overline {f}_{1i}}", shorten <=5pt, shorten >=5pt, Rightarrow, from=2-4, to=1-5]
          \arrow ["{a_{1i}}", from=2-4, to=2-5]
          \arrow [Rightarrow, no head, from=2-4, to=3-4]
          \arrow [Rightarrow, no head, from=2-5, to=3-5]
          \arrow ["{\overline {f}_{0i}}", shorten <=5pt, shorten >=5pt, Rightarrow, from=3-1, to=2-2]
          \arrow ["{a_{0i}}"', from=3-1, to=3-2]
          \arrow ["{\alpha _i}", shorten <=6pt, shorten >=6pt, Rightarrow, from=3-4, to=2-5]
          \arrow ["{a_{0i}}"', from=3-4, to=3-5]
        \end {tikzcd}
        \end{center}\end{enumerate}\end{explication}
\end{definition}
\par{}	We can now express the pseudo-functoriality of restriction. Because we will make reference to the action of restriction on 1-cells, we include their description here. We remark again the a full proof will appear in our \cite{forthcoming-2025-loosebimodules}.
\begin{theorem}[{Restriction of loose bimodules is a pseudo-functor}]\label{djm-00J4}
\par{}	Restriction of loose bimodules extends to a cartesian pseudo-functor
\[\mathsf {Res} : \mathcal {N}\mathsf {iche} \to  {\ell }\mathcal {B}\mathsf {imod}.\]\begin{proof}\par{}	We will only describe the construction itself, foregoing its proof to \cite{forthcoming-2025-loosebimodules}.	\begin{proof}\par{}Consider the following 1-cell \((f_0,\bar  f_0,f_1,\bar  f_1,a) : (F_{00},F_{01},M_0) \to  (F_{10},F_{11},M_1)\) of \(\mathcal {N}\mathsf {iche}\), which we'll abusively denote as \(a\).

  \begin{center}
        \begin {tikzcd}
        {\mathbb {E}_{00}} & {\mathbb {E}_{01}} \\
        {M_{00}} & {M_{01}} \\
        & {\mathbb {E}_{10}} & {\mathbb {E}_{11}} \\
        & {M_{10}} & {M_{11}}
        \arrow ["{F_{00}}"', from=1-1, to=2-1]
        \arrow ["{f_0}"{description, pos=0.7}, curve={height=6pt}, from=1-1, to=3-2]
        \arrow ["{F_{01}}"', from=1-2, to=2-2]
        \arrow [""{name=0, anchor=center, inner sep=0}, "{f_1}"{description}, curve={height=-6pt}, from=1-2, to=3-3]
        \arrow [""{name=1, anchor=center, inner sep=0}, "{M_0}", "\shortmid "{marking}, from=2-1, to=2-2]
        \arrow [""{name=2, anchor=center, inner sep=0}, "{a_0}"{description}, curve={height=6pt}, from=2-1, to=4-2]
        \arrow ["{a_1}"{description}, curve={height=-6pt}, from=2-2, to=4-3]
        \arrow ["{F_{10}}", from=3-2, to=4-2]
        \arrow ["{F_{11}}", from=3-3, to=4-3]
        \arrow [""{name=3, anchor=center, inner sep=0}, "{M_1}", "\shortmid "{marking}, from=4-2, to=4-3]
        \arrow ["{\bar  f_0}", shorten <=3pt, Rightarrow, from=2, to=3-2]
        \arrow ["a", curve={height=-6pt}, shorten <=10pt, shorten >=10pt, Rightarrow, from=1, to=3]
        \arrow ["{\bar  f_1}"', shorten >=3pt, Rightarrow, from=2-2, to=0]
      \end {tikzcd}
  \end{center}

	We define a loose bimodule map \(\mathsf {Res}(a):\mathbb {M}_0(F_{00}, F_{01}) \to  \mathbb {M}_1(F_{10}, F_{11})\) as follows:
  \begin{enumerate}\item{}\textbf{(Endpoints)} On the \(\mathbb {E}_{0i}\), we take \(\mathsf {Res}(a)\) to be \(f_{i}\).
    \item{}\textbf{(Objects of carrier)} Given \((e_{00},e_{01},m_0 : F_{00}e_{00} \mathrel {\mkern 3mu\vcenter {\hbox {$\shortmid $}}\mkern -10mu{\to }} F_{01}e_{01})\) in \(\mathbb {M}_0(F_{00},F_{01})\), recall that we have the tight morphism \(\bar {f}_0e_{00}:a_0F_{00}(e_{00})\to  F_{10}f_0(e_{00})\), 
      which has a conjoint \((\overline {f}_0e_{00})^{<}\) as part of the definition of \(\mathcal {N}\mathsf {iche}\). Similarly, \(\bar {f}_1 e_{01}\) has a companion \((\overline {f}_1e_{01})^{>}\). These loosenings will be used to modify 
      \(am_0: a_0F_{00}e_{00}\mathrel {\mkern 3mu\vcenter {\hbox {$\shortmid $}}\mkern -10mu{\to }} a_1 F_{01}e_{01}\) to have the correct signature. 
      
      We therefore define \(\mathsf {Res}(a)(e_{00},e_{01},m_0)\) as \((f_0(e_{00}),f_1(e_{01}),m_0')\) where \(m_0':F_{01}f_0(e_{00})\mathrel {\mkern 3mu\vcenter {\hbox {$\shortmid $}}\mkern -10mu{\to }} F_{11}f_1(e_{01})\) is the following composite in \(\mathbb {M}_1\):
      
\begin{center}
        \begin {tikzcd}
        {F_{10}f_0e_{00}} & {a_0F_{00}e_{00}} & {a_1F_{01}e_{01}} & {F_{11}f_1e_{01}}
        \arrow ["{(\overline {f}_0e_{00})^{<}}", "\shortmid "{marking}, from=1-1, to=1-2]
        \arrow ["{am_0}", "\shortmid "{marking}, from=1-2, to=1-3]
        \arrow ["{(\overline {f}_1e_{10})^{>}}", "\shortmid "{marking}, from=1-3, to=1-4]
      \end {tikzcd}
      \end{center}
    \item{}\textbf{(Morphisms of carrier)} Consider a square \((\mu _0,\mu _{1},\mu ):(e_{00},e_{01},m_0) \to  (e_{00}',e_{01}',m_0')\) in \(\mathbb {M}_0(F_{00},F_{01})\), 
      where \(\mu \) appears as below:
      
\begin{center}
       \begin {tikzcd}
        {F_{00}e_{00}} & {F_{01}e_{01}} \\
        {F_{00}e_{00}'} & {F_{01}e_{01}'}
        \arrow [""{name=0, anchor=center, inner sep=0}, "{m_0}", "\shortmid "{marking}, from=1-1, to=1-2]
        \arrow ["{F_{00}\mu _0}"', from=1-1, to=2-1]
        \arrow ["{F_{01}\mu _1}", from=1-2, to=2-2]
        \arrow [""{name=1, anchor=center, inner sep=0}, "{m_1}"', "\shortmid "{marking}, from=2-1, to=2-2]
        \arrow ["\mu "{description}, shorten <=4pt, shorten >=4pt, Rightarrow, from=0, to=1]
      \end {tikzcd}
      \end{center}

      Then we define \(\mathsf {Res}(a)(\mu )=(f_0\mu _0,f_1\mu _1,\mu ')\) where \(\mu '\) is the following square:
      
\begin{center}
        \begin {tikzcd}
          {F_{10}f_0e_{00}} & {a_0F_{00}e_{00}} & {a_1F_{01}e_{01}} & {F_{11}f_1e_{01}} \\
          {F_{10}f_0e_{00}'} & {a_0F_{00}e_{00}'} & {a_1F_{01}e_{01}'} & {F_{11}f_1e_{01}'}
          \arrow [""{name=0, anchor=center, inner sep=0}, "{(\overline {f}_0e_{00})^{<}}", "\shortmid "{marking}, from=1-1, to=1-2]
          \arrow ["{F_{10}f_0\mu _0}"', from=1-1, to=2-1]
          \arrow [""{name=1, anchor=center, inner sep=0}, "{am_0}", "\shortmid "{marking}, from=1-2, to=1-3]
          \arrow ["{a_0F_{00}\mu _0}"{description}, from=1-2, to=2-2]
          \arrow [""{name=2, anchor=center, inner sep=0}, "{(\overline {f}_0e_{00})^{>}}", "\shortmid "{marking}, from=1-3, to=1-4]
          \arrow ["{a_1F_{01}\mu _1}"{description}, from=1-3, to=2-3]
          \arrow ["{F_{11}f_1\mu _1}", from=1-4, to=2-4]
          \arrow [""{name=3, anchor=center, inner sep=0}, "{(\overline {f}_0e_{00}')^{<}}"', "\shortmid "{marking}, from=2-1, to=2-2]
          \arrow [""{name=4, anchor=center, inner sep=0}, "{am_1}"', "\shortmid "{marking}, from=2-2, to=2-3]
          \arrow [""{name=5, anchor=center, inner sep=0}, "{(\overline {f}_0e_{00})^{>}}"', "\shortmid "{marking}, from=2-3, to=2-4]
          \arrow ["{(\bar  f_{0,\mu _0})^>}"{description}, draw=none, from=0, to=3]
          \arrow ["{a\mu }"{description}, draw=none, from=1, to=4]
          \arrow ["{(\bar  f_{1,\mu _1})^>}"{description}, draw=none, from=2, to=5]
        \end {tikzcd}
      \end{center}

      where the left (resp. right) square is the conjoint (resp. companion) transpose of the tight naturality square associated to \(\overline {f}_0\mu _0\) (resp. \(\overline {f}_1\mu _1\)).
    \item{}\textbf{(Unitors and laxators)} The unit and lax structure is taken from \(f_{i}\) when in \(\mathbb {E}_{0i}\). 
    There are two remaining cases for the laxators: when we are composing a loose in \(\mathbb {E}_{0i}\) with an \((e_{00},e_{01},m_0 : F_{00}e_{00} \mathrel {\mkern 3mu\vcenter {\hbox {$\shortmid $}}\mkern -10mu{\to }} F_{01}e_{01})\) on one side or the other. 
    By symmetry, it is enough to consider the case of \(\ell  : e_{00}' \mathrel {\mkern 3mu\vcenter {\hbox {$\shortmid $}}\mkern -10mu{\to }} e_{00}\). 
    Then we have \(\mathsf {Res}(a)(\ell )=f_0(\ell ),\mathsf {Res}(a)(e_{00},e_{01},m_0)=(f_0(e_{00}),f_1(e_{01}),m_0'),\) and 
    \(\mathsf {Res}(a)(\ell \odot (e_{00},e_{01},m_0))=\mathsf {Res}(a)(e_{00}',e_{01},F_{00}(\ell )\odot  m_0)=(f_0(e_{00}',f_1(e_{01}),[F_{00}(\ell )\odot  m_0]'))\). 
    We define the laxator mapping from the composite of the former two loose maps to the latter to be \((\mathrm {id}(f_0(e_{00}')),\mathrm {id}(f_1(e_{01}),\lambda _{\ell ,m_0}))\) with \(\lambda _{\ell ,m_0}\) the composite below:

  \begin{center}
    \begin {tikzcd}
      {F_{10}f_0e_{00}'} & {F_{10}f_0e_{00}} & {a_0F_{00}e_{00}} & {a_1F_{01}e_{01}} & {F_{11}f_1e_{01}} \\
      {F_{10}f_0e_{00}'} & {a_0F_{00}e'_{00}} & {a_0F_{00}e_{00}} & {a_1F_{01}e_{01}} & {F_{11}f_1e_{01}} \\
      {F_{10}f_0e_{00}'} & {a_0F_{00}e'_{00}} && {a_1F_{01}e_{01}} & {F_{11}f_1e_{01}}
      \arrow ["{F_{10}f_0\ell }", "\shortmid "{marking}, from=1-1, to=1-2]
      \arrow [equals, from=1-1, to=2-1]
      \arrow ["{(\overline {f}_0e_{00})^{<}}", "\shortmid "{marking}, from=1-2, to=1-3]
      \arrow ["{m_0'}"', "\shortmid "{marking}, shift left=2, curve={height=-24pt}, from=1-2, to=1-5]
      \arrow ["{(\overline {f}_0\ell )^{<}}"{description}, draw=none, from=1-2, to=2-2]
      \arrow ["{am_0}", "\shortmid "{marking}, from=1-3, to=1-4]
      \arrow [equals, from=1-3, to=2-3]
      \arrow ["{(\overline {f}_0e_{00})^{>}}", "\shortmid "{marking}, from=1-4, to=1-5]
      \arrow [equals, from=1-4, to=2-4]
      \arrow [equals, from=1-5, to=2-5]
      \arrow ["{(\overline {f}_0e'_{00})^{<}}", "\shortmid "{marking}, from=2-1, to=2-2]
      \arrow [equals, from=2-1, to=3-1]
      \arrow ["{a_0F_{00}\ell }", "\shortmid "{marking}, from=2-2, to=2-3]
      \arrow [equals, from=2-2, to=3-2]
      \arrow ["{am_0}", "\shortmid "{marking}, from=2-3, to=2-4]
      \arrow ["{(\overline {f}_0e_{00})^{>}}", "\shortmid "{marking}, from=2-4, to=2-5]
      \arrow [equals, from=2-4, to=3-4]
      \arrow [equals, from=2-5, to=3-5]
      \arrow ["{(\overline {f}_0e'_{00})^{<}}", "\shortmid "{marking}, from=3-1, to=3-2]
      \arrow ["{[F_{00}\ell \odot  m_0]'}", "\shortmid "{marking}, shift right=3, curve={height=30pt}, from=3-1, to=3-5]
      \arrow [""{name=0, anchor=center, inner sep=0}, "{a(F_{00}\ell  \odot  m_0)}"', "\shortmid "{marking}, from=3-2, to=3-4]
      \arrow ["{(\overline {f}_0e_{00})^{>}}", "\shortmid "{marking}, from=3-4, to=3-5]
      \arrow ["{a_{F_{00}\ell ,m_0}}"{description}, draw=none, from=2-3, to=0]
    \end {tikzcd}
    \end{center}

    Note here that \(a_{\bullet ,\bullet }\) is the laxator of the double functor \(a\). 
    Since \(\overline {f}_0\) is a \hyperref[djm-00EX]{commuter transformation}, the globular cell \((\overline {f}_0\ell )^{<}\) is a tight isomorphism. 
    Furthermore, since \(a\) is a pseudo-functor, this overall cell is a tight isomorphism. 
    This, and its correlate using \(\overline {f}_1\), shows that \(\mathsf {Res}(a)\) will be a pseudo-functor of loose bimodules. 
    If only \(a\) were pseudo, then this double functor would still be lax on account of the square associated to \(\overline {f}_i\); 
    this is why we must assume that the \(\overline {f}_i\) are commuter transformations in \(\mathsf {Res}\).\end{enumerate}\end{proof}Action on 2-cells\label{djm-00J4-proof-two}\par{}Now suppose we have a 2-cell like so:

\begin{center}
\begin {tikzcd}
	{\mathbb {E}_{0i}} & {\mathbb {E}_{1i}} && {\mathbb {E}_{00}} & {\mathbb {E}_{1i}} \\
	{\mathbb {E}_{0i}} & {\mathbb {E}_{1i}} & {=} & {M_{0i}} & {M_{1i}} \\
	{M_{0i}} & {M_{1i}} && {M_{0i}} & {M_{1i}}
	\arrow ["{f_{1i}}", from=1-1, to=1-2]
	\arrow [Rightarrow, no head, from=1-1, to=2-1]
	\arrow [Rightarrow, no head, from=1-2, to=2-2]
	\arrow ["{f_{1i}}", from=1-4, to=1-5]
	\arrow ["{F_{0i}}"', from=1-4, to=2-4]
	\arrow ["{F_{1i}}", from=1-5, to=2-5]
	\arrow ["{\phi_i}", shorten <=5pt, shorten >=5pt, Rightarrow, from=2-1, to=1-2]
	\arrow ["{f_{0i}}", from=2-1, to=2-2]
	\arrow ["{F_{0i}}"', from=2-1, to=3-1]
	\arrow ["{F_{1i}}", from=2-2, to=3-2]
	\arrow ["{\overline {f}_{1i}}", shorten <=5pt, shorten >=5pt, Rightarrow, from=2-4, to=1-5]
	\arrow ["{a_{1i}}", from=2-4, to=2-5]
	\arrow [Rightarrow, no head, from=2-4, to=3-4]
	\arrow [Rightarrow, no head, from=2-5, to=3-5]
	\arrow ["{\overline {f}_{0i}}", shorten <=5pt, shorten >=5pt, Rightarrow, from=3-1, to=2-2]
	\arrow ["{ a_{0i}}"', from=3-1, to=3-2]
	\arrow ["{\alpha_i}", shorten <=6pt, shorten >=6pt, Rightarrow, from=3-4, to=2-5]
	\arrow ["{a_{0i}}"', from=3-4, to=3-5]
\end {tikzcd}
		\end{center}

		We then define a tight transformation \(\mathsf {Res}(\alpha , \cdots )\) as follows:
		\begin{enumerate}\item{}				On \(\mathbb {E}_{0i}\), \(\mathsf {Res}(\alpha , \cdots )\) is just \(\phi _i\).
			\item{}				The only nontrivial loose case to handle is \(m_0 : F_{00}e_{00} \mathrel {\mkern 3mu\vcenter {\hbox {$\shortmid $}}\mkern -10mu{\to }} F_{01}e_{01}\), which we assign to the square
				
\begin{center}
				\begin {tikzcd}
	{F_{10}f_{00}e_{00}} & {(a_{00})F_{00}e_{00}} & {(a_{01})F_{01}e_{01}} & {F_{11}f_{01}e_{01}} \\
	{F_{10}f_{10}e_{00}} & {(a_{10})F_{00}e_{00}} & {(a_{11})F_{01}e_{01}} & {F_{11}f_{11}e_{01}}
	\arrow ["{(\overline {f}_{00}e_{00})^{<}}", "\shortmid "{marking}, from=1-1, to=1-2]
	\arrow ["{F_{10}\phi _0}"', from=1-1, to=2-1]
	\arrow [""{name=0, anchor=center, inner sep=0}, "{a_0m_0}", "\shortmid "{marking}, from=1-2, to=1-3]
	\arrow ["{\alpha_0  F_{00}}"', from=1-2, to=2-2]
	\arrow ["{(\overline {f}_{01}e_{00})^{>}}", "\shortmid "{marking}, from=1-3, to=1-4]
	\arrow ["{\alpha_1  F_{01}}", from=1-3, to=2-3]
	\arrow ["{F_{11}\phi _1}", from=1-4, to=2-4]
	\arrow ["{(\overline {f}_{10}e_{00})^{<}}"', "\shortmid "{marking}, from=2-1, to=2-2]
	\arrow [""{name=1, anchor=center, inner sep=0}, "{a_1m_0}"', "\shortmid "{marking}, from=2-2, to=2-3]
	\arrow ["{(\overline {f}_{11}e_{00})^{>}}"', "\shortmid "{marking}, from=2-3, to=2-4]
	\arrow ["{\alpha  m}"{description}, shorten <=4pt, shorten >=4pt, Rightarrow, from=0, to=1]
\end {tikzcd}
			\end{center}

				where the left and right squares are the transposes of the commutativity equation defining a 2-cell in \(\mathcal {N}\mathsf {iche}\).\end{enumerate}\end{proof}\end{theorem}

\begin{proposition}[{Source and target of restricted loose bimodules}]\label{ssl-004G}
\par{}The following diagram commutes.
  
\begin{center}
    \begin {tikzcd}
      {\mathcal {N}\mathsf {iche}} & {\ell }\mathcal {B}\mathsf {imod} \\
      {\mathsf {2}\mathcal {C}\mathsf {at}^{\mathsf {colax}}(\Delta [1], \mathcal {D}\mathsf {bl})_{\mathsf {conj}} \times  \mathsf {2}\mathcal {C}\mathsf {at}^{\mathsf {colax}}(\Delta [1], \mathcal {D}\mathsf {bl})_{\mathsf {comp}}} & {\mathcal {D}\mathsf {bl} \times  \mathcal {D}\mathsf {bl}}
      \arrow ["{\mathsf {Res}}", from=1-1, to=1-2]
      \arrow [from=1-1, to=2-1]
      \arrow ["{(()_0, ()_1)}", from=1-2, to=2-2]
      \arrow ["{d_1^* \times  d_1^*}"', from=2-1, to=2-2]
    \end {tikzcd}
  \end{center}

  On objects this means that given a niche, 
  
\begin{center}
        \begin {tikzcd}
          {\mathbb {E}_0} & {\mathbb {E}_1} \\
          {M_0} & {M_1}
          \arrow ["{F_0}"', from=1-1, to=2-1]
          \arrow ["{F_1}",from=1-2, to=2-2]
          \arrow ["M"', "{\shortmid }"{marking}, from=2-1, to=2-2]
        \end {tikzcd}
    \end{center}

    its restriction is a loose bimodule whose source is \(\mathbb {E}_0\) and whose target is \(\mathbb {E}_1\).\begin{proof}\par{}This follows by definition of the 2-functor \(\mathsf {Res}\) given in \hyperref[djm-00J4]{Theorem \ref{djm-00J4}}.\end{proof}\end{proposition}
\par{}	As a corollary of \hyperref[djm-00J4]{Theorem \ref{djm-00J4}}, we conclude that restriction of a loose bimodule \(\mathbb {M} \colon \mathbb {D}_0 \mathrel {\mkern 3mu\vcenter {\hbox {$\shortmid $}}\mkern -10mu{\to }} \mathbb {D}_1\) along lax symmetric monoidal pseudo-double functors \(F_i : \mathbb {E}_i \to  \mathbb {D}_i\) induces a symmetric monoidal structure on the restriction \(\mathbb {M}(F_0, F_1)\) so long as the unitors and laxitors of \(F_0\) (resp. \(F_1\)) are \hyperref[djm-00EX]{conjoint (resp. companion) commuter transformations}. For the remainder of this section, we will expand on this observation.\paragraph{Cartesian and symmetric monoidal structure in \(\mathcal {N}\mathsf {iche}\)}\label{ssl-003M}\par{}Next we will investigate the cartesian and symmetric monoidal structure in \(\mathcal {N}\mathsf {iche}\).
\begin{lemma}[{\(\mathcal {N}\mathsf {iche}\), and \(\mathcal {N}\mathsf {iche}^{\cong }\) are cartesian 2-categories}]\label{djm-00BU}
\par{}The \hyperref[djm-00BS]{domain \(\mathcal {N}\mathsf {iche}\) (resp. \(\mathcal {N}\mathsf {iche}^{\cong }\)) of restriction of loose bimodules} is a cartesian 2-category, and its cartesian products may be constructed componentwise in \({\ell }\mathcal {B}\mathsf {imod}\) and \(\mathcal {D}\mathsf {bl}\).\par{}In particular, the terminal object is the niche \begin{center}
\begin{tikzcd}
	0 & 1 \\
	{\mathbb {L}\mathsf {oose}} & {\mathbb {L}\mathsf {oose}}
	\arrow[hook, from=1-1, to=2-1]
	\arrow[hook, from=1-2, to=2-2]
	\arrow["{\mathsf {Hom}^{l}(\mathbb {L}\mathsf {oose})}"', "\shortmid"{marking}, from=2-1, to=2-2]
\end{tikzcd}
\end{center}

including the endpoints of the \hyperref[djm-0047]{walking loose arrow}.\begin{proof}\par{}As a pullback of cartesian functors between cartesian 2-categories, \(\mathcal {N}\mathsf {iche}\) is cartesian, and its products are constructed componentwise. This relies on knowing that \(\mathsf {2}\mathcal {C}\mathsf {at}^{\mathsf {colax}}(\Delta [1], \mathcal {D}\mathsf {bl})_{\mathsf {comp}}\) and \(\mathsf {2}\mathcal {C}\mathsf {at}^{\mathsf {colax}}(\Delta [1], \mathcal {D}\mathsf {bl})_{\mathsf {conj}}\) are cartesian, but they are since \(\mathsf {2}\mathcal {C}\mathsf {at}^{\mathsf {colax}}(\Delta [1], \mathcal {D}\mathsf {bl})\) is cartesian and a tight transformation into a product is a commuter if and only if its components are.\par{}To see that \(\mathcal {N}\mathsf {iche}^{\cong }\) is also cartesian, it suffices to note that the double functor induced into a product or a pullback is pseudo when both of its components are.\end{proof}\end{lemma}
\par{}We can use \emph{symmetry of internalization} to characterize the symmetric monoidal and cartesian structure in \(\mathcal {N}\mathsf {iche}\).
\begin{lemma}[{Symmetric monoidal objects of \(\mathcal {N}\mathsf {iche}\) and \(\mathcal {N}\mathsf {iche}^{\cong }\)}]\label{djm-00CM}
\par{}A symmetric monoidal object of \(\mathcal {N}\mathsf {iche}^{\cong }\) (resp. \(\mathcal {N}\mathsf {iche}\)) is equivalently:
\begin{enumerate}\item{}		A \hyperref[ssl-0019]{symmetric monoidal loose bimodule} \(M : \mathbb {M} \to  \mathbb {L}\mathsf {oose}\).
	\item{}		Symmetric monoidal double categories \(\mathbb {E}_0\) and \(\mathbb {E}_1\).
	\item{}		Pseudo (resp. lax) symmetric monoidal pseudo-functors \(F_i : \mathbb {E}_i \to   M_i\), where the unitors and laxators of \(F_0\) are , and the unitors and laxators of \(F_1\) are .\end{enumerate}\begin{proof}\par{}First, we note that the projection 2-functors \(M : \mathcal {N}\mathsf {iche} \to  {\ell }\mathcal {B}\mathsf {imod}\) and \(\mathbb {E}_i : \mathcal {N}\mathsf {iche} \to  \mathcal {D}\mathsf {bl}_u\) are cartesian, and therefore they preserve symmetric monoidal structure. It only remains to show that the colax squares \((\otimes _i, \overline {\otimes }_i) : F_i \times  F_i \to  F_i\) and \((1_i, \overline {1}_i) : \{i\} \to  F_i\) assemble into a pseudo (resp. lax) monoidal morphism.

\begin{center}
\begin {tikzcd}
	{\mathbb {E}_i \times  \mathbb {E}_i} & {\mathbb {E}_i} & {\{i\}} & {\mathbb {E}_i} \\
	{ M_i \times M_i} & { M_i} & {\{i\}} & {M_i}
	\arrow ["{\otimes _i}", from=1-1, to=1-2]
	\arrow ["{F_i \times  F_i}"', from=1-1, to=2-1]
	\arrow ["{F_i}", from=1-2, to=2-2]
	\arrow ["{1_i}", from=1-3, to=1-4]
	\arrow [from=1-3, to=2-3]
	\arrow ["{F_i}", from=1-4, to=2-4]
	\arrow ["{\overline {\otimes }_i}", shorten <=6pt, shorten >=6pt, Rightarrow, from=2-1, to=1-2]
	\arrow ["{\otimes_i }"', from=2-1, to=2-2]
	\arrow ["{\overline {1}_i}", shorten <=4pt, shorten >=4pt, Rightarrow, from=2-3, to=1-4]
	\arrow ["{1_i}"', from=2-3, to=2-4]
\end {tikzcd}
\end{center}

 This follows by the \emph{symmetry of internalization} (Theorem 7.4) of \cite{arkor-2024-enhanced}; specifically, a symmetric pseudo-monoid in maps with pseudo (resp. colax) squares is a pseudo (resp. lax) monoidal morphism.\end{proof}\end{lemma}

\begin{lemma}[{Cartesian objects of \(\mathcal {N}\mathsf {iche}^{\cong }\)}]\label{djm-00DD}
\par{}	A cartesian object of \(\mathcal {N}\mathsf {iche}^{\cong }\) consists of 
	\begin{enumerate}\item{}A cartesian loose bimodule \(M\).
		\item{}Cartesian double categories \(\mathbb {E}_0\) and \(\mathbb {E}_1\).
		\item{}Cartesian pseudo functors \(F_i : \mathbb {E}_i \to  M_i\), such that the laxators and unitors of \(F_0\) are conjoints and of \(F_1\) are companions.\end{enumerate}\begin{proof}\par{}	We may again compute the cartesian objects componentwise in \({\ell }\mathcal {B}\mathsf {imod}\) and \(\mathsf {2}\mathcal {C}\mathsf {at}^{\mathsf {ps}}(\Delta [1], \mathcal {D}\mathsf {bl})^2\). By definition, a cartesian loose bimodule is a cartesian object in \({\ell }\mathcal {B}\mathsf {imod}\). On the other hand, by \emph{symmetry of internalization} (Theorem 7.4 of \cite{arkor-2024-enhanced}), a cartesian object in \(\mathsf {2}\mathcal {C}\mathsf {at}^{\mathsf {ps}}(\Delta [1], \mathcal {D}\mathsf {bl})\) is equivalently a cartesian functor between cartesian double categories.\end{proof}\end{lemma}
\paragraph{Restriction preserves cartesian and symmetric monoidal structure}\label{ssl-003N}\par{}Since restriction of loose bimodules is cartesian, it will preserve some 2-algebraic structure. However, since it is only a pseudo-functor, it won't preserve all 2-algebraic structure. Specifically, \emph{flexible} 2-algebraic structure is preserved by cartesian pseudo-functors. Intuitively, 2-algebraic structure is \emph{flexible} when it does not require any equations between 1-cells; this includes symmetric monoidal and cartesian structure. If 2-algebraic structure is flexible, then we can pass it through a pseudo-functor by conjugating every 2-cell by the unitors and laxators of the pseudo-functor.\par{}As a corollary, \hyperref[djm-00J4]{restriction of looes bimodules \(\mathsf {Res}\)} will preserve symmetric monoidal and cartesian structure.\par{}Recall that for a cartesian 2-category \(\mathcal {K}\), 
  \begin{itemize}\item{}\(\mathcal {S}\mathsf {M}^{\mathsf {lax}}{(\mathcal {K})}\) is the 2-category of symmetric monoidal objects in and lax monoidal morphisms.
    \item{}\(\mathcal {C}\mathsf {art}(\mathcal {K})\) is the 2-category of cartesian objects, cartesian maps, and general 2-cells.\end{itemize}
\begin{theorem}[{Restriction of loose bimodules preserves symmetric monoidal and cartesian structure}]\label{djm-00CL}
\par{} extends to pseudo-functors 
\[\mathcal {S}\mathsf {M}(\mathcal {N}\mathsf {iche}) \to  \mathcal {S}\mathsf {M}({\ell }\mathcal {B}\mathsf {imod})\]
and 
\[\mathcal {C}\mathsf {art}(\mathcal {N}\mathsf {iche}) \to  \mathcal {C}\mathsf {art}({\ell }\mathcal {B}\mathsf {imod}).\]
Similarly, for \(\mathsf {Res}^{\cong } : \mathcal {N}\mathsf {iche}^{\cong } \to  {\ell }\mathcal {B}\mathsf {imod}\).\end{theorem}
\par{}	Applying \hyperref[djm-00CL]{this preservation of symmetric monoidal and cartesian structure} to the characterizations of  \hyperref[djm-00CM]{symmetric monoidal} and \hyperref[djm-00DD]{cartesian} of \hyperref[djm-00BS]{\(\mathcal {N}\mathsf {iche}\)} leads to the following conclusions: 
	\begin{itemize}\item{}			Restriction of a symmetric monoidal loose bimodule by pseudo symmetric monoidal pseudo-functors (whose unitors and laxators are conjoints or companions, respectively) yields a symmetric monoidal loose bimodule.
		\item{}			Similarly, restriction of a symmetric monoidal loose bimodule by \emph{lax} symmetric monoidal pseudo-functors (whose unitors and laxators are conjoint or companion \hyperref[djm-00EX]{commuter transformations}, respectively) will yield a symmetric monoidal loose bimodule.
		\item{}			Finally, restriction of a cartesian bimodule by cartesian pseudofunctors yields a cartesian loose bimodule.\end{itemize}
\begin{explication}[{Restriction of a symmetric monoidal loose bimodule}]\label{ssl-003P}
\par{}Let 
  
\begin{center}
    \begin {tikzcd}
      {\mathbb {E}_0} & {\mathbb {E}_1} \\
      {M_0} & {M_1}
      \arrow ["{F_0}"', from=1-1, to=2-1]
      \arrow ["{F_1}",from=1-2, to=2-2]
      \arrow ["M"', "\shortmid "{marking}, from=2-1, to=2-2]
    \end {tikzcd}
  \end{center}

  be a \hyperref[djm-00CM]{symmetric monoidal niche}.\par{}By \hyperref[djm-00CL]{Theorem \ref{djm-00CL}}, the restricted loose bimodule \(\mathbb {M}(F_0, F_1)\) is a \hyperref[ssl-0019]{symmetric monoidal loose bimodule}. Here we explicate the monoidal structure for loose heteromorphisms in \(\mathbb {M}(F_0, F_1)\) and show how it uses the conjoints of the laxator of \(F_0\) and companions of the laxator of \(F_1\).\par{}The laxator of \(F_0\) is depicted in the following diagram.
  
\begin{center}
    \begin {tikzcd}
      {\mathbb {E}_0 \times  \mathbb {E}_0} & {\mathbb {E}_0} \\
      {M_0 \times  M_0} & { M_0}
      \arrow ["{\otimes _0}", from=1-1, to=1-2]
      \arrow ["{F_0 \times  F_0}"', from=1-1, to=2-1]
      \arrow ["{F_0}", from=1-2, to=2-2]
      \arrow ["{\overline {\otimes }_0}", shorten <=7pt, shorten >=7pt, Rightarrow, from=2-1, to=1-2]
      \arrow ["{\otimes_0 }"', from=2-1, to=2-2]
    \end {tikzcd}
  \end{center}

  It is a . So for the identity loose morphism on \(F_0e_0 \otimes  F_0e_0'\), we have the conjoint \(\overline {\otimes }_0^< \colon  F(e_0 \otimes _0 e_0') \to  Fe_0 \otimes  Fe_0'\).\par{}Likewise let \(\overline {\otimes }_1\) be the laxator of \(F_1\). It is a companion commuter transformation. In particular, the companion of the laxator \(\overline {\otimes }_1\) applied to the identity loose morphism on \(F_1 e_1 \otimes  F_1 e_1'\) induces a loose morphism \(\overline {\phi }_1 ^> \colon  F_1 e_1 \otimes  F_1 e_1' \to  F_1(e_1 \otimes _1 e_1')\).\par{}	Now let \(m \colon  F_0 e_{0} \mathrel {\mkern 3mu\vcenter {\hbox {$\shortmid $}}\mkern -10mu{\to }} F_1 e_{1}\) and \(m' \colon  F_0 e_{0}' \mathrel {\mkern 3mu\vcenter {\hbox {$\shortmid $}}\mkern -10mu{\to }} F_1 e_{1}'\) be . Then their monoidal product in the restriction \(\mathbb {M}(F_0, F_1)\) is defined to be the composite
	
\begin{center}
		\begin {tikzcd}
      {F_0(e_{0} \otimes _0 e_{0}')} & {F_0 e_{0} \otimes _0 F_0 e_{0}'} & {F_1e_{1} \otimes _1 F_1e_{1}'} & {F_1(e_{1} \otimes _1 e_{1}')}
      \arrow ["{\overline {\otimes _0}^{<}}", "\shortmid "{marking}, from=1-1, to=1-2]
      \arrow ["{m \otimes  m'}", "\shortmid "{marking}, from=1-2, to=1-3]
      \arrow ["{\overline {\otimes _1}^{>}}", "\shortmid "{marking}, from=1-3, to=1-4]
    \end {tikzcd}
	\end{center}\par{}This formula works in all the cases above; but if the \(F_i\) are pseudo-symmetric monoidal, then the \(\overline {\otimes }_i\) will be isomorphisms, and if everything involved is cartesian, then these will furthermore be cartesian products.\end{explication}

\begin{remark}[{Remark on the preservation of cartesian structure by certain colax 2-functors}]\label{djm-00F2}
\par{}	As remarked in the , we could actually end up with a \emph{normal colax} 2-functor \(\mathsf {Res}^{\mathsf {lax}} : \mathcal {N}\mathsf {iche}^{\mathsf {lax}} \to  {\ell }\mathcal {B}\mathsf {imod}^{\mathsf {lax}}\) if we loosened up the definition of \(\mathcal {N}\mathsf {iche}\). This normal colax 2-functor would continue to preserve products, and moreover it would have the special property that the colaxator associated to composition with a product projection, from either side, is an isomorphism.\par{}	This suggests a general consideration concering colax \(\mathcal {F}\)-functors between \(\mathcal {F}\)-categories. It appears that the correct notion of such a functor is to be a colax 2-functor, but where the counitor is an isomorphism (it is normal) and the colaxator associated to composition with an inert morphism, on either side, is an isomorphism. We may see normality as a nullary case of the latter condition.\par{}	A bit of fiddling has made it seem reasonable to expect that such colax \(\mathcal {F}\)-functors which preserve products would also preserve cartesian objects (this can be straightforwardly verified) and furthermore induce a push-forward on the \(\mathcal {F}\)-categories of cartesian objects (this would benefit from a theoretical explanation). We leave it for an intrepid 2-algebraist to continue these musings.\par{}	We finish this remark by noting that a theorem of the above sort would give us a slight advantage in the double operadic theory of systems: we would know that the restriction of pre-cartesian loose bimodules by pre-cartesian pseudo-functors is still pre-cartesian. This is only a slight improvement on the results of this section; nevertheless, the general question is interesting on its own.\end{remark}

\begin{explication}[{Summary of results of \hyperref[djm-00J3]{Section \ref{djm-00J3}}}]\label{djm-00DM}
\par{}	We may put together the calculations \hyperref[djm-00CM]{Lemma \ref{djm-00CM}} and \hyperref[djm-00DD]{Lemma \ref{djm-00DD}} to conclude the following:
\begin{enumerate}\item{}		Restriction of a symmetric monoidal loose bimodule by pseudo symmetric monoidal pseudo-functors (whose unitors and laxators are conjoints or companions, respectively) yields a symmetric monoidal loose bimodule.
	\item{}		Similarly, restriction of a symmetric monoidal loose bimodule by \emph{lax} symmetric monoidal pseudo-functors (whose unitors and laxators are conjoint or companion \hyperref[djm-00EX]{commuter transformations}, respectively) will yield a symmetric monoidal loose bimodule.
	\item{}		Finally, restriction of a cartesian bimodule by cartesian pseudofunctors yields a cartesian loose bimodule.\end{enumerate}
We will use these constructions liberally to produce examples of systems theories as symmetric monoidal loose right modules.\par{}	Explicitly, if \(m_0 \in  \mathbb {M}(F_0 e_{00}, F_1 e_{10})\) and \(m_1 \in  \mathbb {M}(F_0 e_{01}, F_1 e_{11})\) are loose heteromorphisms of the restricted loose bimodule, then their monoidal product \(m_0 \otimes _{\mathsf {res}} m_1\) in the restriction is defined to be 
	
\begin{center}
		\begin {tikzcd}
	{F_0(e_{00} \otimes _0 e_{01})} & {F_0 e_{00} \otimes _0 F_0 e_{01}} & {F_1e_{10} \otimes _1 F_1e_{11}} & {F_1(e_{10} \otimes _1 e_{11})}
	\arrow ["{\overline {\otimes _0}^{<}}", "\shortmid "{marking}, from=1-1, to=1-2]
	\arrow ["{m_0 \otimes  m_1}", "\shortmid "{marking}, from=1-2, to=1-3]
	\arrow ["{\overline {\otimes _1}^{>}}", "\shortmid "{marking}, from=1-3, to=1-4]
\end {tikzcd}
	\end{center}

	where \(\overline {\otimes }_i\) are the laxators of \(\otimes _i : \mathbb {E}_i \times  \mathbb {E}_i \to  \mathbb {E}_i\), and we've taken the conjoint and companion respectively. This formula works in all the cases above; but if the \(F_i\) are pseudo-symmetric monoidal, then the \(\overline {\otimes }_i\) will be isomorphisms, and if everything involved is cartesian, then these will furthermore be cartesian products.\end{explication}
\section{Symmetric monoidal loose right modules as modules of systems}\label{ssl-0013}\par{}In this section, we define symmetric monoidal loose right modules and describe how they appropriately organize a collection of systems, their interactions, and their maps.  In other words, we will show how the answers to the  questions posed in \hyperref[djm-009K]{Informal definition \ref{djm-009K}} form a symmetric monoidal loose right module. To emphasize this attitude we will often refer to a symmetric monoidal loose right module as a \emph{module of systems over a double category of interactions}.  We also give several examples of this attitude that we will make formal in \hyperref[ssl-0043]{Section \ref{ssl-0043}}.\subsection{Module of systems over interactions}\label{ssl-0047}\par{}We are now ready to give the main attitude of this section: that of a \emph{module of systems over a double category of interactions}. We begin by definining symmetric monoidal loose right modules as symmetric monoidal objects of the \hyperref[djm-004M]{2-category of loose right modules}.
\begin{definition}[{The 2-category of symmetric monoidal loose right  modules}]\label{ssl-003U}
\par{}The 2-category \(\mathcal {S}\mathsf {M}({\ell }\mathcal {M}\mathsf {od}_{\mathsf {r}})\) is the 2-category of symmetric pseudo-monoids and pseudo-morphisms in the cartesian 2-category \({\ell }\mathcal {M}\mathsf {od}_{\mathsf {r}}\) of \hyperref[djm-004M]{loose right modules}. We call the objects of this 2-category \textbf{symmetric monoidal loose right modules}.\end{definition}
\par{}It is helpful to view a symmetric monoidal loose right module as a symmetric monoidal category with an action.
\begin{explication}[{Symmetric monoidal loose right module}]\label{ssl-0044}
\par{}We can think of a symmetric monoidal loose right module \(\mathbb {M} \colon  \bullet  \mathrel {\mkern 3mu\vcenter {\hbox {$\shortmid $}}\mkern -10mu{\to }}  \mathbb {D}_1\) as consisting of:
  \begin{itemize}\item{}Its target, a symmetric monoidal double category \(\mathbb {D}_1\).
    \item{}Its carrier, a symmetric monoidal category \(\mathsf {Car}(\mathbb {M})\).
    \item{}A functor from the carrier to the tight category of the target, \(\mathsf {Car}(\mathbb {M}) \to   \mathsf {Tight}({\mathbb {D}_1})\).
    \item{}An action of the loose category of the target on the carrier \(\mathsf {Loose}({\mathbb {D}_1})\curvearrowright   \mathsf {Car}(\mathbb {M})\).\end{itemize}
  By symmetry of internalization, a symmetric monoidal loose right module is a right pseudo-module internal to symmetric monoidal categories.\end{explication}
\par{}We now interpret a \hyperref[ssl-003U]{symmetric monoidal loose right module} as an \emph{module of systems} over a symmetric monoidal double category of \emph{interfaces and interactions}.
\begin{attitude}[{Module of systems}]\label{ssl-001J}
\par{}We will often consider a \hyperref[ssl-003U]{symmetric monoidal loose right module} to organize the data of systems, their interactions, and their maps. When we do so, we will refer to this module as a \textbf{module of systems}.\end{attitude}

\begin{notation}[{Module of systems}]\label{ssl-003T}
\par{}Given a \hyperref[ssl-001J]{module of systems} \(\mathbb {S} \colon  \bullet  \mathrel {\mkern 3mu\vcenter {\hbox {$\shortmid $}}\mkern -10mu{\to }} \mathbb {I}\), we say that:
  \begin{itemize}\item{}The \hyperref[djm-004A]{target} \(\mathbb {I}\) is the symmetric monoidal double category of \textbf{interfaces} and \textbf{interactions} mediating them.
    \item{}The \hyperref[ssl-0036]{carrier} \(\mathsf {Car}(\mathbb {S})\) is the symmetric monoidal category of \textbf{systems} and \textbf{system maps}.\end{itemize}\par{}Following the phrase a "module over a ring", we often call \(\mathbb {S}\) a \textbf{module of systems over interactions, \(\mathbb {I}\)}.\end{notation}
\par{}To further make sense of how a module of systems organizes the data of systems, their interactions, and their maps, it is helpful to rename the objects, morphisms, and squares of \(\mathbb {L}\mathsf {oose}\) so that the names give attitudes to each objects, morphisms, and squares that they label. These attitudes were introduced in \hyperref[djm-009N]{ontology of double categorical systems theory}.
\begin{attitude}[{Reinterpretting the walking loose arrow for modules of systems}]\label{ssl-001H}
\par{}Let  \(\mathbb {S} \to  \mathbb {L}\mathsf {oose}\) be a \hyperref[ssl-001J]{module of systems}. This functor labels the objects, morphisms, and squares of \(\mathbb {S}\) by the objects, morphisms, and squares of \(\mathbb {L}\mathsf {oose}\).  When considering \(\mathbb {S}\)  as a \hyperref[ssl-001J]{module of systems}, it is helpful to rename the objects, morphisms, and squares of \(\mathbb {L}\mathsf {oose}\) so that they assign meaningful labels to the different components of the module of systems, \(\mathbb {S}\).\par{}First, recall the components of \(\mathbb {L}\mathsf {oose}\) from \hyperref[ssl-001G]{Explication \ref{ssl-001G}}.
\begin{explication}[{}]\label{djm-00JE}
\par{}The \hyperref[djm-0047]{walking loose arrow double category} \(\mathbb {L}\mathsf {oose}\) has two objects, three loose arrows (two identities along with the walking loose arrow \(2_{\ell }\)), two identity tight arrows, and three identity squares. These are all pictured below.

  \begin{center}\begin {tikzcd}
    0 & 0 & 1 & 1 \\
    0 & 0 & 1 & 1
    \arrow ["\shortmid "{marking}, equals, from=1-1, to=1-2]
    \arrow [equals, from=1-1, to=2-1]
    \arrow ["\shortmid "{marking}, from=1-2, to=1-3]
    \arrow [equals, from=1-2, to=2-2]
    \arrow ["\shortmid "{marking}, equals, from=1-3, to=1-4]
    \arrow [equals, from=1-3, to=2-3]
    \arrow [equals, from=1-4, to=2-4]
    \arrow ["\shortmid "{marking}, equals, from=2-1, to=2-2]
    \arrow ["\shortmid "{marking}, from=2-2, to=2-3]
    \arrow ["\shortmid "{marking}, equals, from=2-3, to=2-4]
  \end {tikzcd}\end{center}\end{explication}
\par{}Next we strategically name each component. 

  \begin{itemize}\item{}\textbf{The source} of a module of systems is the symmetric monoidal double category living over the left object \(0\) and its identity morphisms and square.  Since the source of a systems theory is the terminal double category, we will rename the left object \(\bullet \).
    \item{}\textbf{The target} of a module of systems is the symmetric monoidal  double category living over the right object \(1\) and its identity morphisms and square. We adopt the attitude that it represents \emph{interfaces} and \emph{interactions} that relate them. Therefore, we name the right object \(\mathsf {interface}\), its tight identity \(\mathsf {interface\, map}\), its loose identity \(\mathsf {interaction}\), and its identity square \(\mathsf {interaction\, map}\).
    \item{}\textbf{The carrier} of a module of systems is the symmetric monoidal category living over the walking arrow and its identity square.  We adopt the attitude that the loose morphisms in module of systems labeled by the walking arrow will represent \emph{systems} while the squares labeled by the identity square for the walking arrow represent \emph{system maps}. Therefore, we name the walking loose arrow \(\mathsf {system}\) and its identity square \(\mathsf {system\,map}\).\end{itemize}\par{}Collectively, these renamings result in the following new presentation of the walking loose arrow \(\mathbb {L}\mathsf {oose}\). Note that in this presentation we no longer show the identity morphisms and square for the left object \(\bullet \), because in a systems theory they do not label any non-identity morphisms and squares.\begin{figure}[H]\begin{center}
\begin {tikzcd}[column sep = huge]
	\bullet  & {\mathsf {interface}} & {\mathsf {interface}} \\
	\bullet  & {\mathsf {interface}} & {\mathsf {interface}}
	\arrow [""{name=0, anchor=center, inner sep=0}, "{\mathsf {system}}", "\shortmid "{marking}, from=1-1, to=1-2]
	\arrow [Rightarrow, no head, from=1-1, to=2-1]
	\arrow [""{name=1, anchor=center, inner sep=0}, "{\mathsf {interaction}}", "\shortmid "{marking}, from=1-2, to=1-3]
	\arrow [from=1-2, to=2-2]
	\arrow ["{\mathsf {interface\, map}}", from=1-3, to=2-3]
	\arrow [""{name=2, anchor=center, inner sep=0}, "{\mathsf {system}}"', "\shortmid "{marking}, from=2-1, to=2-2]
	\arrow [""{name=3, anchor=center, inner sep=0}, "{\mathsf {interaction}}"', "\shortmid "{marking}, from=2-2, to=2-3]
	\arrow ["{\mathsf {system\, map}}"{description}, draw=none, from=0, to=2]
	\arrow ["{\mathsf {interaction\, map}}"{description}, draw=none, from=1, to=3]
\end {tikzcd}
\end{center}\caption{Ontology of double categorical systems theory}\end{figure}\end{attitude}
\par{}The labels and their interactions are a scheme for interpretting a module of systems.
\begin{explication}[{Labelling perspective for modules of systems}]\label{ssl-001I}
\par{}Given a module of systems \(\mathbb {S} \to  \mathbb {L}\mathsf {oose}\), we refer to the objects living over \(\mathsf {interface}\) as \emph{interfaces}, the loose morphisms living over \(\mathsf {system}\) as \emph{systems}, the loose morphisms living over \(\mathsf {interaction}\) as \emph{interactions}, and so forth for the tight morphisms and squares.\par{}To elaborate the notation in \hyperref[ssl-003T]{Notation \ref{ssl-003T}}, this means that given a module of systems \(\mathbb {S} \colon  \bullet  \mathrel {\mkern 3mu\vcenter {\hbox {$\shortmid $}}\mkern -10mu{\to }} \mathbb {I}\). 
  \begin{itemize}\item{}The target category \(\mathbb {I}\) is the symmetric monoidal double category whose:
      \begin{itemize}\item{}Objects are interfaces.
        \item{}Tight morphisms are interface maps.
        \item{}Loose morphisms are interactions.
        \item{}Squares are interaction maps.\end{itemize}
    \item{}The carrier \(\mathsf {Car}(\mathbb {S})\) is the symmetric monoidal category whose:
      \begin{itemize}\item{}Objects are systems.
        \item{}Morphisms are system maps.\end{itemize}\end{itemize}\par{}We can derive features of a systems theory from the structure of \(\mathbb {L}\mathsf {oose}\). For example, 
  \begin{itemize}\item{}Every system has an interface (its codomain).
    \item{}The composite of a system and an interaction is again a system. This composition is what we mean by \emph{interactions acting on systems}.\end{itemize}\begin{center}
    \begin {tikzcd}
      \bullet  & {\mathsf {interface}} & {\mathsf {interface}}
      \arrow ["{\mathsf {system}}", "\shortmid "{marking}, from=1-1, to=1-2]
      \arrow ["{\mathsf {system}}"', "\shortmid "{marking}, curve={height=24pt}, from=1-1, to=1-3]
      \arrow ["{\mathsf {interaction}}", "\shortmid "{marking}, from=1-2, to=1-3]
    \end {tikzcd}
  \end{center}\par{}Finally, the symmetric monoidal structure  implies that systems, interactions, and their maps can be composed in parallel. This feature is critical for the operadic perspective in which an interaction may act on several independent, component systems, whose result is a single composite system.\end{explication}
\par{}We can now re-interpret \hyperref[ssl-0044]{the explication of a symmetric monoidal loose right module} from the perspective of a module of systems.
\begin{explication}[{Module of systems}]\label{ssl-0045}
\par{}A module of systems \(\mathbb {S} \colon  \bullet  \mathrel {\mkern 3mu\vcenter {\hbox {$\shortmid $}}\mkern -10mu{\to }} \mathbb {I}\) consists of:
  \begin{itemize}\item{}A symmetric monoidal double category of interactions \(\mathbb {I}\).
    \item{}A symmetric monoidal category of systems \(\mathsf {Car}(\mathsf {S})\).
    \item{}A functor assigning systems to their interface, \(\mathsf {Car}(\mathsf {S}) \to  \mathsf {Tight}({\mathbb {I}})\).
    \item{}An action of interactions on systems \(\mathsf {Loose}({\mathbb {I}}) \curvearrowright   \mathsf {Car}(\mathsf {S})\) that respects the interface.\end{itemize}\end{explication}

\begin{remark}[{Restriction of modules of systems}]\label{djm-00J5}
\par{}	If \(\mathbb {S} : \bullet  \mathrel {\mkern 3mu\vcenter {\hbox {$\shortmid $}}\mkern -10mu{\to }} \mathbb {I} \) is a \hyperref[ssl-001J]{module of systems} and \(r : \mathbb {J} \to  \mathbb {I}\) is a lax symmetric monoidal double functor (with unitors and laxators \(\mu \) \hyperref[djm-00EX]{companion commuters}) which we think of as giving a restriction of the double category of interactions from \(\mathbb {I}\) to \(\mathbb {J}\), then for any two systems \(S_1 : \bullet  \mathrel {\mkern 3mu\vcenter {\hbox {$\shortmid $}}\mkern -10mu{\to }} rJ_1\) and \(S_2 : \bullet  \mathrel {\mkern 3mu\vcenter {\hbox {$\shortmid $}}\mkern -10mu{\to }} rJ_2\) are systems, their parallel product as systems over the interface double category \(\mathbb {J}\) is the loose composite:
	
\begin{center}
		\begin {tikzcd}
	\bullet  & {rJ_1 \otimes  rJ_2} & {r(J_1 \otimes  J_2)}
	\arrow ["{S_1 \otimes  S_2}", "\shortmid "{marking}, from=1-1, to=1-2]
	\arrow ["{S_1 \otimes _r S_2}"', "\shortmid "{marking}, bend right, from=1-1, to=1-3]
	\arrow ["{\mu ^{>}}", "\shortmid "{marking}, from=1-2, to=1-3]
\end {tikzcd}
	\end{center}\end{remark}
\subsection{Examples of modules of systems}\label{ssl-001F}\par{}In this section, we will give illustrative examples of modules of systems that are constructed via the doctrines that will be presented in \hyperref[ssl-0041]{Section \ref{ssl-0041}} and \hyperref[djm-00G6]{Section \ref{djm-00G6}}. We give their full definition in \hyperref[ssl-0043]{Section \ref{ssl-0043}} but introduce them here to exemplify how modules of systems answer the questions posed in \hyperref[djm-009K]{Informal definition \ref{djm-009K}}. Additionally, these examples either generalize or restructure categorical systems theories presented in the literature.\subsubsection{Module of open Petri nets over undirected wiring diagrams}\label{ssl-001K}\par{}Here we'll give an example of a module of systems \[\mathbb {C}\mathsf {ospan}(\mathsf {Petri})(\emptyset , \mathbb {C}\mathsf {ospan}(\mathsf {Finset})) \colon  \bullet  \mathrel {\mkern 3mu\vcenter {\hbox {$\shortmid $}}\mkern -10mu{\to }} \mathbb {C}\mathsf {ospan}(\mathsf {Finset})\] by which Petri nets compose via undirected wiring diagrams. This example gives an operadic perspective on the symmetric monoidal double category of open Petri nets introduced in \cite{baez-2021-structured}. This systems theory is generated using the restriction of \hyperref[djm-00GQ]{the port plugging doctrine} to \hyperref[djm-00G5]{a free process}.\par{}We begin by defining Petri nets.\paragraph{Petri nets}\label{ssl-001W}
\begin{definition}[{Petri net}]\label{ssl-001N}
\par{}A Petri net consists of a finite set of species \(S\), a finite set of transitions \(T\) and two functions \(s, t \colon  T \to  \mathbb {N}[S]\) assigning to each transition a source multiset of species and a target multiset of species.\end{definition}

\begin{example}[{Petri net}]\label{ssl-001O}
\par{}The following Petri net has two species \(S\) and \(I\) which represent a susceptible and infected population. The single transition represents an infection event whose source is \((S, I)\)  and  whose target is \((I,I)\). In other words, an infection takes a susceptible and an infected individual and transforms them into two infected individuals.\begin{center}
                          \includegraphics[scale=0.75]{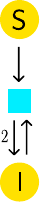}
                        \end{center}\end{example}
\par{}The double categorical perspective taken in \cite{baez-2021-structured} extends the single categorical perspective introduced in \cite{baez-2017-reaction} by introducing the concept of \emph{maps between Petri nets}.
\begin{definition}[{Maps of Petri nets}]\label{ssl-001Z}
\par{}Given Petri nets \[P = (S, T, s \colon  T \to  \mathbb {N}[S], t \colon  T \to  \mathbb {N}[S])\] and \[P' = (S', T', s' \colon  T' \to  \mathbb {N}[S'], t' \colon  T' \to  \mathbb {N}[S']),\] a \emph{map of Petri nets} \(P \to  P'\) consists of maps \(f \colon  S \to  S'\) and \(g \colon  T \to  T'\) satisfying

  \begin{center}
    \begin {tikzcd}
      T & {\mathbb {N}[S]} & T & {\mathbb {N}[S]} \\
      {T'} & {\mathbb {N}[S']} & {T'} & {\mathbb {N}[S']}
      \arrow ["s", from=1-1, to=1-2]
      \arrow ["g"', from=1-1, to=2-1]
      \arrow ["{\mathbb {N}[f]}", from=1-2, to=2-2]
      \arrow ["t", from=1-3, to=1-4]
      \arrow ["g"', from=1-3, to=2-3]
      \arrow ["{\mathbb {N}[f]}", from=1-4, to=2-4]
      \arrow ["{s'}"', from=2-1, to=2-2]
      \arrow ["{t'}"', from=2-3, to=2-4]
    \end {tikzcd}
  \end{center}\end{definition}

\begin{example}[{Maps of Petri nets}]\label{ssl-0020}
\par{}Below depicts a map from a Petri net representing an SIR model of infection to a Petri net representing an SIS model of infection, in which recovered patients are susceptible to the disease once again.\begin{center}
                          \includegraphics[scale=1]{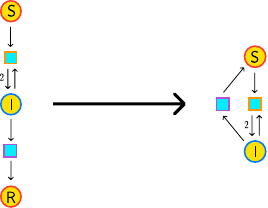}
                        \end{center}\par{}Colors outlining the species and transitions indicate the maps between species and transitions. In particular, the susceptible and recovered populations of the domain Petri net are mapped to the susceptible population of the codomain Petri net. The infected population of the domain Petri net is mapped to the intefected population of the codomain Petri net.\end{example}
\par{}Petri nets and their maps form a category.
\begin{definition}[{Category of Petri nets}]\label{ssl-003X}
\par{}Let \(\mathsf {Petri}\) be the category whose objects are \hyperref[ssl-001N]{Petri nets} and whose morphisms are \hyperref[ssl-001Z]{maps of Petri nets}.\end{definition}
\par{}From \cite{baez-2020-open} this category has small colimits and hence is \hyperref[djm-00AY]{\emph{rex}}.\paragraph{The double category of interfaces and interactions in the module of open Petri nets}\label{ssl-0049}\par{}Recall \hyperref[ssl-001I]{the target of a module of systems} is a symmetric monoidal double category whose tight category consists of interfaces and interface maps and whose loose category consists of interactions and interaction maps. For the module of open Petri nets, the double category of interfaces and interactions is the double category of cospans of finite sets, \(\mathbb {C}\mathsf {ospan}(\mathsf {Finset})\), which we define below.
\begin{definition}[{The double category \(\mathbb {C}\mathsf {ospan}(\mathsf {Finset})\)}]\label{ssl-004C}
\par{}The double category \(\mathbb {C}\mathsf {ospan}(\mathsf {Finset})\) has:
  \begin{itemize}\item{}Objects are finite sets, \(M\).
    \item{}A tight morphism is a map of finite sets, \(f \colon  M \to  M'\).
    \item{}A loose morphism is a cospan of finite sets, \(M \to  J \leftarrow  N\).
    \item{}A square is a commuting diagram:

  \begin{center}
      \begin {tikzcd}
        M & J & N \\
        {M'} & {J'} & {N'}
        \arrow [from=1-1, to=1-2]
        \arrow ["m"', from=1-1, to=2-1]
        \arrow [from=1-2, to=2-2]
        \arrow [from=1-3, to=1-2]
        \arrow ["n", from=1-3, to=2-3]
        \arrow [from=2-1, to=2-2]
        \arrow [from=2-3, to=2-2]
      \end {tikzcd}
    \end{center}\end{itemize}\end{definition}
\par{}Next we highlight how we depict the components of \(\mathbb {C}\mathsf {ospan}(\mathsf {Finset})\) as interfaces, interactions, and their maps labelled by the objects and morphisms of \(\mathbb {L}\mathsf {oose}\).
\begin{example}[{Interfaces in the module of open Petri nets}]\label{ssl-004D}
\par{}In the module of open Petri nets, interfaces and their maps are the tight category of \(\mathbb {C}\mathsf {ospan}(\mathsf {Finset})\). Hence,  an interface is a finite set and an interface map is a map of finite sets.\par{}We depict an interface \(M\) as a box with \(M\) exposed ports. For example, below is the interface \(2\).\begin{center}
                          \includegraphics[scale=0.75]{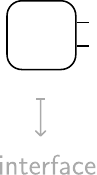}
                        \end{center}\par{}Note that our depiction shows how \(M\) is an object of \(\mathbb {C}\mathsf {ospan}(\mathsf {Petri})(\emptyset , \mathbb {C}\mathsf {ospan}(\mathsf {Finset}))\) living over the object \(\mathsf {interface}\) in \(\mathbb {L}\mathsf {oose}\).\par{}The following diagram depicts an interface map \(3 \to  2\). Colors indicate that the first and third ports (resp. second port) of the domain interface are sent to the first port (resp. second port) of the codomain interface.\begin{center}
                          \includegraphics[scale=0.5]{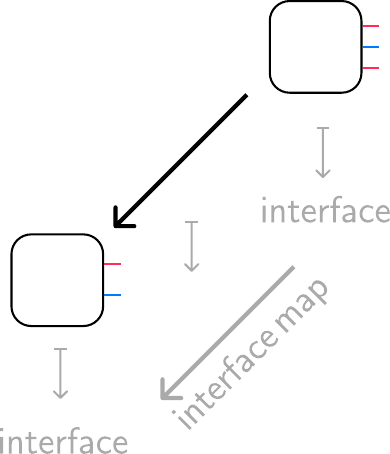}
                        \end{center}\end{example}

\begin{example}[{Interactions in the module of open Petri nets}]\label{ssl-004E}
\par{}In the module of open Petri nets, interactions and their maps are the loose category of \(\mathbb {C}\mathsf {ospan}(\mathsf {Finset})\). Hence an interaction is given by a cospan of finite sets and a map of interactions is a map of cospans.\par{}For example, an interaction that transforms two systems each with interface \(2\) into a single system with interface \(3\) is a cospan \(2 + 2 \to  J \leftarrow  3\). The following cospan is an example of such an interaction.\begin{center}
  \begin {tikzcd}[row sep=tiny]
    \bullet  \\
    \bullet  & \bullet  & \bullet  \\
    & \bullet  & \bullet  \\
    \bullet  & \bullet  & \bullet  \\
    \bullet  \\
    {2 + 2} & 3 & 3
    \arrow [maps to, from=1-1, to=2-2]
    \arrow [maps to, from=2-1, to=3-2]
    \arrow [maps to, from=2-3, to=2-2]
    \arrow [maps to, from=3-3, to=3-2]
    \arrow [maps to, from=4-1, to=3-2]
    \arrow [maps to, from=4-3, to=4-2]
    \arrow [maps to, from=5-1, to=4-2]
    \arrow [from=6-1, to=6-2]
    \arrow [from=6-3, to=6-2]
  \end {tikzcd}
\end{center}\par{}We depict this interaction as the following \emph{undirected wiring diagram} living over the loose identity \(\mathsf {interaction} \colon  \mathsf {interface} \to  \mathsf {interface}\).\begin{center}\includegraphics[scale=0.5]{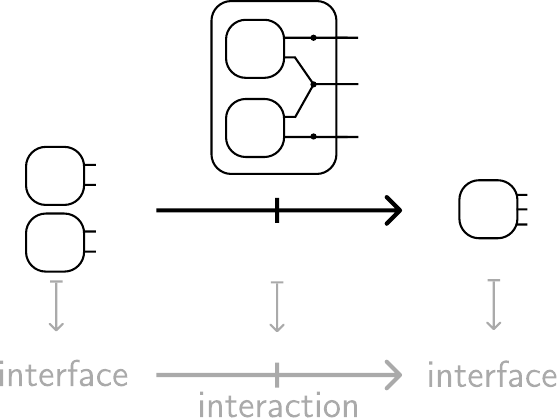}\end{center}\par{}Let \(m \colon  M \to  M'\) and \(n \colon  N \to  N'\) be two \hyperref[ssl-004D]{maps of interfaces} and let \(M \to  J \leftarrow  N\) and \(M' \to  J' \leftarrow  N'\) be two interactions. Recall, that an \emph{interaction map} in the module of open Petri nets is a map \(J \to  J'\) such that the following diagram commutes.\begin{center}
  \begin {tikzcd}
    M & J & N \\
    {M'} & {J'} & {N'}
    \arrow [from=1-1, to=1-2]
    \arrow ["m"', from=1-1, to=2-1]
    \arrow [from=1-2, to=2-2]
    \arrow [from=1-3, to=1-2]
    \arrow ["n", from=1-3, to=2-3]
    \arrow [from=2-1, to=2-2]
    \arrow [from=2-3, to=2-2]
  \end {tikzcd}
\end{center}\par{}Below is an example of an interaction map that is the inclusion of an undirected wiring diagram with 3 junctions into an undirected wiring diagram with 4 junctions.\begin{center}\includegraphics[scale=0.5]{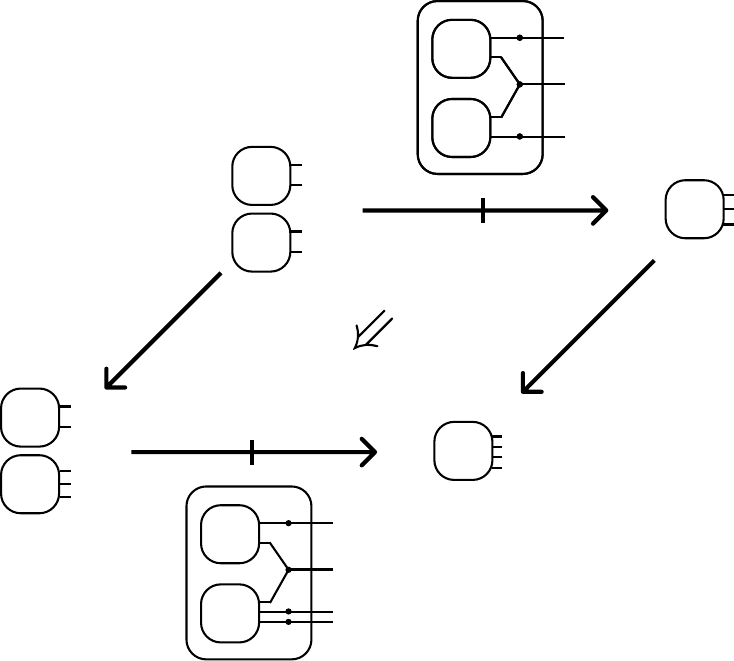}\end{center}\par{}It lives over the following  identity square in \(\mathbb {L}\mathsf {oose}\).\begin{center}
  \begin {tikzcd}[sep=huge]
    & {\mathsf {interface}} & {\mathsf {interface}} \\
    {\mathsf {interface}} & {\mathsf {interface}} & {}
    \arrow [""{name=0, anchor=center, inner sep=0}, "{\mathsf {interaction}}", "\shortmid "{marking}, from=1-2, to=1-3]
    \arrow ["{\mathsf {interface \, map}}"', from=1-2, to=2-1]
    \arrow ["{\mathsf {interface \, map}}", from=1-3, to=2-2]
    \arrow [""{name=1, anchor=center, inner sep=0}, "{\mathsf {interaction}}"', "\shortmid "{marking}, from=2-1, to=2-2]
    \arrow ["{\mathsf {interaction\, map}}"{description}, shorten <=9pt, shorten >=9pt, Rightarrow, from=0, to=1]
  \end {tikzcd}
\end{center}\end{example}
\paragraph{System and system maps in the module of open Petri nets}\label{ssl-0048}\par{}Recall that the carrier of a module of systems is a symmetric monoidal category whose objects are systems and whose morphisms are system maps. In this section, we describe the carrier of the module of open Petri nets.
\begin{example}[{Systems in the module of open Petri nets}]\label{ssl-001M}
\par{}In the module of open Petri nets, a \emph{system} is an \emph{open Petri net}, which we define below.
\begin{definition}[{Open Petri net}]\label{ssl-001Q}
\par{}An \textbf{open Petri net} consists of: 
  \begin{itemize}\item{}A finite set \(M\) representing its interface.
    \item{}A Petri net \((S, T, s \colon  T \to  \mathbb {N}[S], t \colon  T \to  \mathbb {N}[S])\).
    \item{}A map \(p \colon  M \to  S\) which determines the species that each port exposes.\end{itemize}\end{definition}
\par{}In the module of open Petri nets, an open Petri net is a loose morphism whose codomain is its interface and which lives over the walking loose arrow \(\mathsf {system} \colon  \bullet  \mathrel {\mkern 3mu\vcenter {\hbox {$\shortmid $}}\mkern -10mu{\to }} \mathsf {interface}\).\par{}Below is an example of an open Petri net with interface \(2\).  \hyperref[ssl-001O]{This Petri net models an infection event}. The first port exposes the susceptiple population \(S\) and the second port exposes the infected population \(I\).\begin{center}
                          \includegraphics[scale=0.65]{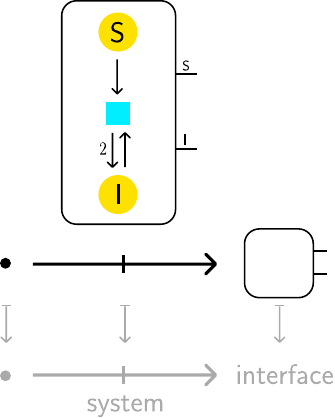}
                        \end{center}\par{}Here is a second example of an open Petri net with interface \(2\). This Petri net models a recovery event. The first port exposes the infected population \(I\) and the second port exposes the recovered population \(R\).\begin{center}
                          \includegraphics[scale=0.65]{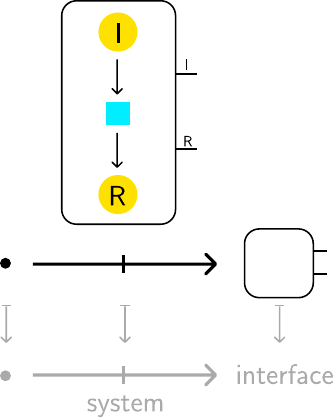}
                        \end{center}\end{example}

\begin{example}[{Monoidal product for open Petri nets}]\label{ssl-001P}
\par{}The module of open Petri nets is equipped with a symmetric monoidal structure that we can use to  "stack" systems and interfaces.\par{}For example, below is the monoidal product of the two open Petri nets representing infection and recovery events. The interface for this monoidal product system is the monoidal product of the interfaces, which in this example is given by coproduct: \(2 + 2\).\begin{center}
                          \includegraphics[scale=0.65]{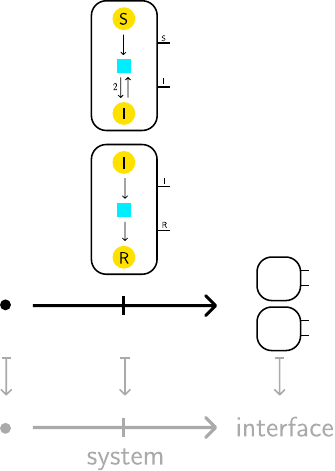}
                        \end{center}\end{example}

\begin{example}[{System maps for open Petri nets}]\label{ssl-001Y}
\par{}In the module of open Petri nets , a \emph{system map} is a map of open Petri nets, which we define below.
\begin{definition}[{Maps of open Petri nets}]\label{ssl-0021}
\par{}A map from an open Petri net with interface \(M\) \[(P = (S, T, s \colon  T \to  \mathbb {N}[S], t \colon  T \to  \mathbb {N}[S]), p \colon  M \to  S)\] to an open Petri net with interface \(M'\) \[(P' = (S', T', s' \colon  T' \to  \mathbb {N}[S'], t' \colon  T' \to  \mathbb {N}[S']), p' \colon  M' \to  S),\] along an interface map \[m \colon  M \to  M'\] consists of a \hyperref[ssl-001Z]{map of Petri nets} \((f \colon  S \to  S', g \colon  T \to  T') \colon  P \to  P'\) such that the following diagram commutes

  \begin{center}
    \begin {tikzcd}
      {M } & S \\
      {M'} & {S'}
      \arrow ["p", from=1-1, to=1-2]
      \arrow ["m"', from=1-1, to=2-1]
      \arrow ["f", from=1-2, to=2-2]
      \arrow ["{p'}"', from=2-1, to=2-2]
    \end {tikzcd}
  \end{center}\end{definition}
\par{}Note that  a system map is \emph{between} two \hyperref[ssl-001M]{systems} and  \emph{along} a \hyperref[ssl-004D]{map between the interfaces} of those systems. Furthermore, it is labelled by the identity square (below called \(\mathsf {system\,map}\)) for the walking loose arrow (below called \(\mathsf {system}\)) in \(\mathbb {L}\mathsf {oose}\):

  \begin{center}
      \begin {tikzcd}[column sep=large]
      \bullet  & {\mathsf {interface}} \\
      \bullet  & {\mathsf {interface}}
      \arrow [""{name=0, anchor=center, inner sep=0}, "{\mathsf {system}}", "\shortmid "{marking}, from=1-1, to=1-2]
      \arrow [Rightarrow, no head, from=1-1, to=2-1]
      \arrow ["{\mathsf {interface\,map}}", from=1-2, to=2-2]
      \arrow [""{name=1, anchor=center, inner sep=0}, "{\mathsf {system}}"', "\shortmid "{marking}, from=2-1, to=2-2]
      \arrow ["{\mathsf {system\, map}}"{description}, draw=none, from=0, to=1]
    \end {tikzcd}
  \end{center}\par{}Below is an example of a system map between an open SIR Petri net and an open SIS Petri net along the interface map exemplified \hyperref[ssl-004D]{here}.\begin{center}
                          \includegraphics[scale=0.5]{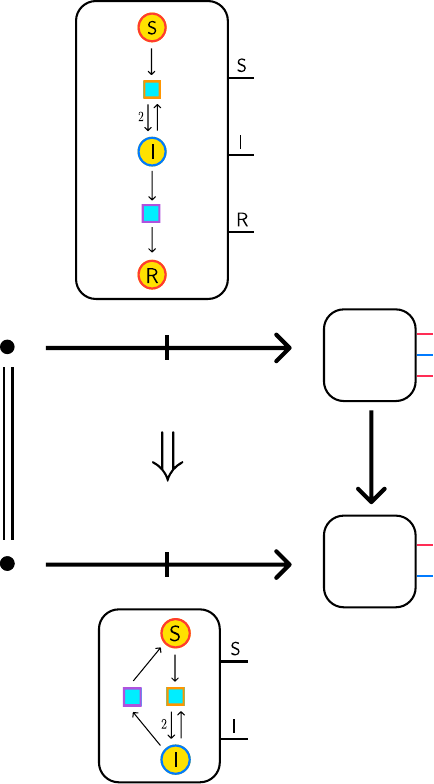}
                        \end{center}\par{}This type of system map has the quality of \emph{abstraction} because it abstracts away the distinction between the susceptible and recovered populations while retaining the overal dynamics of the system.\par{}Next we have  an example of a system map from an open SIR Petri net to a open SIRD Petri net --- in which infected patients either recover or die --- along the indicated map of interfaces.\begin{center}
                          \includegraphics[scale=0.5]{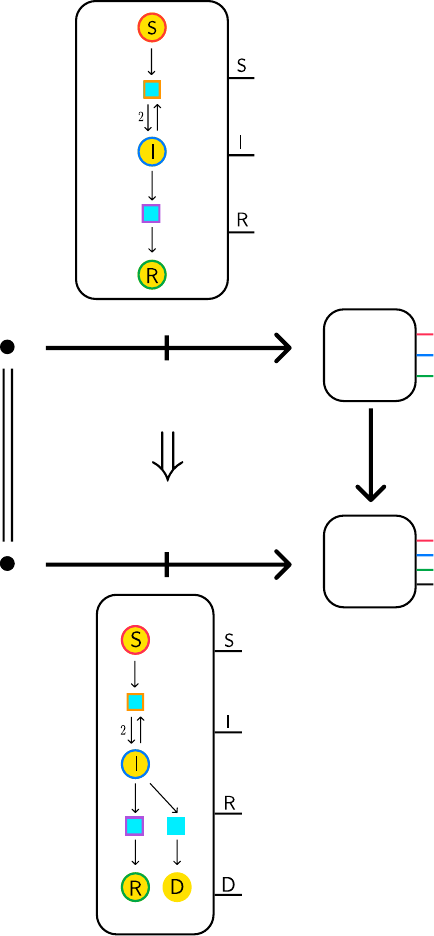}
                        \end{center}\par{}This type of system map has the quality of \emph{embedding} because it shows how the dynamics of the original SIR model are embedded in a larger model.\end{example}

\begin{example}[{The action of interactions on systems in the module of open Petri nets}]\label{ssl-004A}
\par{}Finally, we show how interactions (and their map) act on systems (and their maps) in the module of open Petri nets.
\begin{example}[{The action of undirected wiring diagrams on open Petri nets}]\label{ssl-001T}
\par{}Recall that \(\mathsf {interaction}\) is the loose identity on \(\mathsf {interface}\). Therefore, its composite with the walking loose arrow \[\mathsf {system} \colon  \bullet  \mathrel {\mkern 3mu\vcenter {\hbox {$\shortmid $}}\mkern -10mu{\to }} \mathsf {interface}\] is again the walking loose arrow \(\mathsf {system}\). In other words, the diagram below commutes.

  \begin{center}
    \begin {tikzcd}[column sep = huge]
      \bullet  & {\mathsf {interface}} & {\mathsf {interface}}
      \arrow ["{\mathsf {system}}"', "\shortmid "{marking}, from=1-1, to=1-2]
      \arrow ["{\mathsf {system}}", "\shortmid "{marking}, curve={height=-24pt}, from=1-1, to=1-3]
      \arrow ["{\mathsf {interaction}}"', "\shortmid "{marking}, from=1-2, to=1-3]
    \end {tikzcd}
  \end{center}\par{}We intrepret this fact as "an interaction transforms one system into another system".\par{}Consider the following composite in the module of open Petri nets.\begin{center}
                          \includegraphics[scale=0.5]{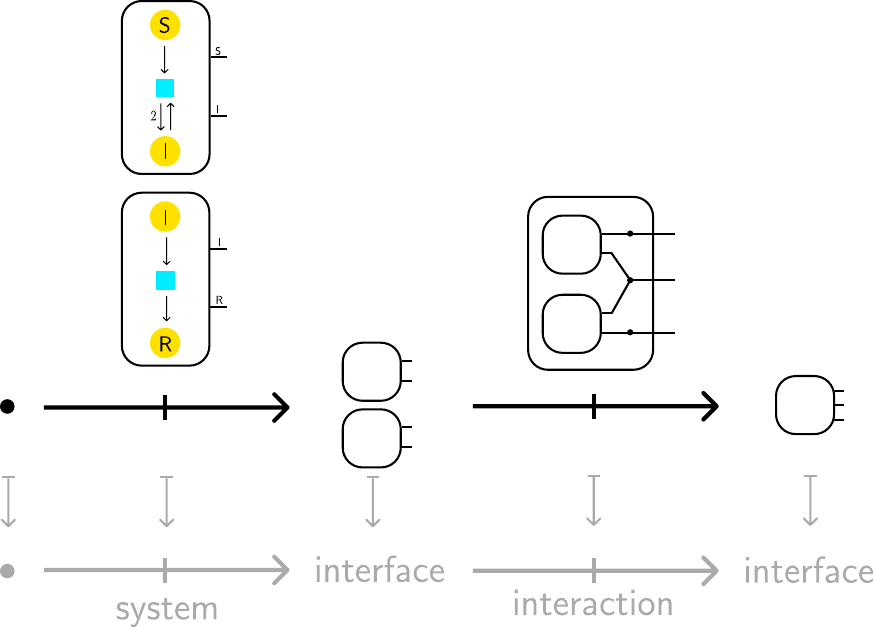}
                        \end{center}\par{}The undirected wiring diagram living over \(\mathsf {interaction}\) transforms the two systems representing infection and recovery events living over \(\mathsf {system}\). As we will see mathematically in a future section, it identifies the infected species exposed by the two open Petri nets.\par{}The composite of this process with these systems is again a system:  an open Petri net representing an SIR infection model.\begin{center}
                          \includegraphics[scale=0.65]{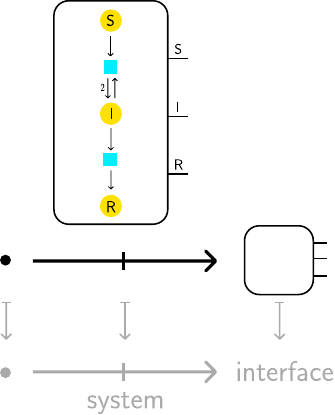}
                        \end{center}\end{example}

\begin{example}[{Maps of composites}]\label{ssl-0026}
\par{}	We saw how \hyperref[ssl-001T]{interactions act on systems} follows from the fact that in \(\mathbb {L}\mathsf {oose}\) the composite
	\[
		\bullet  \xrightarrow {\mathsf {system}} \mathsf {interface} \xrightarrow {\mathsf {interaction}} \mathsf {interface}
	\] is the map \(\bullet  \xrightarrow {\mathsf {system}} \mathsf {interface}\).\par{}	Likewise, maps of interactions act on maps of systems because in \(\mathbb {L}\mathsf {oose}\) the composite of the squares

  \begin{center}
	\begin {tikzcd}[column sep=huge]
		\bullet  & {\mathsf {interface}} & {\mathsf {interface}} \\
		\bullet  & {\mathsf {interface}} & {\mathsf {interface}}
		\arrow [""{name=0, anchor=center, inner sep=0}, "{\mathsf {system}}", from=1-1, to=1-2]
		\arrow [equals, from=1-1, to=2-1]
		\arrow [""{name=1, anchor=center, inner sep=0}, "{\mathsf {interaction}}", from=1-2, to=1-3]
		\arrow [from=1-2, to=2-2]
		\arrow [from=1-3, to=2-3]
		\arrow [""{name=2, anchor=center, inner sep=0}, "{\mathsf {system}}"', from=2-1, to=2-2]
		\arrow [""{name=3, anchor=center, inner sep=0}, "{\mathsf {interaction}}"', from=2-2, to=2-3]
		\arrow ["{\mathsf {system\, map}}"{description}, shorten <=4pt, shorten >=4pt, Rightarrow, from=0, to=2]
		\arrow ["{\mathsf {interaction \,map}}"{description}, shorten <=4pt, shorten >=4pt, Rightarrow, from=1, to=3]
	\end {tikzcd}
\end{center}

is the square

  \begin{center}
		\begin {tikzcd}[column sep=huge]
			\bullet  & {\mathsf {interface}} \\
			\bullet  & {\mathsf {interface}}
			\arrow [""{name=0, anchor=center, inner sep=0}, "{\mathsf {system}}", from=1-1, to=1-2]
			\arrow [equals, from=1-1, to=2-1]
			\arrow [from=1-2, to=2-2]
			\arrow [""{name=1, anchor=center, inner sep=0}, "{\mathsf {system}}"', from=2-1, to=2-2]
			\arrow ["{\mathsf {system\, map}}"{description}, shorten <=4pt, shorten >=4pt, Rightarrow, from=0, to=1]
		\end {tikzcd}
\end{center}\par{}	Below is an example of a map of systems (left) acted on by a map of processes (right). For simplicity we have omitted the colors indicating the details of these maps and leave them to be inferred by the reader.\begin{center}
                          \includegraphics[scale=0.5]{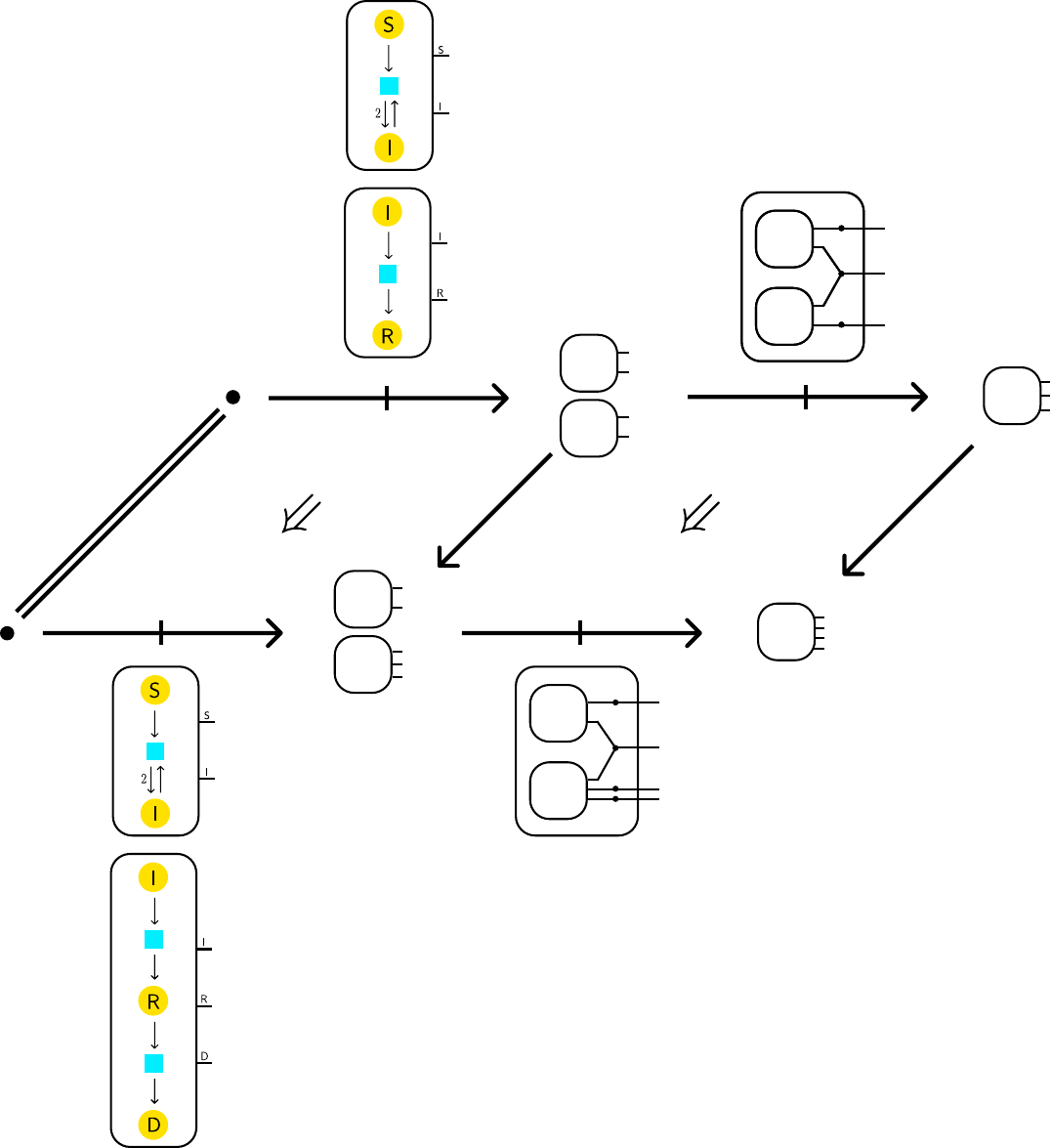}
                        \end{center}\par{}	This action produces the following map of systems.\begin{center}
                          \includegraphics[scale=0.5]{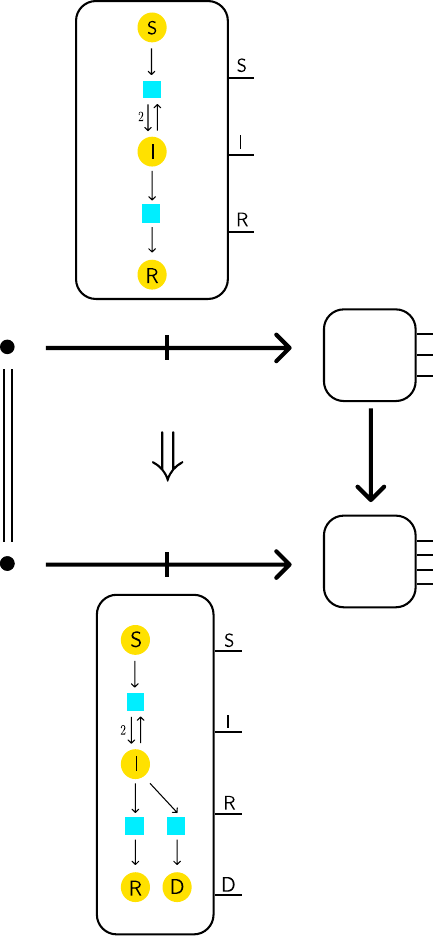}
                        \end{center}\end{example}
\end{example}
\subsubsection{Module of Moore machines over lenses}\label{ssl-0029}\par{}Here we'll give an example of a module of systems 
  \[\mathbb {M}\mathsf {oore}(\mathsf {proj}_{\mathsf {Set}}, (-)^2) \colon  \bullet  \mathrel {\mkern 3mu\vcenter {\hbox {$\shortmid $}}\mkern -10mu{\to }} \mathbb {L}\mathsf {ens}(\mathsf {proj}_{\mathsf {Set}})\]
  by which deterministic Moore machines composing via lenses. For more on Moore machines, the monograph \cite{jaz-2021-book} goes into great detail about these sorts of systems and their behaviors in its first few chapter. In particular, the kinds of Moore machines we will cover here (the traditional ones) are covered in Chapter 1 of \cite{jaz-2021-book}, their non-deterministic variants are covered in Chapter 2, and their maps and behaviors are covered in Chapter 3. In this paper, we will see Moore machines and their variants again in the guise of \hyperref[djm-00G6]{\emph{generalized Moore machines}}.\par{}The construction of this module of systems is defined in \hyperref[djm-00EE]{Section \ref{djm-00EE}}. In this section we exemplify its systems, interactions, and their maps.\paragraph{The double category of interfaces and interactions in the module of deterministic Moore machines}\label{ssl-005F}
\begin{explication}[{The double category \(\mathbb {L}\mathsf {ens}(\mathsf {proj}_{\mathsf {Set}})\)}]\label{ssl-005G}
\par{}The double category of interactions for the module of deterministic Moore machines is the \hyperref[djm-00B6]{lense double category} \(\mathbb {L}\mathsf {ens}(\mathsf {proj}_{\mathsf {Set}})\) where  \(\mathsf {proj}_{\mathsf {Set}}\) is the \hyperref[djm-00E4]{simple fibration} associated with \(\mathsf {Set}\), which has:
  \begin{itemize}\item{}An object is a pair of sets which we denote \({A^\# \choose  A}\).
    \item{}A loose morphism is a pair of set maps \(f \colon  A \to  B\) and \(f^\# \colon  A \times  B^\# \to  A^\#\) which we notate as a lens \({f^\# \choose  f} \colon   {A^{\sharp } \choose  A} \leftrightarrows  {B^{\sharp } \choose  B}\). See also \hyperref[djm-00E6]{Example \ref{djm-00E6}}.
    \item{}A tight morphism is a pair of set maps  \(f \colon  A \to  B\) and \(f^\# \colon  A^\# \times  B \to  B^\#\) which we notate as a  chart \({f^\# \choose  f} \colon   {A^{\sharp } \choose  A} \rightrightarrows  {B^{\sharp } \choose  B}\).\end{itemize}\end{explication}

\begin{example}[{Interactions in the module of deterministic Moore machines}]\label{ssl-005K}
\par{}Since interactions are lenses interaction that transforms a \hyperref[ssl-002A]{Moore machine} with interface \({I \choose  O}\) into a Moore machine with interface \({I' \choose  O'}\) consists of:
  \begin{itemize}\item{}A map of outputs \(f \colon  O \to  O'\).
    \item{}A map of inputs in the reverse order \(f^\# \colon  I' \times  O \to  I\).\end{itemize}\par{}Let \(2\) be the two element set \(\{0, 1\}\). There is an interaction \[{2 \choose  2} \times  {2 \choose  2} \mathrel {\mkern 3mu\vcenter {\hbox {$\shortmid $}}\mkern -10mu{\to }} {2 \choose  2 \times  2}\] where 
  \begin{itemize}\item{}The map on outputs is the swap \((b_1, b_2) \mapsto  (b_2,b_1)\). 
    \item{}The map of inputs is given by  \((a, (b_1,b_2) )\mapsto  (a, b_1)\).\end{itemize}\par{}This interaction is depicted by the directed wiring diagram below where each wire carries a copy of \(2\).\begin{center}
                              \includegraphics[scale=0.5]{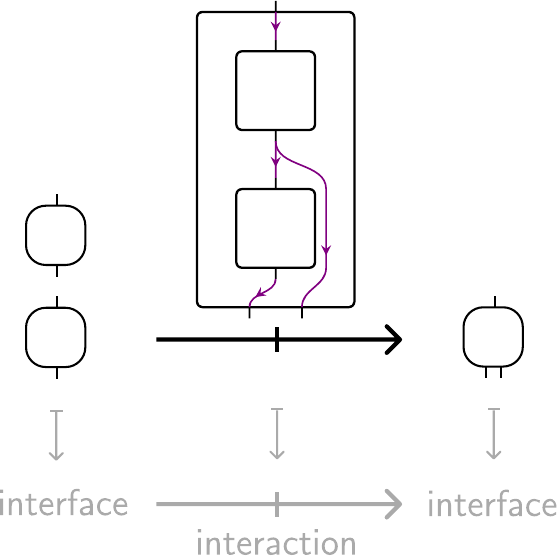}
                            \end{center}\label{ssl-002F}\par{}This interaction transforms two Moore machines (called the components) with interface \({2 \choose  2}\) into a single Moore machine (called the composite) with interface \({2 \choose  2 \times  2}\).  
  \begin{itemize}\item{}The composite machine outputs the pair of outputs generated by the components.
    \item{}Given an input to the composite machine, that input is sent directly to the first machine and the output of the first component machine is used as the input to the second component machine.\end{itemize}
  We will show how this interaction acts on systems in \hyperref[ssl-002B]{Example \ref{ssl-002B}}.\end{example}
\paragraph{Systems and the action of interactions in the module of deterministic Moore machines}\label{ssl-005I}
\begin{definition}[{Open deterministic Moore machines}]\label{ssl-002A}
\par{}A deterministic Moore machine with interface \({I \choose  O}\) consists of:
  \begin{itemize}\item{}A set of states \(S\).
    \item{}An update \(u \colon  I \times  S \to  S\). 
    \item{}A readout \(r \colon  S \to  O\).\end{itemize}\par{}Note that a deterministic Moore machine is equivalently a \hyperref[ssl-005G]{lens} \[{u \choose  r} \colon  {S \choose  S} \leftrightarrows  {I \choose  O}\]\end{definition}

\begin{example}[{Mod 2 counter}]\label{ssl-005L}
\par{}Below is a \hyperref[ssl-002A]{deterministic Moore machine} with interface \({2 \choose  2}\). It has two states, its update function counts modulo 2, and its readout is the identity.\begin{center}
                              \includegraphics[scale=0.65]{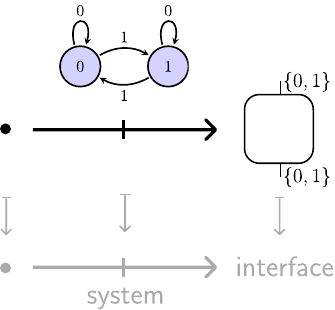}
                            \end{center}\label{ssl-002I}\end{example}

\begin{example}[{Mod 4 counter}]\label{ssl-005M}
\par{}Below is a \hyperref[ssl-002A]{deterministic Moore machine} with interface \({2 \choose  2}\). It has 4 states, its update function counts modulo 4, and its readout is the parity of the state.\begin{center}
                              \includegraphics[scale=0.65]{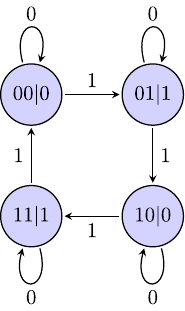}
                            \end{center}\label{ssl-005N}\end{example}

\begin{example}[{The action of lenses on deterministic Moore machines}]\label{ssl-002B}
\par{}Consider composing two \hyperref[ssl-002A]{mod 2 counters} according to the interaction described in \hyperref[ssl-005K]{Example \ref{ssl-005K}}. The result is a mod 4 counter.\begin{center}
                              \includegraphics[scale=0.5]{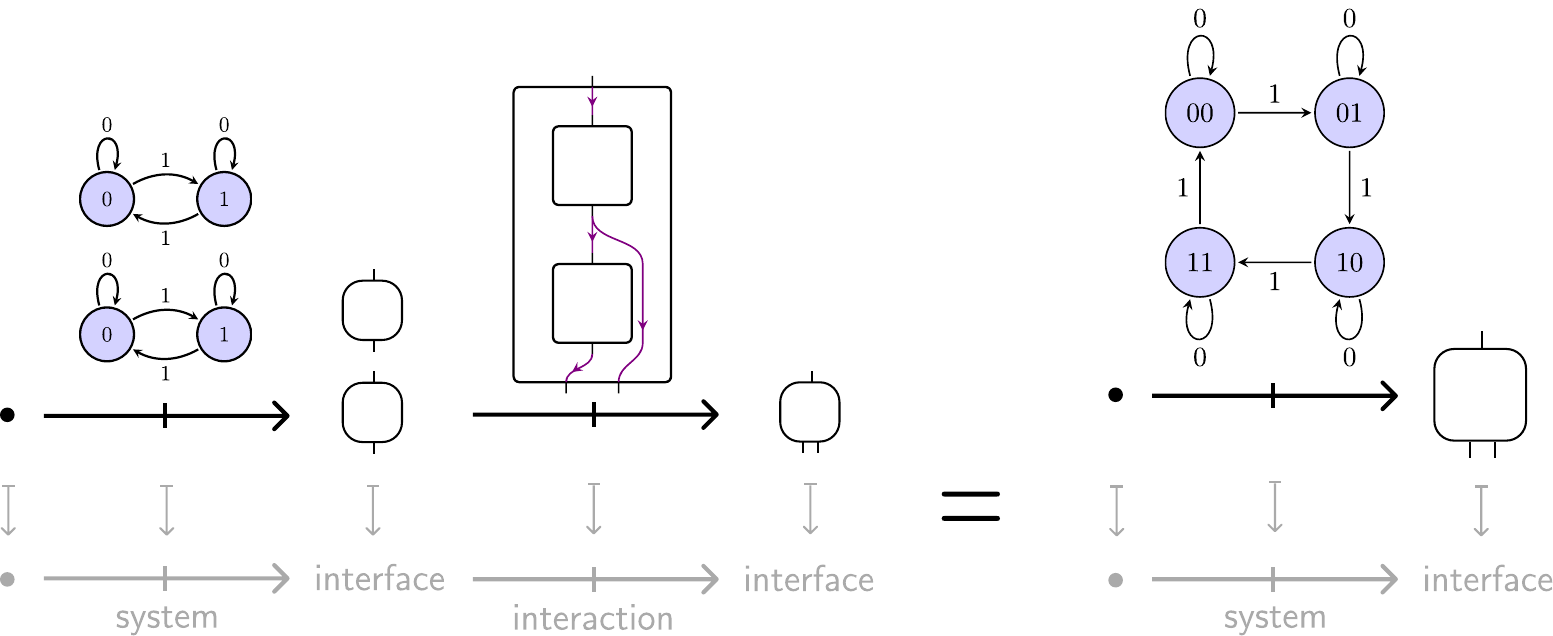}
                            \end{center}\label{ssl-002G}\par{}A state of the composite system is the product of the states of the component systems. The first component system is tracking the \(1\)s digit in the binary representation of a number and the second component system is tracking the \(2\)s digit.\par{}So far we have seen two examples of Moore machines: \hyperref[ssl-002I]{the mod 2 counter} and \hyperref[ssl-002G]{the mod 4 counter}. In both of these examples, the read out is an isomorphism, which means that the whole state is exposed by the interface.\par{}Consider composing \hyperref[ssl-002J]{the mod 4 counter} with an interaction that \({2 \choose  2 \times  2} \mathrel {\mkern 3mu\vcenter {\hbox {$\shortmid $}}\mkern -10mu{\to }} {2 \choose  2}\) that forgets the first output. The result is a mod 4 counter whose \hyperref[ssl-005M]{output exposes only the parity of each state}. The states \(00\) and \(10\) output \(0\) while the states \(01\) and \(11\) output \(1\).\begin{center}
                              \includegraphics[scale=0.5]{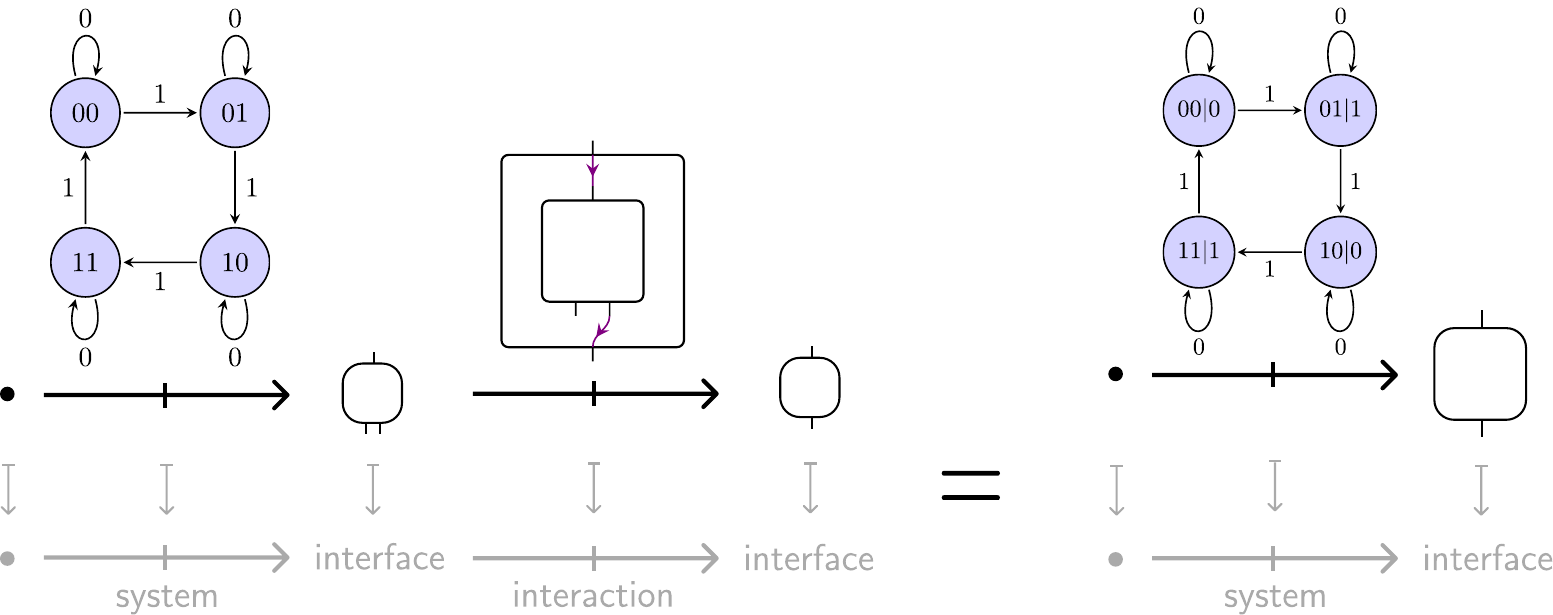}
                            \end{center}\label{ssl-002J}\end{example}
\paragraph{System maps in the module of deterministic Moore machines}\label{ssl-005H}\par{}In this section, we will show how maps of Moore machines can represent trajectories.
\begin{example}[{Trajectories in deterministic Moore machine}]\label{ssl-002O}
\par{}A trajectory in a Moore machine is a sequence of inputs and a sequence of states that evolves according to the sequence of inputs. For example, in the \hyperref[ssl-002I]{mod 2 counter}, one example trajectory consists of the input sequence \[ \cdots  \; 1 \; 1 \; 1 \;1\;1\;1 \; \cdots  \] and the state sequence \[ \cdots \;  0 \; 1 \;0 \; 1 \;0  \;  \cdots \] because if the mod 2 counter receives a sequence of input it will cycle through the four states. 

  Another example of a trajectory consists of the input sequence  \[ \cdots  \; 0\; 0 \; 0 \; 0 \; 0 \;  \cdots  \] and the state sequence \[\cdots  \; 0\; 0 \; 0 \; 0 \; 0 \; \cdots \] because if the system starts in state \(0\) and receives only \(0\) inputs, then it will remain in state \(0\). Note that the input sequence  \[ \cdots  \; 0\; 0 \; 0 \; 0 \; 0 \; \cdots  \] and the state sequence \[\cdots  \; 1\; 1 \; 1 \; 1\; 1 \; \cdots .\]

  This example shows that the input sequence alone does not uniquely determine a trajectory.\end{example}
\par{}In the remainder of this section we will show how a trajectory in a Moore machine is represented by system maps of a certain sort.
\begin{example}[{The timeline Moore machine}]\label{ssl-002Q}
\par{}Consider the \hyperref[ssl-002A]{Moore machine} below. Its interface is \({1 \choose  \mathbb {Z}}\) and states \(\mathbb {Z}\). The single input increments the state and the readout outputs the entire state.\begin{center}
                              \includegraphics[scale=0.65]{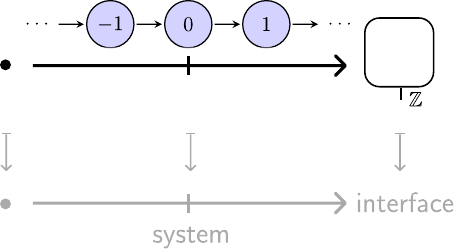}
                            \end{center}\label{ssl-002K}\par{}This Moore machine functions like a timeline. Each state is a time and the input incremements the clock.\end{example}

\begin{example}[{Interface maps of Moore machines}]\label{ssl-002P}
\par{}A map of interfaces \({I_1 \choose  O_1} \to  {I_2 \choose  O_2}\) consists of 
  \begin{itemize}\item{}A map of outputs \( O _1\to  O_2\).
    \item{}A map of inputs \( I _1\times  O_1 \to  I_2\).\end{itemize}\par{}Given a Moore machine with interface \({I \choose  O}\). We can index the trajectories by the behavior of the trajectory on its interface. The behavior of a trajectory on its interface is a sequence of inputs and outputs. This behavior is captured by a map of interface \({1 \choose  \mathbb {Z}} \to  {I \choose  O}\).\par{}For example, below is a map of interfaces \({1 \choose  \mathbb {Z}} \to  {2 \choose  2}\) consisting of the input sequence of all \(1\)s and output sequence that alternates between \(0\) and \(1\). This map defines a behavior on the interface \({2 \choose  2}\) that is: output alternating \(0\)s and \(1\)s when receiving the input \(1\).\begin{center}
                              \includegraphics[scale=0.65]{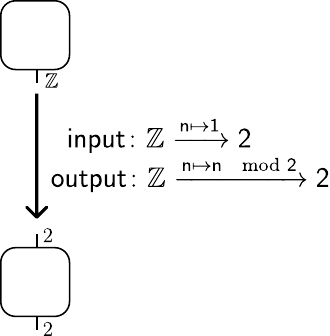}
                            \end{center}\label{ssl-002L}\end{example}

\begin{example}[{Maps of Moore machines}]\label{ssl-002H}
\par{}A map of Moore machines along a given \hyperref[ssl-002P]{interface map} consists of a map of states such that at each state the behavior of the systems on their interfaces are compatible according to the correspondence defined by the interface map.\par{}Explicitly, for \(i = 1, 2\) suppose we have \hyperref[ssl-002A]{Moore machines} \[(S_i, u_i \colon  I_i \times  S_i \to  S_i, r_i \colon  S_i \to  O_i)\]  with interfaces \({I_i \choose  O_i}\). A \emph{map of Moore machines} along  \hyperref[ssl-002P]{the map of interfaces} \[(h \colon  O_1 \to  O_2, h_\# \colon  I_1 \times  O_1 \to  I_2)\] consists of a map of states \(s \colon  S_1 \to  S_2\) satisfying:

  \begin{center}
    \begin {tikzcd}
      {S_1} & {O_1} & {I_1 \times  S_1} & {S_1} \\
      {S_2} & {O_2} & {I_1 \times  O_1 \times  S_2} \\
      && {I_2 \times  S_2} & {S_2}
      \arrow ["{r_1}", from=1-1, to=1-2]
      \arrow ["s"', from=1-1, to=2-1]
      \arrow ["h"', from=1-2, to=2-2]
      \arrow ["{u_1}", from=1-3, to=1-4]
      \arrow ["{I_1 \times  (r_1, s)}"', from=1-3, to=2-3]
      \arrow ["s", from=1-4, to=3-4]
      \arrow ["{r_2}"', from=2-1, to=2-2]
      \arrow ["{h_\# \times  S_2}"', from=2-3, to=3-3]
      \arrow ["{u_2}"', from=3-3, to=3-4]
    \end {tikzcd}
  \end{center}\par{}Below is an example of a system map from the \hyperref[ssl-002K]{timeline Moore machine} to \hyperref[ssl-002I]{the mod 2 counter} along \hyperref[ssl-002L]{the interface map} depicted above.  The map of states defines a sequence of states (in other words, a trajectory) that produces the desired behavior on the interface.\begin{center}
                              \includegraphics[scale=0.49]{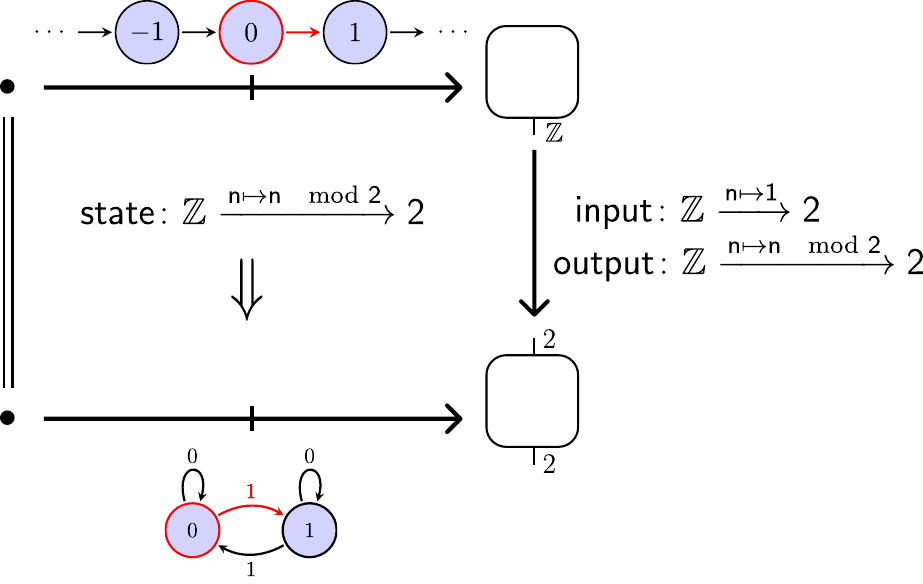}
                            \end{center}\label{ssl-002R}\par{}In this example the map of states was uniquely defined by the interface map. However, the next example shows that this is not always the case.\par{}Using \hyperref[ssl-002L]{the same interface map}, we can instead define a trajectory of the \hyperref[ssl-002J]{mod 4 counter that outputs only the parity of the state} as a system map out of \hyperref[ssl-002K]{the timeline Moore machine}. There are two such system maps because states \(00\) and state \(10\) produce the same behavior on the interface.\begin{center}
                            \includegraphics[width=0.45\textwidth]{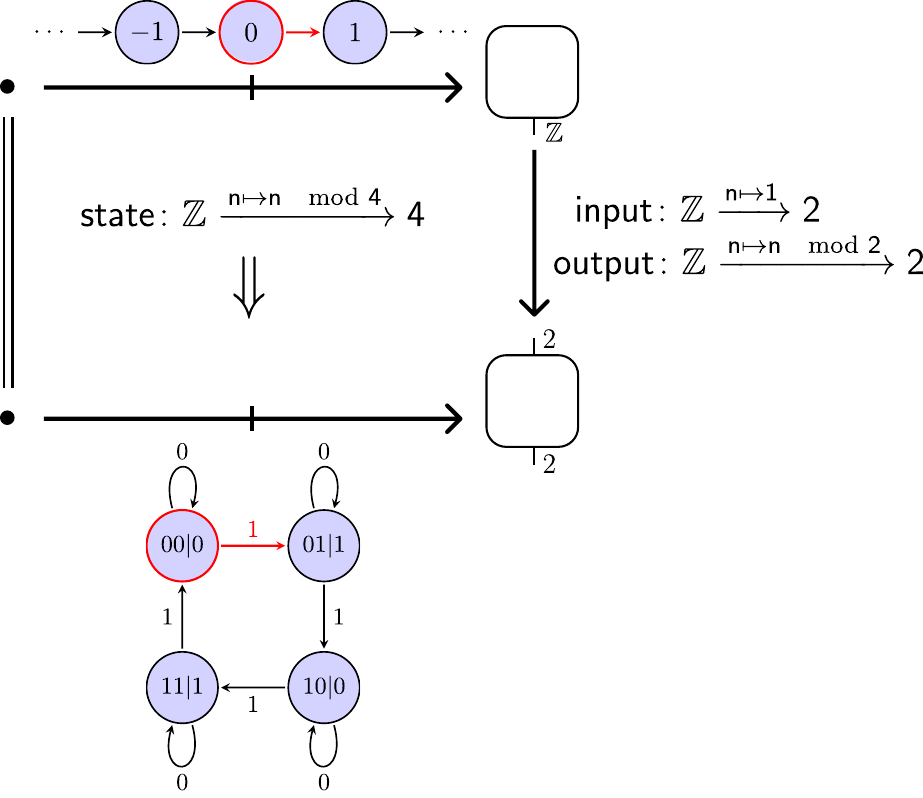}\label{ssl-002M}\hfill\includegraphics[width=.45\textwidth]{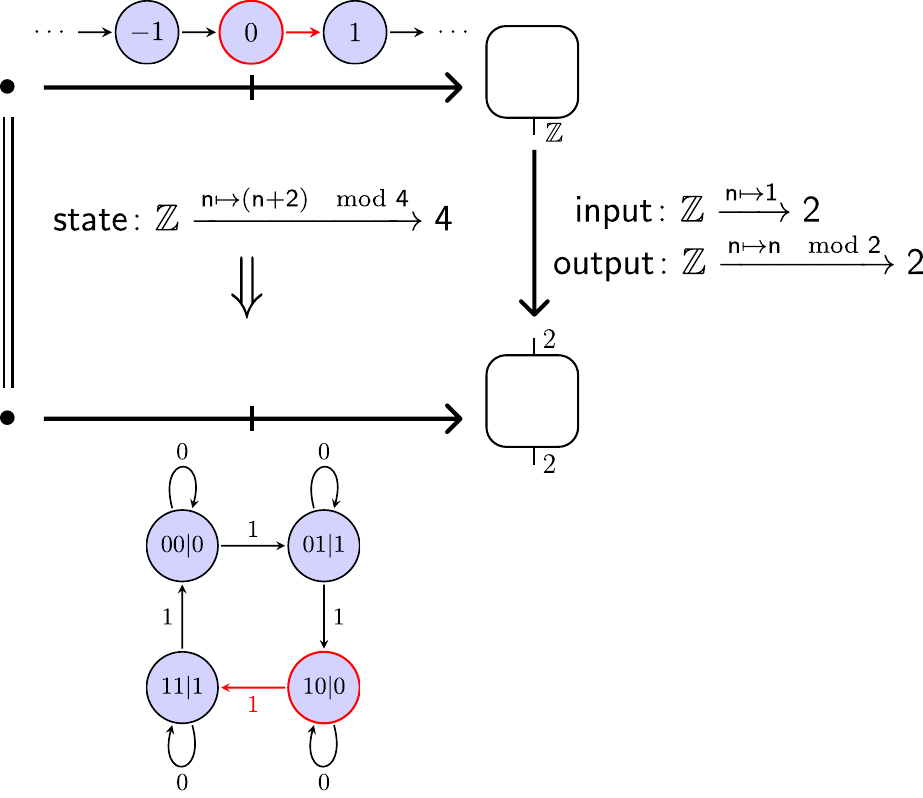}
                            \end{center}\label{ssl-002N}\end{example}
\section{Constructing modules of systems via doctrines}\label{ssl-003W}\par{}	In this section, we will introduce the notion of a \emph{doctrine of systems theories} which is a formula for constructing particular classes of modules of systems from more primitive data than the module itself. In the following sections we will give examples of doctrines which reconstruct modules of systems that are prevalent in the literature.\subsection{Doctrines of systems theories}\label{ssl-003Y}\par{}We begin by giving a definition of a doctrines of systems theory. This definition is intentionally lightweight.
\begin{definition}[{Doctrine of systems theories}]\label{djm-00A8}
\par{}A \textbf{doctrine of systems theories} consists of a \hyperref[djm-00AN]{cartesian 2-functor} \(\pi _{\mathcal {D}\mathsf {}} : \mathcal {D}\mathsf {}_{\mathsf {sys}} \to  \mathcal {D}\mathsf {}_{\mathsf {inter}}\) together with a (strictly) commuting square of cartesian pseudo-functors\begin{figure}[H] 
	\begin{center}
		\begin {tikzcd}
			{\mathcal {D}\mathsf {}_{\mathsf {sys}}} & {{\ell }\mathcal {M}\mathsf {od}_{\mathsf {r}}} \\
			{\mathcal {D}\mathsf {}_{\mathsf {inter}}} & {\mathcal {D}\mathsf {bl}}
			\arrow ["{\mathbb {S}}", from=1-1, to=1-2]
			\arrow ["{\pi _{\mathcal {D}\mathsf {}}}"', from=1-1, to=2-1]
			\arrow ["{(-)_1}",from=1-2, to=2-2]
			\arrow ["{\mathbb {I}}"', from=2-1, to=2-2]
		\end {tikzcd}
\end{center}\caption{Doctrine of systems theories}\label{ssl-0056} \end{figure} into the \hyperref[djm-004M]{2-categories of loose right modules} and double categories.\par{}	We often notate such a doctrine by the triple, \[(\pi _{\mathcal {D}\mathsf {}} \colon  \mathcal {D}\mathsf {}_{\mathsf {sys}} \to  \mathcal {D}\mathsf {}_{\mathsf {inter}}, \mathbb {S}, \mathbb {I}) .\]\par{}	A doctrine induces pseudofunctors \[\mathcal {S}\mathsf {M}(\mathcal {D}\mathsf {}_{\mathsf {sys}}) \to  \mathcal {S}\mathsf {M}({\ell }\mathcal {M}\mathsf {od}_{\mathsf {r}})\] and \[\mathcal {C}\mathsf {art}(\mathcal {D}\mathsf {}_{\mathsf {sys}}) \to  \mathcal {C}\mathsf {art}({\ell }\mathcal {M}\mathsf {od}_{\mathsf {r}}),\] and similarly on \(\mathcal {D}\mathsf {}_{\mathsf {inter}}\) which produce the actual modules of systems.\par{}	For a symmetric monoidal object \(T\) of \(\mathcal {D}\mathsf {}_{\mathsf {sys}}\), 
	\begin{itemize}\item{}We call \(T\) a \textbf{systems theory} in the doctrine \(\mathcal {D}\mathsf {}_{\mathsf {sys}}\).
		\item{}We call \(\pi _{\mathcal {D}}(T)\) the \textbf{interaction theory} for the systems theory \(T\).
		\item{}The systems theory \(T\) defines a module of systems \(\mathbb {S}(T)\) over thte double category of interactions \(\mathbb {I}(\pi _{\mathcal {D}}(T))\).\end{itemize}\end{definition}

\begin{explication}[{Doctrine of systems theories}]\label{ssl-004F}
\par{}	We think of a doctrine \(\mathbb {S} \colon  \mathcal {D}\mathsf {}_{\mathsf {sys}} \to  {\ell }\mathcal {M}\mathsf {od}_{\mathsf {r}}\) as a way of answering the questions of \hyperref[djm-009K]{Informal definition \ref{djm-009K}} about what it means \emph{to be a system}. 
	\begin{enumerate}\item{}			An object \(T \in  \mathcal {S}\mathsf {M}(\mathcal {D}\mathsf {}_{\mathsf {sys}})\) is a \emph{theory of \(\mathcal {D}\)-systems}. The object \(T\) itself is really the underlying \emph{data} needed to specify precise answers to the questions of \hyperref[djm-009K]{Informal definition \ref{djm-009K}}. The object \(\pi _{\mathcal {D}}(T) \in  \mathcal {S}\mathsf {M}(\mathcal {D}\mathsf {}_{\mathsf {inter}})\) is the data needed to specify the interactions.
		\item{}			The pseudo-functor \(\mathbb {S}\) then takes this data \(T\) and produces a module of systems \[\mathbb {S}(T) \colon  \bullet  \mathrel {\mkern 3mu\vcenter {\hbox {$\shortmid $}}\mkern -10mu{\to }} \mathbb {I}(\pi _{\mathcal {D}}(T))\] in which systems are acted on by the symmetric monoidal double category \(\mathbb {I}(\pi _{\mathcal {D}}(T))\) of interactions.\end{enumerate}\end{explication}

\begin{remark}[{On the minimal notion of doctrine}]\label{djm-00A9}
\par{}Our definition of doctrine is rather minimal. We choose to go with this minimal definition so that it can act as an organizing principle, rather than as a mathematical object of study in its own right (which might require us to discover further structure carried by doctrines of systems theories in particular, above their 2-functoriality). That is to say, we will construct examples of doctrines as defined in \hyperref[djm-00A8]{Definition \ref{djm-00A8}}, but we will not investigate any higher category of doctrines.\par{}Nevertheless, even with such a minimal definition of doctrine, the notion will help us organize the vast array of systems theories in use by applied category theorists. We may make use of morphisms of doctrines given by precomposition to express that one doctrine is a special case of another. In the upcoming \hyperref[djm-00G5]{Section \ref{djm-00G5}}, we will perform a construction at the doctrine level to restrict a doctrine to wiring diagrams or free processes.\end{remark}
\subsection{Doctrine of initial processes}\label{ssl-003Z}\par{}Given a symmetric monoidal double category \((\mathbb {D}, \otimes , 1)\), its  \[\mathsf {Hom}^{l}(\mathbb {D}) \colon  \mathbb {D} \mathrel {\mkern 3mu\vcenter {\hbox {$\shortmid $}}\mkern -10mu{\to }} \mathbb {D}\] is symmetric monoidal as well. To turn this loose bimodule into a module of systems, we collapse the left action. We do this essentially by giving a \hyperref[djm-00CM]{symmetric monoidal niche},
  
\begin{center}
    \begin {tikzcd}
      \bullet  & {\mathbb {D}} \\
      {\mathbb {D}} & {\mathbb {D}}
      \arrow ["1"', from=1-1, to=2-1]
      \arrow [equals, from=1-2, to=2-2]
      \arrow ["{\mathsf {Hom}^{l}(\mathbb {D})}"', "\shortmid "{marking}, from=2-1, to=2-2]
    \end {tikzcd}
  \end{center}

  and taking its restriction \[\mathsf {Hom}^{l}(\mathbb {D})(1, \mathbb {D}) \colon  \bullet  \mathrel {\mkern 3mu\vcenter {\hbox {$\shortmid $}}\mkern -10mu{\to }} \mathbb {D}.\]\par{}The constraint that this niche is symmetric monoidal implies that: 
  \begin{itemize}\item{}The map \(\bullet  \to  \mathbb {D}\) must map the single object of \(\bullet \) to \(1\).
    \item{}The canonical isomorphism \(1 \otimes  1 \cong  1\) is a \hyperref[djm-00EX]{conjoint commuter transform}.\end{itemize}\par{}This first of these features implies that the systems in the module of systems \(\mathsf {Hom}^{l}(\mathbb {D})(1, \mathbb {D})\) are loose morphisms \(x \colon  1 \mathrel {\mkern 3mu\vcenter {\hbox {$\shortmid $}}\mkern -10mu{\to }} d\) in \(\mathbb {D}\). These systems are acted on by loose morphisms in \(\mathbb {D}\).\par{}The second of these features allows us to take the monoidal product of systems in the module of systems theory that is the restriction of this niche. In particular, given two systems  \(x \colon  1 \mathrel {\mkern 3mu\vcenter {\hbox {$\shortmid $}}\mkern -10mu{\to }} d\) and \(x' \colon  1 \mathrel {\mkern 3mu\vcenter {\hbox {$\shortmid $}}\mkern -10mu{\to }} d'\), their monoidal product \emph{in \(\mathbb {D}\)} is the loose morphism \(x \otimes  x' \colon  1 \otimes  1 \mathrel {\mkern 3mu\vcenter {\hbox {$\shortmid $}}\mkern -10mu{\to }} d \otimes  d'\). Notably, this is not a system in the systems theory \(\mathsf {Hom}^{l}(\mathbb {D})(1, \mathbb {D})\), because its domain is not \(1\). However, we can apply the conjoint commuter transform of the isomorphism \(1 \otimes  1 \cong  1\) to modify \(x \otimes  x'\) so that its domain is \(1\) and hence is a system in this systems theory.

  \begin{center}
    \begin {tikzcd}
      1 & {1 \otimes  1} & {d \otimes  d'}
      \arrow ["{\cong  }", "\shortmid "{marking}, tail reversed, no head, from=1-1, to=1-2]
      \arrow ["{x \otimes  x'}", "\shortmid "{marking}, from=1-2, to=1-3]
    \end {tikzcd}
  \end{center}\par{}In the remainder of this section, we define a \textbf{doctrine of initial processes}.  A \emph{theory} in this doctrine will be a symmetric monoidal double category \((\mathbb {D}, \otimes , 1)\) whose canonical isomorphism \(1 \otimes  1 \cong  1\) is a \hyperref[djm-00EX]{conjoint commuter transform} and it will define the module of systems that is the restriction of the niche shown above.\par{}We begin by defining a 2-category of pointed double categories, whose symmetric monoidal objects objects will be the \emph{theories} of the doctrine of initial processes.
\begin{definition}[{Pointed double category}]\label{djm-00GK}
\par{}	A \textbf{pointed double category} is a double category \(\mathbb {D}\) equipped with an object \(d_{\bullet } \in  \mathbb {D}\). Equivalently, it is a double functor \(d_{\bullet } : \bullet  \to  \mathbb {D}\) from the terminal double category.\end{definition}

\begin{definition}[{The cartesian 2-category of pointed double functors}]\label{ssl-003V}
\par{}	The \textbf{2-category \(\mathcal {D}\mathsf {bl}_{\bullet }\) of pointed double categories} and double functors which preserve the point up to a \hyperref[djm-00EY]{conjoint tight isomorphism} is defined to be the following pullback:
	
\begin{center}
		\begin {tikzcd}
			\mathcal {D}\mathsf {bl}_{\bullet } & {{\mathsf {2}\mathcal {C}\mathsf {at}^{\mathsf {colax}}(\Delta [1], \mathcal {D}\mathsf {bl})_{\mathsf {conj}}}} \\
			\mathcal {D}\mathsf {bl} & {\mathcal {D}\mathsf {bl} \times  \mathcal {D}\mathsf {bl}}
			\arrow [from=1-1, to=1-2]
			\arrow [from=1-1, to=2-1]
			\arrow ["{(d_1^*, d_0^*)}", from=1-2, to=2-2]
			\arrow [""{name=0, anchor=center, inner sep=0}, "{(\bullet , \mathrm {id})}"', from=2-1, to=2-2]
			\arrow ["\lrcorner "{anchor=center, pos=0.125}, draw=none, from=1-1, to=0]
		\end {tikzcd}
	\end{center}

	which, as a pullback of cartesian 2-functors, is a cartesian 2-category.\end{definition}

\begin{lemma}[{Symmetric monoidal and cartesian objects of \(\mathcal {D}\mathsf {bl}_{\bullet }\)}]\label{djm-00GM}
\par{}	A symmetric monoidal object of \hyperref[djm-00GK]{the 2-category \(\mathcal {D}\mathsf {bl}_{\bullet }\) of pointed double categories} is equivalent to a symmetric monoidal double category \(\mathbb {D}\) for which the canonical isomorphism \(1 \otimes  1 \cong  1\) is a \hyperref[djm-00EX]{conjoint commuter cell}.\par{}	Similarly, a cartesian object of \(\mathcal {D}\mathsf {bl}_{\bullet }\) is a cartesian double category for which the unique isomorphism \(1 \times  1 \to  1\) is a \hyperref[djm-00EX]{conjoint commuter cell}.\begin{proof}\par{}	First, note that the point of a symmetric monoidal pointed double category is its identity object.\par{}	The remainder follows by \emph{symmetry of internalization} (Theorem 7.4 of \cite{arkor-2024-enhanced}): forgetting for a moment the conjoint commuter condition in \(\mathcal {D}\mathsf {bl}_{\bullet }\), we would have that \(\mathcal {S}\mathsf {M}(\mathcal {D}\mathsf {bl}_{\bullet }) \cong  \ast  \downarrow ^{\mathsf {ps}} \mathcal {S}\mathsf {M}(\mathcal {D}\mathsf {bl})\) is the pseudo-slice under the terminal symmetric monoidal category. But a symmetric monoidal double functor from the terminal symmetric monoidal double category must be (up to isomorphism given by the unitor of that functor) be an inclusion of the monoidal unit; since \(\mathcal {D}\mathsf {bl}_{\bullet }\) does in fact require the conjoint commuter condition, we see that the laxator of the inclusion of the monoidal unit must be a conjoint commuter.\par{}	A similar argument applies in the cartesian case.\end{proof}\end{lemma}
\par{}	We can now describe the doctrine of initial processes. First, by defining \emph{right niches} and their restriction to looser right modules and then giving a doctrine that factors through this restriction.
\begin{definition}[{The 2-category of right niches}]\label{ssl-0040}
\par{}A \textbf{right niche} is a niche whose left leg is a pseudo double functor out of the terminal double category:
  
\begin{center}
    \begin {tikzcd}
      \bullet  & {\mathbb {E}_1} \\
      {M_0} & {M_1}
      \arrow ["{F_0}"', from=1-1, to=2-1]
      \arrow ["{F_1}", from=1-2, to=2-2]
      \arrow ["M"', "\shortmid "{marking}, from=2-1, to=2-2]
    \end {tikzcd}
  \end{center}\par{}There is a 2-category of right niches \(\mathcal {N}\mathsf {iche}_r\) which is the pullback of the following 2-categories. 
  
\begin{center}
  \begin {tikzcd}
    {\mathcal {N}\mathsf {iche}_r} & \mathcal {N}\mathsf {iche} & {\ell }\mathcal {B}\mathsf {imod} \\
    {\mathcal {D}\mathsf {bl}_{\bullet } \times  \mathsf {2}\mathcal {C}\mathsf {at}^{\mathsf {colax}}(\Delta [1], \mathcal {D}\mathsf {bl})_{\mathsf {comp}}} & {\mathsf {2}\mathcal {C}\mathsf {at}^{\mathsf {colax}}(\Delta [1], \mathcal {D}\mathsf {bl})_{\mathsf {conj}} \times  \mathsf {2}\mathcal {C}\mathsf {at}^{\mathsf {colax}}(\Delta [1], \mathcal {D}\mathsf {bl})_{\mathsf {comp}}} & {\mathcal {D}\mathsf {bl} \times  \mathcal {D}\mathsf {bl}}
    \arrow [from=1-1, to=1-2]
    \arrow [from=1-1, to=2-1]
    \arrow [from=1-2, to=1-3]
    \arrow [from=1-2, to=2-2]
    \arrow ["{(()_0, ()_1)}", from=1-3, to=2-3]
    \arrow [""{name=0, anchor=center, inner sep=0}, from=2-1, to=2-2]
    \arrow [""{name=1, anchor=center, inner sep=0}, "{d_0^* \times  d_0^*}"', from=2-2, to=2-3]
    \arrow ["\lrcorner "{anchor=center, pos=0.125}, draw=none, from=1-1, to=0]
    \arrow ["\lrcorner "{anchor=center, pos=0.125}, draw=none, from=1-2, to=1]
  \end {tikzcd}
  \end{center}\par{}Since \(\mathcal {N}\mathsf {iche}_r\) is the pullback of cartesian 2-categories along cartesian 2-functors, it is cartesian as is the restriction \[\mathsf {Res} \colon  \mathcal {N}\mathsf {iche}_r \to  {\ell }\mathcal {M}\mathsf {od}_{\mathsf {r}}.\]\end{definition}

\begin{proposition}[{Restriction of right niches}]\label{ssl-004H}
\par{}There is a cartesian 2-functor \(\mathsf {Res} \colon  \mathcal {N}\mathsf {iche} \to  {\ell }\mathcal {M}\mathsf {od}_{\mathsf {r}}\) that makes the following diagram commute:
  
\begin{center}
    \begin {tikzcd}
      {\mathcal {N}\mathsf {iche}_r} & {\ell }\mathcal {M}\mathsf {od}_{\mathsf {r}} \\
      \mathcal {N}\mathsf {iche} & {\ell }\mathcal {B}\mathsf {imod}
      \arrow ["\mathsf {Res}", from=1-1, to=1-2]
      \arrow [from=1-1, to=2-1]
      \arrow [from=1-2, to=2-2]
      \arrow ["\mathsf {Res}"', from=2-1, to=2-2]
    \end {tikzcd}
  \end{center}

  On objects this implies that the restriction of a right niche is a loose right module.\begin{proof}\par{}The outer diagram in the following commutes by \hyperref[ssl-004G]{Proposition \ref{ssl-004G}}. 
    
\begin{center}
      \begin {tikzcd}
        {\mathcal {N}\mathsf {iche}_r} & \mathcal {N}\mathsf {iche} \\
        {\mathcal {D}\mathsf {bl}_{\bullet } \times  \mathsf {2}\mathcal {C}\mathsf {at}^{\mathsf {colax}}(\Delta [1], \mathcal {D}\mathsf {bl})_{\mathsf {comp}}} & {\ell }\mathcal {M}\mathsf {od}_{\mathsf {r}} & {\ell }\mathcal {B}\mathsf {imod} \\
        & \mathcal {D}\mathsf {bl} & {\mathcal {D}\mathsf {bl} \times  \mathcal {D}\mathsf {bl}}
        \arrow [from=1-1, to=1-2]
        \arrow [from=1-1, to=2-1]
        \arrow ["\mathsf {Res}"{description}, dashed, from=1-1, to=2-2]
        \arrow ["\mathsf {Res}", from=1-2, to=2-3]
        \arrow ["{ d_1^* \circ  \pi _2}", from=2-1, to=3-2]
        \arrow [from=2-2, to=2-3]
        \arrow [from=2-2, to=3-2]
        \arrow ["\lrcorner "{anchor=center, pos=0.125}, draw=none, from=2-2, to=3-3]
        \arrow ["{(()_0, ()_1)}", from=2-3, to=3-3]
        \arrow ["{(\bullet , \mathrm {id})}"', from=3-2, to=3-3]
      \end {tikzcd}
    \end{center}

    So \(\mathsf {Res} \colon  \mathcal {N}\mathsf {iche} \to  {\ell }\mathcal {M}\mathsf {od}_{\mathsf {r}}\) is induced by the universal property of pullback. Furthermore, as the pullback of cartesian 2-functors, it is cartesian as well.\end{proof}\end{proposition}

\begin{proposition}[{The right niche of a pointed double category}]\label{ssl-004I}
\par{}There exists a cartesian 2-functor \(p \colon  \mathcal {D}\mathsf {bl}_{\bullet } \to  \mathcal {N}\mathsf {iche}_r\) that takes a pointed double category \(d_\bullet  \colon  \bullet  \to   \mathbb {D}\) to the right niche 
  
\begin{center}
    \begin {tikzcd}
      \bullet  & {\mathbb {D}} \\
      {\mathbb {D}} & {\mathbb {D}}
      \arrow ["{d_\bullet }"', from=1-1, to=2-1]
      \arrow ["\mathrm {id}", from=1-2, to=2-2]
      \arrow ["{\mathsf {Hom}^{l}(\mathbb {D})}"', from=2-1, to=2-2]
    \end {tikzcd}
  \end{center}\begin{proof}\par{}Let \(\pi  \colon  \mathcal {D}\mathsf {bl}_{\bullet } \to  \mathcal {D}\mathsf {bl}\) be the functor that forgets the point of a pointed double category. Note that on objects, \(s_0^* \colon  \mathcal {D}\mathsf {bl} \to  \mathsf {2}\mathcal {C}\mathsf {at}^{\mathsf {colax}}(\Delta [1], \mathcal {D}\mathsf {bl})_{\mathsf {comp}}\) takes a double category to the identity double functor on it. Then, the following diagram of cartesian 2-functors commutes and induces the desired cartesian 2-functor.

  \begin{center}
      \begin {tikzcd}
        \mathcal {D}\mathsf {bl}_{\bullet } & \mathcal {D}\mathsf {bl} \\
        & {\mathcal {N}\mathsf {iche}_r} & {\ell }\mathcal {B}\mathsf {imod} \\
        & {\mathcal {D}\mathsf {bl}_{\bullet } \times  \mathsf {2}\mathcal {C}\mathsf {at}^{\mathsf {colax}}(\Delta [1], \mathcal {D}\mathsf {bl})_{\mathsf {comp}}} & {\mathcal {D}\mathsf {bl} \times  \mathcal {D}\mathsf {bl}}
        \arrow ["\pi ", from=1-1, to=1-2]
        \arrow [dashed, from=1-1, to=2-2]
        \arrow ["{(\mathrm {id}, s_0^*\pi )}"', from=1-1, to=3-2]
        \arrow ["\mathsf {Hom}^{l}", from=1-2, to=2-3]
        \arrow [from=2-2, to=2-3]
        \arrow [from=2-2, to=3-2]
        \arrow ["{(()_0, ()_1)}", from=2-3, to=3-3]
        \arrow [""{name=0, anchor=center, inner sep=0}, "{\pi  \times  d_0^*}"', from=3-2, to=3-3]
        \arrow ["\lrcorner "{anchor=center, pos=0.125}, draw=none, from=2-2, to=0]
      \end {tikzcd}
    \end{center}\end{proof}\end{proposition}

\begin{definition}[{Doctrine of initial processes}]\label{djm-00GL}
\par{}	The \textbf{doctrine of initial processes} is given by the following commuting square of cartesian 2-functors:
	
\begin{center}
		\begin {tikzcd}
			\mathcal {D}\mathsf {bl}_{\bullet } & {\mathcal {N}\mathsf {iche}_r} & {\ell }\mathcal {M}\mathsf {od}_{\mathsf {r}} \\
			\mathcal {D}\mathsf {bl} && \mathcal {D}\mathsf {bl}
			\arrow ["p", from=1-1, to=1-2]
			\arrow [from=1-1, to=2-1]
			\arrow ["{\mathsf {Res}}", from=1-2, to=1-3]
			\arrow [from=1-3, to=2-3]
			\arrow [equals, from=2-1, to=2-3]
		\end {tikzcd}
	\end{center}

	where \(p \colon  \mathcal {D}\mathsf {bl}_{\bullet } \to  \mathcal {N}\mathsf {iche}_r\) is defined in \hyperref[ssl-004I]{Proposition \ref{ssl-004I}}.\end{definition}
\par{}Applying \(\mathcal {S}\mathsf {M}\) to this doctrine, we see that \hyperref[djm-00GM]{symmetric monoindal objects of \(\mathcal {D}\mathsf {bl}_{\bullet }\)} are theories that produce modules of systems. Finally, we verify that these theories are indeed what we claimed at the beginning of this section.
\begin{explication}[{Doctrine of initial processes}]\label{ssl-004K}
\par{}Explicitly, the \hyperref[djm-00GL]{doctrine of initial processes} sends a pointed double category \(d_{\bullet } : \bullet  \to  \mathbb {D}\) to the loose bimodule \(\mathsf {Hom}^{l}(\mathbb {D})(d_{\bullet }, \mathbb {D}) : \bullet  \mathrel {\mkern 3mu\vcenter {\hbox {$\shortmid $}}\mkern -10mu{\to }} \mathbb {D}\) of its \hyperref[djm-005F]{loose Hom bimodule}. In particular, the objects of its \hyperref[ssl-0036]{carrier} are loose morphisms \(d_{\bullet } \mathrel {\mkern 3mu\vcenter {\hbox {$\shortmid $}}\mkern -10mu{\to }} d\). The loose morphisms of its \hyperref[djm-004A]{target} are loose morphisms of \(\mathbb {D}\) which act by composition.\par{}Recall, that systems theories are symmetric monoidal objects of the doctrine's domain. By \hyperref[djm-00GM]{Lemma \ref{djm-00GM}} the symmetric monoidal objects of \(\mathcal {D}\mathsf {bl}_{\bullet }\) are, essentially, just the symmetric monoidal double categories (satisfying a mild condition) pointed by their unit.\par{}	Therefore, when we start with a symmetric monoidal double category \(\mathbb {P}\) (which we refer to as a \emph{process theory}) and pass it through the doctrine of initial processes, we end up with a module of systems over \(\mathbb {P}\) whose systems are  loose homs out of its monoidal unit, which we call \emph{initial processes}.\end{explication}
\section{Span and cospan doctrines via adequate triples}\label{ssl-0041}\par{}	In this section, we will further specify the \hyperref[ssl-003Z]{doctrine of initial processes} includes many of our examples, we now turn to more specific doctrines which have a more uniform flavor. The remaining constructions will all use the same underlying categorical technology: \emph{spans}.\par{}	Some systems theories --- those that compose by \emph{sharing variables} like the Hamiltonian and Langranian mechanics of \cite{baez-2021-open} (see our similar construction, following 1.3.2 of \cite{schreiber-20xx-differential} in \hyperref[djm-00FU]{Example \ref{djm-00FU}}) and the Schultz-Spivak-Vasilakopolou approach to Willems' style behavioral control theory through sheaves (\emph{machines} in \cite{schultz-2019-dynamical}) --- are naturally described by spans and pullback. But many others --- those which compose by \emph{gluing} together subsystems such as Petri nets \cite{baez-2020-open} and stock-flow diagrams \cite{baez-2022-categorical} --- are better handled by cospans, as well explored by the decorated and structured cospan literature (see, e.g. \cite{courser-2020-open}). Cospans are, however, just spans in the opposite category.\par{}	We begin by defining a doctrine of spans of an adequate triple, which factors through the \hyperref[djm-00GL]{doctrine of initial processes}. We then define doctrines of spans of lex categories and cospans of rex categories, which factor through the doctrine of spans of adequate triples. Respectively, these correspond to the doctrines of variables sharing and port-plugging.\subsection{Doctrine of spans of adequate triples}\label{ssl-004V}\par{}	We begin by defining the cartesian category of \hyperref[djm-00GN]{pointed adequate triples}. Symmetric monoidal objects in this category will be theories in the doctrine of spans of adequate triples.
\begin{definition}[{Pointed adequate triple}]\label{djm-00GN}
\par{}	A \textbf{pointed adequate triple} consists of an \hyperref[djm-0087]{adequate triple} \((\mathsf {C}, (L, R))\) equipped with an object \(c \in  \mathsf {C}\) thought of as a map from the terminal adequate triple \(c : (\bullet , (\mathsf {all}, \mathsf {all})) \to  (\mathsf {C}, (L, R))\).\par{}	The 2-category \(\mathcal {A}\mathsf {dTr}_{\bullet }\) of pointed adequate triples consists of maps of adequate triples which preserve the point up to an isomorphism in the left class. That is, we may define \(\mathcal {A}\mathsf {dTr}_{\bullet }\) as the full sub-2-category of the pseudo-slice \(\bullet  \downarrow ^{\mathsf {ps}} \mathcal {A}\mathsf {dTr}\) under the terminal adequate triple, itself constructed by the following pullback:
	
\begin{center}
		\begin {tikzcd}
			\bullet  \downarrow ^{\mathsf {ps}} \mathcal {A}\mathsf {dTr} & {\mathsf {2}\mathcal {C}\mathsf {at}^{\mathsf {ps}}(\Delta [1], \mathcal {A}\mathsf {dTr})} \\
			\mathcal {A}\mathsf {dTr} & {\mathcal {A}\mathsf {dTr} \times  \mathcal {A}\mathsf {dTr}}
			\arrow [from=1-1, to=1-2]
			\arrow [from=1-1, to=2-1]
			\arrow ["{(d_1^{\ast }, d_0^{\ast })}", from=1-2, to=2-2]
			\arrow [""{name=0, anchor=center, inner sep=0}, "{(\bullet , \mathrm {id})}"', from=2-1, to=2-2]
			\arrow ["\lrcorner "{anchor=center, pos=0.125}, draw=none, from=1-1, to=0]
		\end {tikzcd}
	\end{center}

	consisting of those 1-cells 
	
\begin{center}
		\begin {tikzcd}
	\bullet  & \bullet  \\
	\mathsf {C} & \mathsf {D}
	\arrow [equals, from=1-1, to=1-2]
	\arrow ["c"', from=1-1, to=2-1]
	\arrow ["d", from=1-2, to=2-2]
	\arrow ["\ell ", shorten <=4pt, shorten >=4pt, Rightarrow, from=2-1, to=1-2]
	\arrow ["f"', from=2-1, to=2-2]
\end {tikzcd}
	\end{center}

	whose colaxator isomorphism \(\ell \) is in the left class.\end{definition}

\begin{observation}[{Extending the span construction to pointed adequate triples}]\label{ssl-004W}
\par{}The \hyperref[djm-0088]{span} cartesian 2-functor \(\mathbb {S}\mathsf {pan} : \mathcal {A}\mathsf {dTr} \to  \mathcal {D}\mathsf {bl}\) extends to a cartesian 2-functor \(\mathbb {S}\mathsf {pan} : \mathcal {A}\mathsf {dTr} \to  \mathcal {D}\mathsf {bl}\) which maps the adequate triple \((\mathsf {C}, (L, R))\) with point \(c \in  \mathsf {C}\) to the double category \(\mathbb {S}\mathsf {pan}(\mathsf {C}, (L, R))\) with point \(c\) as well.\par{}We note that \(\mathbb {S}\mathsf {pan}_{\bullet }\) does indeed land in \(\mathcal {D}\mathsf {bl}_{\bullet }\) since we assumed that the isomorphisms witnessing preservation of the pointings were in the left class, and these are all conjoints in \(\mathbb {S}\mathsf {pan}(\mathsf {C}, (L, R))\).\end{observation}
\par{}	Finally, we may then restrict the  to \hyperref[djm-00GN]{pointed adequate triples} along the \hyperref[djm-0088]{span construction} to get a doctrine of spans in an adequate triple.
\begin{definition}[{Doctrine of spans in a pointed adequate triple}]\label{djm-00GO}
\par{}	By restricting the , we get the following doctrine of spans in a pointed adequate triple.
	
\begin{center}
		\begin {tikzcd}
			\mathcal {A}\mathsf {dTr}_{\bullet } & \mathcal {D}\mathsf {bl}_{\bullet } & {\mathcal {N}\mathsf {iche}_r} & {\ell }\mathcal {M}\mathsf {od}_{\mathsf {r}} \\
			\mathcal {A}\mathsf {dTr} & \mathcal {D}\mathsf {bl} && \mathcal {D}\mathsf {bl}
			\arrow ["{\mathbb {S}\mathsf {pan}_{\bullet }}", from=1-1, to=1-2]
			\arrow [from=1-1, to=2-1]
			\arrow ["p", from=1-2, to=1-3]
			\arrow [from=1-2, to=2-2]
			\arrow ["\mathsf {Res}", from=1-3, to=1-4]
			\arrow ["()_1", from=1-4, to=2-4]
			\arrow ["{\mathbb {S}\mathsf {pan}}"', from=2-1, to=2-2]
			\arrow [equals, from=2-2, to=2-4]
		\end {tikzcd}
	\end{center}\end{definition}

\begin{explication}[{Doctrine of spans in a pointed adequate triple}]\label{ssl-004X}
\par{}In the \hyperref[djm-00GO]{doctrine of spans in a pointed adequate triple}, a theory is a symmetric monoidal pointed adequate triple. As in \hyperref[djm-00GM]{the case of pointed double categories}, a symmetric monoidal object of \(\mathcal {A}\mathsf {dTr}_{\bullet }\) consists of a \hyperref[djm-00AQ]{symmetric monoidal adequate triple} \((\mathsf {C}, (L, R))\) whose point is the monoidal object \(1\) of \(\mathsf {C}\).\par{}This theory produces the module of systems \[\mathbb {S}\mathsf {pan}(\mathsf {C}, (L,R))(1, -) \colon  \bullet  \mathrel {\mkern 3mu\vcenter {\hbox {$\shortmid $}}\mkern -10mu{\to }} \mathbb {S}\mathsf {pan}(\mathsf {C}, (L, R)).\] This module of systems has:
  \begin{itemize}\item{}An interface is an object of \(\mathsf {C}\).
    \item{}A system with interface \(c\) is a span \(1 \xleftarrow {\ell  \in  L} x \xrightarrow {r \in  R} c\).
    \item{}An interaction is a span in the adequate triple which acts on systems by pullback.
    \item{}Maps of systems and maps of interactions are maps of spans.\end{itemize}\end{explication}
\subsection{Span and cospan doctrines for lex and rex categories}\label{ssl-0051}\par{}	We now give two specific examples of this doctrine: the \emph{variable sharing} doctrine of spans in a \hyperref[djm-001Z]{lex category} and the \emph{port-plugging} doctrine of cospans in a \hyperref[djm-00AY]{rex category}.\par{}We begin by defining pointed lex and rex categories, whose symmetric monoidal objects will be \emph{theories} in these doctrines.
\begin{lemma}[{The pointed adequate triple of a lex (resp. rex) category}]\label{ssl-0050}
\par{}The cartesian 2-functor \(\mathcal {L}\mathsf {ex} \to  \mathcal {A}\mathsf {dTr}\) defined in \hyperref[djm-008C]{Construction \ref{djm-008C}} extends to a cartesian 2-functor \(\mathcal {L}\mathsf {ex} \to  \mathcal {A}\mathsf {dTr}_{\bullet }\) which on objects sends a lex category \(\mathsf {C}\) to the adquate triple \((\mathsf {C}, (\mathsf {all}, \mathsf {all}))\) pointed by the initial object \(1 \in  \mathsf {C}\).\par{}Likewise, the cartesian 2-functor \({\mathcal {R}\mathsf {ex}}^{\mathsf {co}} \to  \mathcal {A}\mathsf {dTr}\) defined in \hyperref[djm-008D]{Construction \ref{djm-008D}} extends to a cartesian 2-functor \({\mathcal {R}\mathsf {ex}}^{\mathsf {co}}  \to  \mathcal {A}\mathsf {dTr}_{\bullet }\) which on objects sends a rex category \(\mathsf {C}\) to the  adquate triple \((\mathsf {C}, (\mathsf {all}, \mathsf {all}))\) pointed by the terminal object \(0 \in  \mathsf {C}\).\begin{proof}\par{}This follows from the definition of \hyperref[djm-00GN]{pointed adequate triples} and the fact that lex categories preserve initial objects and hence the lex functor from the terminal lex category to a lex category \(\mathsf {C}\) picks out the initial object of \(\mathsf {C}\). Likewise for rex functors from the terminal rex category to a rex category \(\mathsf {C}\) picks out the terminal object of \(\mathsf {C}\).\end{proof}\end{lemma}

\begin{observation}[{Symmetric monoidal pointed lex and rex categories}]\label{ssl-0052}
\par{}We note that \(\mathcal {L}\mathsf {ex}\) has biproducts, so that every lex category is uniquely symmetric monoidal with the symmetric monoidal product given by cartesian product.\par{}Likewise, \(\mathcal {R}\mathsf {ex}\) has biproducts and so every rex category is uniquely symmetric monoidal with the symmetric monoidal product given by coproduct.\end{observation}
\par{}	We can restrict the \hyperref[djm-00GO]{adequate triple doctrine} along the inclusion \(\mathcal {L}\mathsf {ex} \to  \mathcal {A}\mathsf {dTr}\) of \hyperref[djm-008C]{Construction \ref{djm-008C}}. We will call the resulting doctrine the \emph{variable sharing} doctrine.
\begin{definition}[{Variable sharing doctrine}]\label{djm-00GP}
\par{}	We define the \textbf{variable sharing} \hyperref[djm-00A8]{doctrine} to be the restriction of the \hyperref[djm-00GO]{doctrine of spans of adequate triples} to \hyperref[djm-008C]{adequate triples induced by lex categories}, we get the following doctrine:
	
\begin{center}
		\begin {tikzcd}
      \mathcal {L}\mathsf {ex} & \mathcal {A}\mathsf {dTr}_{\bullet } & \mathcal {D}\mathsf {bl}_{\bullet } & \mathcal {N}\mathsf {iche} & {\ell }\mathcal {M}\mathsf {od}_{\mathsf {r}} \\
      \mathcal {L}\mathsf {ex} & \mathcal {A}\mathsf {dTr} & \mathcal {D}\mathsf {bl} && \mathcal {D}\mathsf {bl}
      \arrow [from=1-1, to=1-2]
      \arrow [equals, from=1-1, to=2-1]
      \arrow ["{\mathbb {S}\mathsf {pan}_{\bullet }}", from=1-2, to=1-3]
      \arrow [from=1-2, to=2-2]
      \arrow ["p", from=1-3, to=1-4]
      \arrow [from=1-3, to=2-3]
      \arrow ["\mathsf {Res}", from=1-4, to=1-5]
      \arrow ["_1", from=1-5, to=2-5]
      \arrow [from=2-1, to=2-2]
      \arrow ["{\mathbb {S}\mathsf {pan}}"', from=2-2, to=2-3]
      \arrow [equals, from=2-3, to=2-5]
    \end {tikzcd}
	\end{center}\end{definition}

\begin{definition}[{Port-plugging doctrine}]\label{djm-00GQ}
\par{}	We define the \textbf{port-plugging} \hyperref[djm-00A8]{doctrine} to be the restriction of the \hyperref[djm-00GO]{span doctrine of a pointed adequate triple} to \hyperref[djm-008D]{adequate triples induced by rex categories}:
	
\begin{center}
    \begin {tikzcd}
      {\mathcal {R}\mathsf {ex}^{\mathsf {co}}} & \mathcal {L}\mathsf {ex} & \mathcal {A}\mathsf {dTr}_{\bullet } & \mathcal {D}\mathsf {bl}_{\bullet } & \mathcal {N}\mathsf {iche} & {\ell }\mathcal {M}\mathsf {od}_{\mathsf {r}} \\
      {\mathcal {R}\mathsf {ex}^{\mathsf {co}}} & \mathcal {L}\mathsf {ex} & \mathcal {A}\mathsf {dTr} & \mathcal {D}\mathsf {bl} && \mathcal {D}\mathsf {bl}
      \arrow ["^{\mathsf {op}}", from=1-1, to=1-2]
      \arrow [equals, from=1-1, to=2-1]
      \arrow [from=1-2, to=1-3]
      \arrow [equals, from=1-2, to=2-2]
      \arrow ["{\mathbb {S}\mathsf {pan}_{\bullet }}", from=1-3, to=1-4]
      \arrow [from=1-3, to=2-3]
      \arrow ["p", from=1-4, to=1-5]
      \arrow [from=1-4, to=2-4]
      \arrow ["\mathsf {Res}", from=1-5, to=1-6]
      \arrow ["{_1}"{description}, from=1-6, to=2-6]
      \arrow ["^{\mathsf {op}}", from=2-1, to=2-2]
      \arrow [from=2-2, to=2-3]
      \arrow ["{\mathbb {S}\mathsf {pan}}"', from=2-3, to=2-4]
      \arrow [equals, from=2-4, to=2-6]
    \end {tikzcd}
		\end{center}

	This sends a rex category \(\mathsf {C}\) to the loose right module \[\mathbb {C}\mathsf {ospan}(\mathsf {C})^{\mathsf {op}}(\emptyset , -) \colon  \bullet  \mathrel {\mkern 3mu\vcenter {\hbox {$\shortmid $}}\mkern -10mu{\to }} \mathbb {C}\mathsf {ospan}(\mathsf {C})^{\mathsf {op}}\] of cospans out of the initial object. We can further dualize (in the tight direction) to get \[\mathbb {C}\mathsf {ospan}(\mathsf {C})(\emptyset , -) \colon  \bullet  \mathrel {\mkern 3mu\vcenter {\hbox {$\shortmid $}}\mkern -10mu{\to }} \mathbb {C}\mathsf {ospan}(\mathsf {C})\] as a pseudo-functor of \(\mathsf {C} \in  \mathcal {R}\mathsf {ex}\) (rather than \(\mathcal {R}\mathsf {ex}^{\mathsf {co}}\)).\end{definition}
\section{The doctrine of generalized Moore machines}\label{djm-00G6}\par{}	In \hyperref[ssl-0041]{Section \ref{ssl-0041}}, we saw how systems that compose via \emph{variable sharing} or \emph{port plugging} are defined by theories in  doctrines defined by spans and cospans.\par{}	Other systems --- those which are \emph{generalized Moore machines}, such as systems of ODEs and POMDPs and other "automata" --- compose best via \emph{lenses}, as described for example in \cite{vagner-2014-algebras} (see also \cite{jaz-2021-book} for a pedagogical overview). However, lenses may also be described as spans of a particular sort (see \hyperref[djm-00B6]{Definition \ref{djm-00B6}}).\par{}	Roughly, a generalized Moore machine consists of an internal state, an input, and an output. Outputs are determined by the state and states change according to input.  In the doctrine of generalized Moore machines that we will describe in this section, a theory (i.e. the data needed to define a module of systems) is a \textbf{tangency}, which consists of:
	\begin{itemize}\item{}A notion of \emph{space}.
		\item{}A notion of \emph{bundle} over spaces, which may be pulled back.
		\item{}An assignment of a "tangent" bundle to each space.\end{itemize}
	A tangency is the data necessary to construct a module of systems for which systems are \emph{generalized Moore machines} and interactions are \emph{lenses} (in the generalized sense of Spivak \cite{spivak-2019-generalized}).\par{}	With this data in place, we can now sketch the definition a generalized Moore machine following the discussion in \cite{jaz-2021-double}. A generalized Moore machine consists of:
	\begin{itemize}\item{}			A way things may be in the form of a space \(S\) of \emph{states}.
		\item{}			A way things may change in the form of an element \(TS\) in the fiber over \(S\). In many examples of a tangency, such an element is in fact a bundle \(TS \to  S\). Then, for a state \(s \in  S\), the fiber \(T_s S\) represents possible changes from \(s\).
		\item{}			An \emph{interface} in the form of:
			\begin{itemize}\item{}						A space \(O\) of \emph{observations} that may be made of the system or \emph{orientations} which the system may take in its environment.
				\item{}						An element \(I\) in the fiber over \(O\) of \emph{inputs}. In many examples of a tangency, such an element is in fact a bundle \(I \to  O\). For an orientation \(o \in  O\), and the space \(I_o\) represents \emph{parameters} to the dynamics of the system which are available in  that orientation.					\end{itemize}
		\item{}			An \emph{observation} \(e \colon  S \to  O\) of the system, which exposes some aspects of state.
		\item{}			A \emph{(parameterized) dynamics} \(u \colon  e^{\ast }I \to  TS\). In many examples of a tangency, these dynamics assign a state \(s \in  S\) and a valid parameter \(p \in  I_{e(s)}\), to a change \(u_s(p) \in  T_s S\).\end{itemize}\subsection{Tangencies}\label{ssl-004L}\par{}In this section, we give a definition of symmetric monoidal tangencies, which are  the theories in the doctrine of generalized Moore machines.
\begin{definition}[{The 2-category of tangencies}]\label{djm-00BK}
\par{}	A \textbf{tangency} is a \hyperref[djm-008F]{cartesian fibration} \(\pi  : E \to  B\) equipped with a section \(T : B \to  E\) (with \(\pi  T = \mathrm {id}_B\)).\par{}		The 2-category \(\mathcal {T}\mathsf {an}\) of tangencies is defined as the following pullback of 2-categories:
	
\begin{center}
		\begin {tikzcd}[column sep = large]
			{\mathcal {T}\mathsf {an}} && {\mathsf {2}\mathcal {C}\mathsf {at}^{\mathsf {colax}}(\Delta [2], \mathcal {C}\mathsf {at})} \\
			\mathcal {F}\mathsf {ib} & {2\mathcal {C}\mathsf {at}(\Delta [1], \mathcal {C}\mathsf {at}) \times  \mathcal {C}\mathsf {at}} & {\mathsf {2}\mathcal {C}\mathsf {at}^{\mathsf {colax}}(\Delta [1], \mathcal {C}\mathsf {at}) \times  \mathsf {2}\mathcal {C}\mathsf {at}^{\mathsf {colax}}(\Delta [1], \mathcal {C}\mathsf {at})}
			\arrow [from=1-1, to=1-3]
			\arrow [from=1-1, to=2-1]
			\arrow ["{(d_0^{\ast }, d_1^{\ast })}", from=1-3, to=2-3]
			\arrow [""{name=0, anchor=center, inner sep=0}, "{(U, d_0^{\ast })}"', from=2-1, to=2-2]
			\arrow ["{(\iota , s_0^*)}"', from=2-2, to=2-3]
			\arrow ["\lrcorner "{anchor=center, pos=0.125}, draw=none, from=1-1, to=0]
		\end {tikzcd}
	\end{center}

	where \(U \colon  \mathcal {F}\mathsf {ib} \to  2\mathcal {C}\mathsf {at}(\Delta [1], \mathcal {C}\mathsf {at})\) forgets the structure of the fibration, \(\iota  \colon  2\mathcal {C}\mathsf {at}(\Delta [1], \mathcal {C}\mathsf {at}) \to  \mathsf {2}\mathcal {C}\mathsf {at}^{\mathsf {colax}}(\Delta [1], \mathcal {D}\mathsf {bl})\) is the inclusion, and \(s_0^* \colon  \mathcal {C}\mathsf {at} \to  \mathsf {2}\mathcal {C}\mathsf {at}^{\mathsf {colax}}(\Delta [1], \mathcal {D}\mathsf {bl})\) takes a category to its identity functor.\end{definition}

\begin{explication}[{The 2-category of tangencies}]\label{ssl-004N}
\par{}Here we explicitly state the objects, morphisms, and 2-cells of \hyperref[djm-00BK]{\(\mathcal {T}\mathsf {an}\)}.\par{}
                  \begin{itemize}\item{}An object is a tangency. In other words, a pair \((\pi  \colon  E \to  B, T \colon  B \to  E)\) where \(\pi \) is a cartesian fibration and \(T\) is a section of \(\pi \).
    \item{}A morphism of tangencies \[((f, \overline {f}), \phi ) \colon   (\pi _1, T_1) \to  (\pi _2, T_2)\] consists of a \hyperref[djm-008F]{cartesian functor} \((f, \overline {f}) : \pi _1 \to  \pi _2\) together with a colax morphism of sections \(\phi  : \overline {f} \circ  T_1 \Rightarrow  T_2 \circ  f\) so that \(\pi _2 \phi  = \mathrm {id}_{f}\):
      
\begin{center}
        \begin {tikzcd}
          {B_1} & {B_2} && {B_1} & {B_2} \\
          {E_1} & {E_2} & {=} \\
          {B_1} & {B_2} && {B_1} & {B_2}
          \arrow ["f", from=1-1, to=1-2]
          \arrow ["{T_1}"', from=1-1, to=2-1]
          \arrow ["{T_2}", from=1-2, to=2-2]
          \arrow ["f", from=1-4, to=1-5]
          \arrow [Rightarrow, no head, from=1-4, to=3-4]
          \arrow [Rightarrow, no head, from=1-5, to=3-5]
          \arrow ["\phi ", shorten <=2pt, shorten >=2pt, Rightarrow, from=2-1, to=1-2]
          \arrow ["{\overline {f}}"', from=2-1, to=2-2]
          \arrow ["{\pi _1}"', from=2-1, to=3-1]
          \arrow ["{\pi _2}", from=2-2, to=3-2]
          \arrow ["f"', from=3-1, to=3-2]
          \arrow ["f"', from=3-4, to=3-5]
        \end {tikzcd}
      \end{center}

      These compose by whiskering:
      
\begin{center}
        \begin {tikzcd}
          {B_1} & {B_3} && {B_1} & {B_2} & {B_3} \\
          {E_1} & {E_3} & {:=} & {E_1} & {E_2} & {E_3} \\
          {B_1} & {B_3} && {B_1} & {B_2} & {B_3}
          \arrow ["{g\circ  f}", from=1-1, to=1-2]
          \arrow ["{T_1}"', from=1-1, to=2-1]
          \arrow ["{T_3}", from=1-2, to=2-2]
          \arrow ["f", from=1-4, to=1-5]
          \arrow ["{T_1}"', from=1-4, to=2-4]
          \arrow ["g", from=1-5, to=1-6]
          \arrow ["{T_2}"{description}, from=1-5, to=2-5]
          \arrow ["{T_3}", from=1-6, to=2-6]
          \arrow ["{\phi  ; \psi }"{description}, shorten <=2pt, shorten >=2pt, Rightarrow, from=2-1, to=1-2]
          \arrow ["{\overline {g}\circ  \overline {f}}"', from=2-1, to=2-2]
          \arrow [from=2-1, to=3-1]
          \arrow [from=2-2, to=3-2]
          \arrow ["\phi "{description}, shorten <=2pt, shorten >=2pt, Rightarrow, from=2-4, to=1-5]
          \arrow ["{\overline {f}}"', from=2-4, to=2-5]
          \arrow ["{\pi _1}"', from=2-4, to=3-4]
          \arrow ["\psi "{description}, shorten <=2pt, shorten >=2pt, Rightarrow, from=2-5, to=1-6]
          \arrow ["{\overline {g}}"', from=2-5, to=2-6]
          \arrow ["{\pi _2}", from=2-5, to=3-5]
          \arrow ["{\pi _3}", from=2-6, to=3-6]
          \arrow ["{g\circ  f}"', from=3-1, to=3-2]
          \arrow ["f"', from=3-4, to=3-5]
          \arrow ["g"', from=3-5, to=3-6]
        \end {tikzcd}
      \end{center}
    \item{}A 2-cell of tangency morphisms is a 2-cell \((\alpha , \overline {\alpha }) : (f, \overline {f}) \Rightarrow  (g, \overline {g})\) of \(\mathcal {F}\mathsf {ib}\) for which
        
\begin{center}
          \begin {tikzcd}
            {B_1} & {B_2} && {B_1} & {B_2} \\
            {B_1} & {B_2} & {=} & {E_1} & {E_2} \\
            {E_1} & {E_2} && {E_1} & {E_2}
            \arrow ["g", from=1-1, to=1-2]
            \arrow [equals, from=1-1, to=2-1]
            \arrow [equals, from=1-2, to=2-2]
            \arrow ["g", from=1-4, to=1-5]
            \arrow ["{T_1}"', from=1-4, to=2-4]
            \arrow ["{T_2}", from=1-5, to=2-5]
            \arrow ["\alpha "{description}, shorten <=2pt, shorten >=2pt, Rightarrow, from=2-1, to=1-2]
            \arrow ["f", from=2-1, to=2-2]
            \arrow ["{T_1}"', from=2-1, to=3-1]
            \arrow ["{T_2}", from=2-2, to=3-2]
            \arrow ["\psi "{description}, shorten <=2pt, shorten >=2pt, Rightarrow, from=2-4, to=1-5]
            \arrow ["{\overline {g}}", from=2-4, to=2-5]
            \arrow [equals, from=2-4, to=3-4]
            \arrow [equals, from=2-5, to=3-5]
            \arrow ["\phi "{description}, shorten <=2pt, shorten >=2pt, Rightarrow, from=3-1, to=2-2]
            \arrow ["{\overline {f}}"', from=3-1, to=3-2]
            \arrow ["{\overline {\alpha }}"{description}, shorten <=2pt, shorten >=2pt, Rightarrow, from=3-4, to=2-5]
            \arrow ["{\overline {f}}"', from=3-4, to=3-5]
          \end {tikzcd}
        \end{center}\end{itemize}
                \end{explication}

\begin{attitude}[{Tangencies as states and possible changes}]\label{djm-00F0}
\par{}	A \hyperref[djm-00BK]{tangency} is a cartesian fibration \(\pi  : E \to  B\) equipped with a section \(T : B \to  E\). We think of a tangency as follows:
	\begin{enumerate}\item{}\(B\) is a category of \emph{state spaces}.
		\item{}			For every state space \(S \in  B\), the fiber of \(\pi  : E \to  B\) over \(S\) is a category \(E_S\)  of \emph{bundles} over \(S\), which can be pulled-back along maps of spaces.
		\item{}			The section \(T : B \to  E\) assigns to every state space its "tangent bundle", or bundle of possible changes.\end{enumerate}

	A \emph{system} in a tangency is then a "generalized Moore machine", or a \hyperref[djm-00B6]{\emph{lens}} \({u \choose  E} \colon  {TS \choose  S} \mathrel {\mkern 3mu\vcenter {\hbox {$\shortmid $}}\mkern -10mu{\to }} {I \choose  O}\) which corresponds to the following span:
	
\begin{center}
		\begin {tikzcd}
			TS & {e^{\ast }I} & I & E \\
			S & S & O & B
			\arrow [dashed, maps to, from=1-1, to=2-1]
			\arrow ["u"', from=1-2, to=1-1]
			\arrow ["{\mathsf {lift}(e)}", from=1-2, to=1-3]
			\arrow [dashed, maps to, from=1-2, to=2-2]
			\arrow ["\lrcorner "{anchor=center, pos=0.125}, draw=none, from=1-2, to=2-3]
			\arrow [dashed, maps to, from=1-3, to=2-3]
			\arrow ["\pi ", from=1-4, to=2-4]
			\arrow [equals, from=2-1, to=2-2]
			\arrow ["e"', from=2-2, to=2-3]
		\end {tikzcd}
	\end{center}

	We think of this as consisting of a morphism \(e : S \to  O\) in \(B\) which \emph{exposes} some variables of state, and a parameterized update \(u : e^{\ast } I \to  TS\). Using \hyperref[ssl-004S]{lens notation}, we notate such a system a lens \[{u \choose  e} \colon  {TS \choose  S} \mathrel {\mkern 3mu\vcenter {\hbox {$\shortmid $}}\mkern -10mu{\to }} {I \choose  O}.\]\par{}	We refer the reader to \cite{jaz-2021-double} for further discussion.\end{attitude}

\begin{attitude}[{Morphisms of tangencies as changing the state space}]\label{ssl-004R}
\par{}Let \((\pi _1, T_1)\) and \((\pi _2, T_2)\) be tangencies and let \[((f, \overline  f) \colon  \pi _1 \to  \pi _2, \phi  \colon  \overline  f \circ  T_1 \Rightarrow  T_2 \circ  f)\] be \hyperref[ssl-004N]{a morphism of tangencies}. 
  We interpret this morphism as transforming the state space \(S\) in the tangency \((\pi _1, T_1)\) into the state space \(f(S)\) in the tangency \((\pi _2, T_2)\). Simply applying \(\overline {f}\) to the bundle of changes \(T_1S\) over \(S\) is \emph{not} the bundle of changes  \(T_2(fS)\) over \(f(S)\). However, applying \(\phi _S\) to \(\overline {f}(T_1S)\) maps it into the bundle of changes \(T_2(fS)\) over \(f(S)\).\par{}We can interpret this transformation as a lens. For \(S \in  B_1\), the morphism \(\phi _S\) in \(E_2\) defines the following lens in \(\mathbb {L}\mathsf {ens}(E_2)\)

  \begin{center}
    \begin {tikzcd}
      {T_2(fS)} & {\overline {f}(T_1S)} & {\overline {f}(T_1S)} & {E_2} \\
      {f(S)} & {f(S)} & {f(S)} & {B_2}
      \arrow [dashed, maps to, from=1-1, to=2-1]
      \arrow ["{\phi _S}"', from=1-2, to=1-1]
      \arrow [equals, from=1-2, to=1-3]
      \arrow [dashed, maps to, from=1-2, to=2-2]
      \arrow [dashed, maps to, from=1-3, to=2-3]
      \arrow ["{\phi _2}", from=1-4, to=2-4]
      \arrow [equals, from=2-1, to=2-2]
      \arrow [equals, from=2-2, to=2-3]
    \end {tikzcd}
  \end{center}\par{}Using \hyperref[ssl-004S]{lens notation}, \(\phi \) induces lenses \[{\phi _S \choose  \mathrm {id}} \colon  {T_2(fS) \choose  f(S)} \mathrel {\mkern 3mu\vcenter {\hbox {$\shortmid $}}\mkern -10mu{\to }} {\overline {f}(T_1S) \choose  f(S)}.\]\par{}See Proposition 4.5.1.11 of the manuscript \cite{jaz-2021-book} to see the Euler method as an example of a morphism of tangencies from a tangency of ODEs to a tangency for discrete-time, continuous-space Moore machines. We expect that the Runge-Kutta method could also be expressed as such morphism, allowing us to derive the compositionality theorem proved by Ngotiaoco \cite{ngotiaoco-2017-compositionality} via the pseudo-functoriality of the tangency \hyperref[djm-00A8]{doctrine}.\end{attitude}

\begin{remark}[{Tangencies in \cite{jaz-2021-book}}]\label{djm-00BP}
\par{}	The definition of a tangency in the context of categorical systems theory first appears as a "dynamical systems doctrine" in Definition 1.1 of \cite{jaz-2021-double}. It appears further in \cite{jaz-2021-book} as a "theory of dynamical systems" Definition 3.5.0.4 and Definition 4.5.1.2 respectively.\par{}		Here, we have changed the terminology to "tangency" to accomodate the other theories of dynamical systems we find in other doctrines. We also finally add the symmetric monoidal structure necessary for the parallel product of systems; see \hyperref[djm-00BJ]{Lemma \ref{djm-00BJ}}.\end{remark}
\par{}Finally, we show that \(\mathcal {T}\mathsf {an}\) is cartesian and characterize its symmetric monoidal objects. These are the theories in the doctrine of generalized Moore machines.
\begin{lemma}[{Tangencies form a cartesian 2-category}]\label{djm-00BL}
\par{}The 2-category \(\mathcal {T}\mathsf {an}\) of \hyperref[djm-00BK]{tangencies} form a cartesian 2-category with products constructed in \(\mathcal {F}\mathsf {ib}\) (and therefore ultimately in \(\mathcal {C}\mathsf {at}\)).\begin{proof}\par{}	The product of \((\pi _1, T_1)\) and \((\pi _2, T_2)\) is \((\pi _1 \times  \pi _2, T_1 \times  T_2)\). We equip the	projections with the identity transformations \(p_i(T_1 \times  T_2) = T_ip_i\). It is straightforward to verify the 2-dimensional universal property of the product, since the requisite equations may be checked componentwise.\end{proof}\end{lemma}

\begin{lemma}[{Symmetric monoidal tangency}]\label{djm-00BJ}
\par{}By \emph{symmetry of internalization} (Theorem 7.4 of \cite{arkor-2024-enhanced}), a symmetric monoidal tangency \((\pi , T) \in  \mathcal {S}\mathsf {M}(\mathcal {T}\mathsf {an})\) is equivalently a symmetric monoidal fibration \(\pi  : E \to  B\) equipped with a lax monoidal section \(T : B \to  E\).\end{lemma}
\subsection{Tangencies as systems theories}\label{ssl-004M}\par{}	In this section, we construct a loose bimodule from a tangency. We then  this loose bimodule to form a module of systems.
\begin{lemma}[{From tangencies to adequate triples}]\label{djm-00EL}
\par{}	Let \((\pi  \colon  E \to  B, T \colon  B \to  E)\) be a \hyperref[djm-00BK]{tangency}. The section \(T\) induces a map of \hyperref[djm-0087]{adequate triples} by 
	\[
		T \colon  (B, (\mathsf {id}, \mathsf {id})) \to  (E, (\mathsf {vert}, \mathsf {cart}))
	\]
	from the  to the \hyperref[djm-008E]{adequate triple associated to the fibration \(\pi \)}. This construction gives a cartesian 2-functor
	\[
		\mathsf {T} : \mathcal {T}\mathsf {an} \to  \mathsf {2}\mathcal {C}\mathsf {at}^{\mathsf {colax}}(\Delta [1], \mathcal {A}\mathsf {dTr}).
 	\]
	into maps of adequate triples and colax transformations.\begin{proof}\par{}	The underlying 2-functor into \(\mathcal {C}\mathsf {at}\) is just a projection:
	\[
		\mathcal {T}\mathsf {an} \to  \mathsf {2}\mathcal {C}\mathsf {at}^{\mathsf {colax}}(\Delta [2], \mathcal {C}\mathsf {at}) \xrightarrow {d_2^{\ast }} \mathsf {2}\mathcal {C}\mathsf {at}^{\mathsf {colax}}(\Delta [1], \mathcal {C}\mathsf {at})	
	\]
		It only remains to show that this lifts along the forgetful functor \(\mathsf {2}\mathcal {C}\mathsf {at}^{\mathsf {colax}}(\Delta [1], \mathcal {A}\mathsf {dTr}) \to  \mathsf {2}\mathcal {C}\mathsf {at}^{\mathsf {colax}}(\Delta [1], \mathcal {C}\mathsf {at})\). Since \(T : B \to  E\) gives a map of adequate triples \(T : (B, (\mathsf {id}, \mathsf {id})) \to  (E, (\mathsf {vert}, \mathsf {cart}))\), and since both the assignments \(B \mapsto  (B (\mathsf {id}, \mathsf {id}))\) and \(\pi  \mapsto  (E, (\mathsf {vert}, \mathsf {cart}))\) are 2-functorial, and since the forgetful functor \(\mathcal {A}\mathsf {dTr} \to  \mathcal {C}\mathsf {at}\) is locally fully faithful, we have demonstrated that the above 2-functor lifts to \(\mathsf {2}\mathcal {C}\mathsf {at}^{\mathsf {colax}}(\Delta [1], \mathcal {A}\mathsf {dTr})\).\par{}	Since the projection is evidently cartesian, and since \hyperref[djm-00AK]{products in \(\mathcal {A}\mathsf {dTr}\) are constructed in \(\mathcal {C}\mathsf {at}\)}, the 2-functor \(\mathsf {T}\) is cartesian.\end{proof}\end{lemma}

\begin{lemma}[{Span construction of tangency lands in commuter tranformations}]\label{djm-00F3}
\par{}	The composite 2-functor
	
\begin{center}
		\begin {tikzcd}[column sep = huge]
			{\mathcal {T}\mathsf {an}} & {\mathsf {2}\mathcal {C}\mathsf {at}^{\mathsf {colax}}(\Delta [1], \mathcal {A}\mathsf {dTr})} & {\mathsf {2}\mathcal {C}\mathsf {at}^{\mathsf {colax}}(\Delta [1], \mathcal {D}\mathsf {bl})}
			\arrow ["{\mathsf {T}}", from=1-1, to=1-2]
			\arrow ["{\mathsf {2}\mathcal {C}\mathsf {at}^{\mathsf {colax}}(\Delta [1], \mathbb {S}\mathsf {pan})}", from=1-2, to=1-3]
		\end {tikzcd}
	\end{center}

	lands in \(\mathsf {2}\mathcal {C}\mathsf {at}^{\mathsf {colax}}(\Delta [1], \mathcal {D}\mathsf {bl})_{\mathsf {conj}}\). We will refer to the resulting composite as 
	\[
		\mathbb {S}\mathsf {pan}(\mathsf {T}) : \mathcal {T}\mathsf {an} \to  \mathsf {2}\mathcal {C}\mathsf {at}^{\mathsf {colax}}(\Delta [1], \mathcal {D}\mathsf {bl})_{\mathsf {conj}}.
	\]\begin{proof}\par{}	We need to show that for any 1-cell \(\phi  : T_2f \Rightarrow  \overline {f} T_1 : \pi _1 \to  \pi _2\) of \hyperref[djm-00BK]{tangencies}, the resulting tight transformation \(\mathbb {S}\mathsf {pan}(\phi )\) is a \hyperref[djm-00EX]{conjoint commuter transformation}. Note that by definition, \(\phi \) is vertical; it is therefore in the left class of the adequate triple \((E_2, (\mathsf {vert}, \mathsf {cart}))\) and is \hyperref[djm-00F4]{therefore component-wise a conjoint} in \(\mathbb {S}\mathsf {pan}(E_2, (\mathsf {vert}, \mathsf {cart}))\). It remains to show that it is a commuter; but the only loose morphisms in \(\mathbb {S}\mathsf {pan}(B_1, (\mathsf {id}, \mathsf {id}))\) are identities, and the transpose of the loose identity of a conjoint is always a tight isomorphism.\end{proof}\end{lemma}

\begin{explication}[{From tangencies to spans}]\label{ssl-004P}
\par{}Let's examine the action of the functor \(\mathbb {S}\mathsf {pan}(\mathsf {T}) : \mathcal {T}\mathsf {an} \to  \mathsf {2}\mathcal {C}\mathsf {at}^{\mathsf {colax}}(\Delta [1], \mathcal {D}\mathsf {bl})_{\mathsf {conj}}\) on objects.\par{}Let \((\pi  \colon  E \to  B, T \colon  B \to  E)\) be a tangency. We perform the \hyperref[djm-0088]{span construction} to get a map \[\mathbb {S}\mathsf {pan}(T) : \mathbb {S}\mathsf {pan}(B, (\mathsf {id}, \mathsf {id})) \to  \mathbb {S}\mathsf {pan}(E, (\mathsf {vert}, \mathsf {cart})).\] Note that: 
  \begin{itemize}\item{}By \hyperref[djm-00EM]{Observation \ref{djm-00EM}}, \(\mathbb {S}\mathsf {pan}(B,( \mathsf {id}, \mathsf {id}))\) is the loosely discrete double category \(\mathbb {T}\mathsf {}(B)\) on \(B\) which has \(B\) as its tight category and only identity loose arrows.
    \item{}As per \hyperref[djm-00B6]{Definition \ref{djm-00B6}} (see also \hyperref[djm-00B7]{Explication \ref{djm-00B7}}), \(\mathbb {S}\mathsf {pan}(E, (\mathsf {vert}, \mathsf {cart}))\) is the double category \(\mathbb {L}\mathsf {ens}(E)\) of \hyperref[djm-00B6]{Spivak lenses}.\end{itemize} 
   
   Therefore, \(\mathbb {S}\mathsf {pan}(T)\) is given just by the action of \(T\) in the tight categories; it sends a "base space" \(b \in  B\) to its "tangent bundle" \(Tb \in  E\).\par{}To simplify terminology, let's refer to \(\mathbb {S}\mathsf {pan}(T)\) by \(T\) so that under \(\mathbb {S}\mathsf {pan}(\mathsf {T})\), the tangency \((\pi , T)\) maps to the double functor:
			\[T : \mathbb {T}\mathsf {}(B) \to  \mathbb {L}\mathsf {ens}(E).\]\end{explication}

\begin{lemma}[{From tangencies to niches}]\label{ssl-004O}
\par{}There is a cartesian 2-functor \(\mathcal {T}\mathsf {an} \to  \mathcal {N}\mathsf {iche}\) that takes a tangency \((\pi  \colon  E \to  B, T \colon  B \to  E)\) to the niche

  \begin{center}
    \begin {tikzcd}[column sep = large]
      {\mathbb {T}\mathsf {}(B)} & {\mathbb {L}\mathsf {ens}(E)} \\
      {\mathbb {L}\mathsf {ens}(E)} & {\mathbb {L}\mathsf {ens}(E)}
      \arrow ["T"', from=1-1, to=2-1]
      \arrow [equals, from=1-2, to=2-2]
      \arrow ["{ \mathbb {L}\mathsf {ens}(E) (-,-)}"', "\shortmid "{marking}, from=2-1, to=2-2]
    \end {tikzcd}
  \end{center}\begin{proof}\par{}The composite \[\mathcal {T}\mathsf {an}  \xrightarrow {\mathbb {S}\mathsf {pan}(\mathsf {T})} \mathsf {2}\mathcal {C}\mathsf {at}^{\mathsf {colax}}(\Delta [1], \mathcal {D}\mathsf {bl})_{\mathsf {conj}} \xrightarrow {d_0^*} \mathcal {D}\mathsf {bl}\] 
    takes a tangency with total space \(E\) to the double category \(\mathbb {L}\mathsf {ens}(E)\) of \hyperref[djm-00B6]{Spviak lenses}.\par{}The outer diagram in the following commutes.

  \begin{center}
      \begin {tikzcd}
        \mathcal {T}\mathsf {an} & \mathcal {D}\mathsf {bl} \\
        \mathsf {2}\mathcal {C}\mathsf {at}^{\mathsf {colax}}(\Delta [1], \mathcal {D}\mathsf {bl})_{\mathsf {conj}} & {\mathcal {N}\mathsf {iche}} & {\ell }\mathcal {B}\mathsf {imod} \\
        & {\mathsf {2}\mathcal {C}\mathsf {at}^{\mathsf {colax}}(\Delta [1], \mathcal {D}\mathsf {bl})_{\mathsf {conj}} \times  \mathsf {2}\mathcal {C}\mathsf {at}^{\mathsf {colax}}(\Delta [1], \mathcal {D}\mathsf {bl})_{\mathsf {comp}}} & {\mathcal {D}\mathsf {bl} \times  \mathcal {D}\mathsf {bl}}
        \arrow ["{d_0^* \circ  \mathbb {S}\mathsf {pan}(\mathsf {T})}", from=1-1, to=1-2]
        \arrow ["{\mathbb {S}\mathsf {pan}(\mathsf {T})}"', from=1-1, to=2-1]
        \arrow [dashed, from=1-1, to=2-2]
        \arrow ["\mathsf {Hom}^{l}", from=1-2, to=2-3]
        \arrow ["{(\mathrm {id}, s_0^*d_0^*)}"', from=2-1, to=3-2]
        \arrow [from=2-2, to=2-3]
        \arrow [from=2-2, to=3-2]
        \arrow ["{()_0, () _1)}", from=2-3, to=3-3]
        \arrow [""{name=0, anchor=center, inner sep=0}, "{()_1 \times () _1}"', from=3-2, to=3-3]
        \arrow [""{name=0p, anchor=center, inner sep=0}, phantom, from=3-2, to=3-3, start anchor=center, end anchor=center]
        \arrow ["\lrcorner "{anchor=center, pos=0.125, rotate=45}, draw=none, from=2-2, to=0p]
      \end {tikzcd}
    \end{center}

    Hence by the universal property of pullbacks, we have the desired map \(\mathcal {T}\mathsf {an} \to  \mathcal {N}\mathsf {iche}\).\end{proof}\end{lemma}

\begin{explication}[{From tangencies to loose bimodules}]\label{djm-00F5}
\par{}	Let \(\mathsf {Sys} \colon  \mathcal {T}\mathsf {an} \to  {\ell }\mathcal {B}\mathsf {imod}\) be the composite of the cartesian 2-functor \(\mathcal {T}\mathsf {an} \to  \mathcal {N}\mathsf {iche}\) defined in \hyperref[ssl-004O]{Lemma \ref{ssl-004O}} with the  \(\mathsf {Res} \colon  \mathcal {N}\mathsf {iche} \to  {\ell }\mathcal {B}\mathsf {imod}\) turns a tangency \((\pi  \colon  E \to  B , T \colon  B \to  E)\) into the loose bimodule \[\mathbb {L}\mathsf {ens}(T-, -) \colon  \mathbb {T}\mathsf {}(B) \mathrel {\mkern 3mu\vcenter {\hbox {$\shortmid $}}\mkern -10mu{\to }} \mathbb {L}\mathsf {ens}(E).\]\par{}	We see that a heteromorphism with domain \(S \in  B\) and codomain \({I \choose  O} \in  \mathbb {L}\mathsf {ens}(E)\) of the resulting bimodule is a \hyperref[djm-00B6]{lens} \({u \choose  e} \colon  {TS \choose  S} \mathrel {\mkern 3mu\vcenter {\hbox {$\shortmid $}}\mkern -10mu{\to }} {I \choose  O}\). Explicitly, such a lens consists of \(e \colon  S \to  O\)  in \(B\) and \(u \colon  e^*(I) \to  TS\) in \(E\) and these form the span
	
\begin{center}
		\begin {tikzcd}
			TS & {e^{\ast }I} & I \\
			S & S & {O}
			\arrow [dashed, maps to, from=1-1, to=2-1]
			\arrow ["u"', from=1-2, to=1-1]
			\arrow ["{\mathsf {lift}(e)}", from=1-2, to=1-3]
			\arrow [dashed, maps to, from=1-2, to=2-2]
			\arrow ["\lrcorner "{anchor=center, pos=0.125}, draw=none, from=1-2, to=2-3]
			\arrow [dashed, maps to, from=1-3, to=2-3]
			\arrow [equals, from=2-1, to=2-2]
			\arrow ["e"', from=2-2, to=2-3]
		\end {tikzcd}
		\end{center}

		in \(E\).
		Such a morphism is precisely a \emph{system} in the tangency as described in \hyperref[djm-00F0]{Attitude \ref{djm-00F0}}.\par{}	It is also worth understanding the functoriality of this construction. Let \((\pi _1, T_1)\) and \((\pi _2, T_2)\) be tangencies and let \[((f, \overline  f) \colon  \pi _1 \to  \pi _2, \phi  \colon  \overline  f \circ  T_1 \Rightarrow  T_2 \circ  f)\] be \hyperref[ssl-004N]{a morphism of tangencies}. Recall the \hyperref[ssl-004R]{attitude for morphisms of tangencies} The morphism on tangencies induces a (covariant) action on systems which for \(S \in  B\) and \(I \in  E\) sends the system \({u \choose  e} \colon  {T_1S \choose  S} \mathrel {\mkern 3mu\vcenter {\hbox {$\shortmid $}}\mkern -10mu{\to }} {I \choose  O}\) in \(\mathbb {L}\mathsf {ens}(E_1)\) to the composite 
	
\begin{center}
		\begin {tikzcd}
			{{T_2(fS) \choose  fS}} & {{\overline {f}T_1S \choose  fS}} & {{\overline {f}I \choose  fO}}
			\arrow ["{{\phi _S \choose  \mathrm {id}}}", "\shortmid "{marking}, from=1-1, to=1-2]
			\arrow ["{{\overline {f}u \choose  fe}}", "\shortmid "{marking}, from=1-2, to=1-3]
		\end {tikzcd}
	\end{center}
 
	in \(\mathbb {L}\mathsf {ens}(E_2)\). This lens corresponds to the following span:
	
\begin{center}
		\begin {tikzcd}
			{T_2fS} & {\overline {f}T_1S} & {\overline {f}(e^{\ast }I)} & {\overline {f}I} & {E_2} \\
			fS & fS & fS & fO & {B_2}
			\arrow [dashed, maps to, from=1-1, to=2-1]
			\arrow ["{\phi _S}"', from=1-2, to=1-1]
			\arrow [dashed, maps to, from=1-2, to=2-2]
			\arrow ["{\overline {f}u}"', from=1-3, to=1-2]
			\arrow ["{\overline {f}\mathsf {lift}(e)}", from=1-3, to=1-4]
			\arrow [dashed, maps to, from=1-3, to=2-3]
			\arrow ["\lrcorner "{anchor=center, pos=0.125}, draw=none, from=1-3, to=2-4]
			\arrow [dashed, maps to, from=1-4, to=2-4]
			\arrow ["{\pi _2}", from=1-5, to=2-5]
			\arrow [equals, from=2-1, to=2-2]
			\arrow [equals, from=2-2, to=2-3]
			\arrow ["fe"', from=2-3, to=2-4]
		\end {tikzcd}
	\end{center}\end{explication}
\par{}Finally, we are ready to define the doctrine of generalized Moore machines.
\begin{definition}[{Doctrine of generalized Moore machines}]\label{djm-00HA}
\par{}	We define the \hyperref[djm-00A8]{doctrine} of \textbf{generalized Moore machines} as follows:
	
\begin{center}
		\begin {tikzcd}
			{\mathcal {T}\mathsf {an}} & {\ell }\mathcal {M}\mathsf {od}_{\mathsf {r}} \\
			\mathcal {F}\mathsf {ib} & \mathcal {D}\mathsf {bl}
			\arrow ["{\mathbb {M}\mathsf {oore}}", from=1-1, to=1-2]
			\arrow [from=1-1, to=2-1]
			\arrow [from=1-2, to=2-2]
			\arrow ["{\mathbb {L}\mathsf {ens}}"', from=2-1, to=2-2]
		\end {tikzcd}
	\end{center}

	where \(\mathbb {M}\mathsf {oore}\) is the composite of the \hyperref[djm-00J2]{collapse of a loose bimodule into a right bimodule} of  \(\mathsf {Sys}\), the \hyperref[djm-00F5]{loose bimodule of systems associated to a tangency}. The 2-functor \(\mathbb {L}\mathsf {ens} : \mathcal {F}\mathsf {ib} \to  \mathcal {D}\mathsf {bl}\) is the \hyperref[djm-00B6]{lens construction}.\end{definition}

\begin{explication}[{Doctrine of generalized Moore machines}]\label{ssl-004T}
\par{}In the \hyperref[djm-00HA]{doctrine of generalized Moore machines}, a systems theory is a \hyperref[djm-00BJ]{symmetric monoidal tangency} \((\pi  \colon  E \to  B, T \colon  B \to  E)\). It produces the module of systems \(\mathbb {M}\mathsf {oore}(\pi , T)\) over the double category of interactions \(\mathbb {L}\mathsf {ens}(E)\).\par{}
                  \begin{itemize}\item{}An interface in this module of systems is an object \({I \choose  O}\) in \(\mathbb {L}\mathsf {ens}(E)\).
    \item{}A system with  interface \({I \choose  O}\) consists of a state space \(S \in  B\) and a lens \({u \choose  e} \colon  {TS \choose  S} \mathrel {\mkern 3mu\vcenter {\hbox {$\shortmid $}}\mkern -10mu{\to }} {I \choose  O}\) in which \(e\) represents an observation and \(u\) represents parameterized dynamics.
    \item{}An interaction is a lens that acts on systems via lens composition.
    \item{}Maps of systems and maps of interactions are given by charts.\end{itemize}
                \par{}	Next we investigate the monoidal product of systems in this doctrine. For \(i = 1,2\), let \({u_i \choose  r_i} \colon  {TS_i \choose  S_i} \mathrel {\mkern 3mu\vcenter {\hbox {$\shortmid $}}\mkern -10mu{\to }} {I_i \choose  O_i}\) be systems in \(\mathbb {M}\mathsf {oore}(\pi , T)\). Since \((\pi , T)\) is a symmetric monoidal object of \(\mathcal {T}\mathsf {an}\), by \emph{symmetry of internalization} \(T \colon  B \to  E\) has laxator \(\mu _{S_1, S_2} : TS_1 \otimes  TS_2 \Rightarrow  T(S_1 \otimes  S_2)\) which defines a lens \[{\mu _{S_1, S_2} \choose  \mathrm {id} } \colon  {T(S_1 \otimes  S_2) \choose  S_1 \otimes  S_2} \mathrel {\mkern 3mu\vcenter {\hbox {$\shortmid $}}\mkern -10mu{\to }} {TS_1 \otimes  TS_2 \choose  S_1 \otimes  S_2}.\] The monoidal product of the systems \({u_1 \choose  r_1}\) and \({u_2 \choose  r_2}\) is the composite of lenses
  
\begin{center}
    \begin {tikzcd}
      {{T(S_1 \otimes  S_2) \choose  S_1 \otimes  S_2}} & {{TS_1 \otimes  TS_2 \choose  S_1 \otimes  S_2}} & {{I_1 \otimes  I_2 \choose  O_1 \otimes  O_2}}
      \arrow ["{{\mu _{S_1, S_2} \choose  \mathrm {id}}}", "\shortmid "{marking}, from=1-1, to=1-2]
      \arrow ["{{u_1 \otimes  u_2 \choose  e_1 \otimes  e_2}}", "\shortmid "{marking}, from=1-2, to=1-3]
    \end {tikzcd}  
  \end{center}

  which corresponds to the following span

  \begin{center}
		\begin {tikzcd}[column sep=large]
      {T(S_1 \otimes  S_2)} & {TS_1 \otimes  TS_2} & {e_1^{\ast }I_1 \otimes  e_2^{\ast }I_2} & {I_1 \otimes  I_2} \\
      {S_1 \otimes  S_2} & {S_1 \otimes  S_2} & {S_1 \otimes  S_2} & {O_1 \otimes  O_2}
      \arrow [dashed, maps to, from=1-1, to=2-1]
      \arrow ["{\mu _{S_1,S_2}}"', from=1-2, to=1-1]
      \arrow [dashed, maps to, from=1-2, to=2-2]
      \arrow ["{u_1 \otimes  u_2}"', from=1-3, to=1-2]
      \arrow ["{\mathsf {lift}(e_1) \otimes  \mathsf {lift}(e_2)}", from=1-3, to=1-4]
      \arrow [dashed, maps to, from=1-3, to=2-3]
      \arrow ["\lrcorner "{anchor=center, pos=0.125}, draw=none, from=1-3, to=2-4]
      \arrow [dashed, maps to, from=1-4, to=2-4]
      \arrow [equals, from=2-1, to=2-2]
      \arrow [equals, from=2-2, to=2-3]
      \arrow ["{e_1 \otimes  e_2}"', from=2-3, to=2-4]
    \end {tikzcd}
  \end{center}\end{explication}
\par{}We conclude this section with an example of a tangency and the module of systems that it defines.
\begin{example}[{Tangent bundle tangency}]\label{djm-00II}
\par{}	A \textbf{submersion} \(p \colon  M \to  N\) of manifolds  is a smooth map whose pushforward \(T_x p : T_x M \to  T_{p(x)}N\) on tangent spaces is surjective. Submersions stable under pullback, and so equip the category \(\mathsf {Man}\) of manifolds with the structure of a \hyperref[djm-00E0]{display map category}. Given a submersion \(p \colon  M \to  N\), the fibers \(M_y := p^{-1}(y)\) are also manifolds. There is an \hyperref[djm-00E1]{associated fibration} is \(\pi  : \mathsf {Subm} \to  \mathsf {Man}\) which sends a submersion \(p \colon  M \to  N\) to the manifold \(N\). It has a section \(T \colon  \mathsf {Man} \to  \mathsf {Subm}\) which on objects
	assigns a manifold \(M\) to  its tangent bundle \(TM \to  M\). Together,  \((\pi , T)\) is a \hyperref[djm-00BK]{tangency}.\par{}	A system in the resulting module of systems then consists of a diagram of manifolds like so:
	
\begin{center}
		\begin {tikzcd}
			TS & {e^{\ast }I} & I \\
			S & S & O
			\arrow [from=1-1, to=2-1]
			\arrow ["u"', from=1-2, to=1-1]
			\arrow [from=1-2, to=1-3]
			\arrow [from=1-2, to=2-2]
			\arrow ["\lrcorner "{anchor=center, pos=0.125}, draw=none, from=1-2, to=2-3]
			\arrow [from=1-3, to=2-3]
			\arrow [equals, from=2-1, to=2-2]
			\arrow ["e"', from=2-2, to=2-3]
		\end {tikzcd}
	\end{center}

	which we interpret as follows
	\begin{enumerate}\item{}			The interface of this system is a submersion \(I \to  O\). The manifold \(O\) is a space of observations. The submersion defines a manifold \(I_o\) for every observation \(o \in  O\).
		\item{}			We have a manifold \(S\) whose points are \emph{states}.
		\item{}			The map \(e : S \to  O\) exposes an observation for every state.
		\item{}			There are smooth maps \(u : I_{e(s)} \to  T_sS\) which are organized into a map of bundles \(u : e^{\ast }I \to  TS\). These maps define a tangent vector \(u(s, i) \in  T_sS\) parameterized by \(i \in  I\). These maps express the system of ordinary differential equations
			\[
				\frac {ds}{dt} = u(s, t)
			\]\end{enumerate}\par{}	These systems interact through lenses which allow us to set the parameters of some systems by exposed observations of other systems.\par{}	Recall that a map out of \hyperref[ssl-002Q]{the timeline Moore machine} is a \hyperref[ssl-002O]{trajectory} for \hyperref[ssl-0029]{deterministic Moore machines}. Similarly, in this module of systems, there is a ‘clock system’ defined by the ODE \(\frac {dt}{dt} = 1\) which is expressed as the following diagram:
		
\begin{center}
		\begin {tikzcd}
	{T\mathbb {R}} & \mathbb {R} & \mathbb {R} \\
	\mathbb {R} & \mathbb {R} & \mathbb {R}
	\arrow [from=1-1, to=2-1]
	\arrow ["1"', from=1-2, to=1-1]
	\arrow [equals, from=1-2, to=1-3]
	\arrow [equals, from=1-2, to=2-2]
	\arrow ["\lrcorner "{anchor=center, pos=0.125}, draw=none, from=1-2, to=2-3]
	\arrow [equals, from=1-3, to=2-3]
	\arrow [equals, from=2-1, to=2-2]
	\arrow ["\mathrm {id}"', equals, from=2-2, to=2-3]
\end {tikzcd}
	\end{center}

	System maps out of this clock system are solutions of the codomain system.\end{example}
\subsection{Doctrine of open coalgebras}\label{djm-00EE}\par{}	In this section, we will define a \hyperref[djm-00A8]{doctrine} of \textbf{open coalgebras}. Where a coalgebra for an endofunctor \(F\) is a map \(d : S \to  FS\), and open coalgebra will be a kind of \hyperref[djm-00HA]{generalized Moore machine}.
\begin{definition}[{Open coalgebras of an endofunctor}]\label{djm-00HD}
\par{}	Let \(\mathsf {C}\) be a cartesian category and let \(F : \mathsf {C} \to  \mathsf {C}\) be an endofunctor. Let \(I\) and  \(O\) be objects of \(\mathsf {C}\). An \textbf{open coalgebra} for \(f\) with interface \({I \choose  O}\) consists of an object \(S \in  C\) together with functions \[
		\left \{\begin {aligned}
			e &: S \to  O \\
			u &: S \times  I \to  FS
		\end {aligned}\right .
	\]
	If \(\mathsf {C}\) is cartesian closed, an open coalgebra for \(F\) with interface \({I \choose  O}\) is equivalently a coalgebra for the endofunctor \(X \mapsto  O \times  (FX)^I\).\end{definition}
\par{}	We will now define a doctrine of open coalgebras.
\begin{definition}[{2-Category of endofunctors on cartesian categories}]\label{djm-00EF}
\par{}	Define the 2-category \(\mathcal {E}\mathsf {nd}^{\mathsf {colax}}(\mathcal {C}\mathsf {art})\) of \textbf{endofunctors on cartesian categories} to be the following pullback of 2-categories:
	
\begin{center}
		\begin {tikzcd}
	{\mathcal {E}\mathsf {nd}^{\mathsf {colax}}(\mathcal {C}\mathsf {art})} & {\mathsf {2}\mathcal {C}\mathsf {at}^{\mathsf {colax}}(\mathsf {Arr}, \mathcal {C}\mathsf {at})} \\
	{\mathcal {C}\mathsf {art}} & {\mathcal {C}\mathsf {at}^2}
	\arrow [from=1-1, to=1-2]
	\arrow [from=1-1, to=2-1]
	\arrow ["\lrcorner "{anchor=center, pos=0.125}, draw=none, from=1-1, to=2-2]
	\arrow ["{(() _0, () _1)}", from=1-2, to=2-2]
	\arrow ["\Delta ",from=2-1, to=2-2]
\end {tikzcd}
	\end{center}

	Explicitly:
	\begin{enumerate}\item{}			The objects of \(\mathcal {E}\mathsf {nd}^{\mathsf {colax}}(\mathcal {C}\mathsf {art})\) are cartesian endofunctors \(F : \mathsf {C} \to  \mathsf {C}\) of cartesian categories \(\mathsf {C}\).
		\item{}			A morphism is a colax square
			
\begin{center}
				\begin {tikzcd}
	\mathsf {C} & \mathsf {D} \\
	\mathsf {C} & \mathsf {D}
	\arrow ["a", from=1-1, to=1-2]
	\arrow ["F"', from=1-1, to=2-1]
	\arrow ["G", from=1-2, to=2-2]
	\arrow ["\alpha ", shorten <=4pt, shorten >=4pt, Rightarrow, from=2-1, to=1-2]
	\arrow ["A"', from=2-1, to=2-2]
\end {tikzcd}
			\end{center}

			where \(A \colon  \mathsf {C} \to  \mathsf {D}\) is a \emph{cartesian} functor.
		\item{}			A 2-cell is a natural transformation \(\varphi  : A \Rightarrow  B\) to that the following equation holds:
			
\begin{center}
				\begin {tikzcd}
	\mathsf {C} & \mathsf {D} && \mathsf {C} & \mathsf {D} \\
	\mathsf {C} & \mathsf {C}\mathsf {D} & {=} & \mathsf {C} & \mathsf {D} \\
	\mathsf {C} & \mathsf {D} && \mathsf {C} & \mathsf {D}
	\arrow ["a", from=1-1, to=1-2]
	\arrow [equals, from=1-1, to=2-1]
	\arrow [equals, from=1-2, to=2-2]
	\arrow ["A", from=1-4, to=1-5]
	\arrow ["F"', from=1-4, to=2-4]
	\arrow ["G", from=1-5, to=2-5]
	\arrow ["\varphi ", shorten <=2pt, shorten >=2pt, Rightarrow, from=2-1, to=1-2]
	\arrow ["b", from=2-1, to=2-2]
	\arrow ["F"', from=2-1, to=3-1]
	\arrow ["G", from=2-2, to=3-2]
	\arrow ["\alpha ", shorten <=2pt, shorten >=2pt, Rightarrow, from=2-4, to=1-5]
	\arrow ["A", from=2-4, to=2-5]
	\arrow [equals, from=2-4, to=3-4]
	\arrow [equals, from=2-5, to=3-5]
	\arrow ["\beta ", shorten <=2pt, shorten >=2pt, Rightarrow, from=3-1, to=2-2]
	\arrow ["B"', from=3-1, to=3-2]
	\arrow ["\varphi ", shorten <=2pt, shorten >=2pt, Rightarrow, from=3-4, to=2-5]
	\arrow ["B"', from=3-4, to=3-5]
\end {tikzcd}
			\end{center}\end{enumerate}\end{definition}

\begin{construction}[{Tangency associated to an endofunctor on a cartesian category}]\label{djm-00EG}
\par{}	Let \(F : \mathsf {C} \to  \mathsf {C}\) be an endofunctor of a cartesian category \(\mathsf {C}\). Define the section \(- \times  F- : \mathsf {C} \to  \mathsf {proj}_\mathsf {C}\) of the \hyperref[djm-00E4]{simple fibration} of \(\mathsf {C}\) by the following composite:
	
\begin{center}
	\begin {tikzcd}
		\mathsf {C} & {\mathsf {C}\times  \mathsf {C}} & {\mathsf {C} \times  \mathsf {C}} & {\mathsf {proj}_\mathsf {C}}
		\arrow ["\Delta ", from=1-1, to=1-2]
		\arrow ["{\mathsf {C} \times  F}", from=1-2, to=1-3]
		\arrow ["\times ", from=1-3, to=1-4]
	\end {tikzcd}
	\end{center}

	where \(\mathsf {C} \times  \mathsf {C} \xrightarrow {\times } \mathsf {proj}_\mathsf {C}\) is the functor which on objects is defined by the projection \[(c, x) \mapsto  (c \times  x \to  c)\] and on morphisms is defined by:  
	
\begin{center}
		\begin {tikzcd}[row sep = tiny]
			&& {c_1 \times  Fc_1} & {c_2 \times  Fc_2} \\
			{\alpha  : c_1 \to  c_2} & \mapsto  \\
			&& {c_1} & {c_2}
			\arrow ["{\alpha  \times  f\alpha }", from=1-3, to=1-4]
			\arrow [from=1-3, to=3-3]
			\arrow [from=1-4, to=3-4]
			\arrow ["\alpha "', from=3-3, to=3-4]
		\end {tikzcd}
	\end{center}\par{}	This gives a cartesian 2-functor \(\mathcal {E}\mathsf {nd}^{\mathsf {colax}}(\mathcal {C}\mathsf {art}) \to  \mathcal {T}\mathsf {an}\).\par{}	Note that it sends a 1-cell \((A, \alpha ) : (\mathsf {C}, F) \to  (\mathsf {D}, G)\) to the following 2-cell:\begin{center}
		\begin {tikzcd}[row sep = large]
			\mathsf {C} & \mathsf {D} \\
			{\mathsf {C}\times  \mathsf {C}} & {\mathsf {D}\times  \mathsf {D}} \\
			{\mathsf {C} \times  \mathsf {C}} & {\mathsf {D} \times  \mathsf {D}} \\
			{\mathsf {proj}_\mathsf {C}} & {\mathsf {proj}_\mathsf {D}}
			\arrow ["A", from=1-1, to=1-2]
			\arrow ["\Delta "', from=1-1, to=2-1]
			\arrow ["\Delta ", from=1-2, to=2-2]
			\arrow ["{A \times  A}", from=2-1, to=2-2]
			\arrow ["{\mathsf {C} \times  F}"', from=2-1, to=3-1]
			\arrow ["{\mathsf {D} \times  G}", from=2-2, to=3-2]
			\arrow ["{\mathrm {id}_A \times  \alpha }", shorten <=4pt, shorten >=4pt, Rightarrow, from=3-1, to=2-2]
			\arrow ["{A \times  A}"', from=3-1, to=3-2]
			\arrow ["\times "', hook, from=3-1, to=4-1]
			\arrow ["\times ", hook, from=3-2, to=4-2]
			\arrow ["{\mathsf {proj}_A}"', from=4-1, to=4-2]
		\end {tikzcd}
	\end{center}\label{djm-00EG-fig}\begin{proof}\par{}	Verification is straightforward: the 1- and 2-cells of \(\mathcal {E}\mathsf {nd}^{\mathsf {colax}}(\mathcal {C}\mathsf {art})\) appear in the middle of \hyperref[djm-00EG-fig]{Figure \ref{djm-00EG-fig}} almost entirely unaltered.\end{proof}\end{construction}
\par{}	In the doctrine of open coalgebras, as systems theory will be a symmetric monoidal object of \(\mathcal {E}\mathsf {nd}^{\mathsf {colax}}(\mathcal {C}\mathsf {art})\). By \emph{symmetry of internalization}, it is straightforward to compute what these are.
\begin{lemma}[{Symmetric monoidal object of \(\mathcal {E}\mathsf {nd}^{\mathsf {colax}}(\mathcal {C}\mathsf {art})\)}]\label{djm-00HE}
\par{}	A symmetric monoidal object of the \hyperref[djm-00EF]{2-category \(\mathcal {E}\mathsf {nd}^{\mathsf {colax}}(\mathcal {C}\mathsf {art})\) of endofunctors on cartesian categories} consists of a lax symmetric monoidal structure on the endofunctor \(F : \mathsf {C} \to  \mathsf {C}\) with respect to the cartesian product of \(\mathsf {C}\).\begin{proof}\par{}	Explicitly, \emph{symmetry of internalization} says that a symmetric monoidal object of \(\mathcal {E}\mathsf {nd}^{\mathsf {colax}}(\mathcal {C}\mathsf {art})\) consists of an endofunctor \(F : \mathsf {C} \to  \mathsf {C}\) on a cartesian category together with a symmetric monoidal structure \(1 : \bullet  \to  \mathsf {C}\) and \(\otimes  :  \times  \mathsf {C} \to  \mathsf {C}\) --- both cartesian functors --- for which \(F\) is lax symmetric monoidal. But if \(1 : \bullet  \to  \mathsf {C}\) is a cartesian functor, it must be the inclusion of the terminal object; and for \(\otimes  : \mathsf {C} \times  \mathsf {C} \to  \mathsf {C}\) to be cartesian means that we have a coherent interchange isomorphism 
	\[(x \times  y) \otimes  (z \times  w) \xrightarrow {\sim } (x \otimes  z) \times  (y \otimes  w).\]
	By the Eckman-Hilton argument, or more generally by the Baez-Dolan stabilization hypothesis (see Corollary 6.2.9 of \cite{gepner-2015-enriched}), this implies that the monoidal product and cartesian product must coincide: \(\otimes  \cong  \times \). Therefore, the data of a symmetric monoidal object in \(\mathcal {E}\mathsf {nd}^{\mathsf {colax}}(\mathcal {C}\mathsf {art})\) consists of a lax monoidal structure on \(F\) with respect to the cartesian product of \(\mathsf {C}\).\end{proof}\end{lemma}

\begin{definition}[{Doctrine of open coalgebras}]\label{djm-00HC}
\par{}	We define the \hyperref[djm-00A8]{doctrine} of \hyperref[djm-00HD]{\textbf{open coalgebras}} to be the restriction of the \hyperref[djm-00HA]{doctrine of generalized Moore machines} along the cartesian 2-functor \(\mathcal {E}\mathsf {nd}^{\mathsf {colax}}(\mathcal {C}\mathsf {art}) \to  \mathcal {T}\mathsf {an}\) of \hyperref[djm-00EG]{Construction \ref{djm-00EG}}:
	
\begin{center}
		\begin {tikzcd}
			{\mathcal {E}\mathsf {nd}^{\mathsf {colax}}(\mathcal {C}\mathsf {art})} & {\mathcal {T}\mathsf {an}} & {\ell }\mathcal {M}\mathsf {od}_{\mathsf {r}} \\
			\mathcal {C}\mathsf {art} & \mathcal {F}\mathsf {ib} & \mathcal {D}\mathsf {bl}
			\arrow [from=1-1, to=1-2]
			\arrow [from=1-1, to=2-1]
			\arrow ["{\mathbb {M}\mathsf {oore}}", from=1-2, to=1-3]
			\arrow [from=1-2, to=2-2]
			\arrow [from=1-3, to=2-3]
			\arrow ["\mathsf {Simp}"', from=2-1, to=2-2]
			\arrow ["{\mathbb {L}\mathsf {ens}}"', from=2-2, to=2-3]
		\end {tikzcd}
	\end{center}

	where \(\mathsf {Simp} : \mathcal {C}\mathsf {art} \to  \mathcal {F}\mathsf {ib}\) is sends a cartesian category \(\mathsf {C}\) to its \hyperref[djm-00E4]{simple fibration} \(\mathsf {proj}_\mathsf {C}\).\end{definition}

\begin{notation}[{Lenses of simple fibrations}]\label{ssl-005Q}
\par{}Let \(\mathsf {C}\) be a cartesian category. An object of \(\mathbb {L}\mathsf {ens}(\mathsf {proj}_\mathsf {C})\) is a pair \({O \times  I \choose  O}\) indicating that \(O \times  I\) lives in the fiber over \(O\) in the simple fibration \(\mathsf {proj}_\mathsf {C}\). We often abbreviate this object to simply \({I \choose  O}\).\end{notation}

\begin{explication}[{Doctrine of open coalgebras}]\label{ssl-005P}
\par{}	From \hyperref[djm-00HE]{Lemma \ref{djm-00HE}}, a systems theory in the doctrine of open coalgebras consists of  an endofunctor \(F : \mathsf {C} \to  \mathsf {C}\) with a lax monoidal structure with respect to the cartesian product of \(\mathsf {C}\). In the module of systems produced by this systems theory 
  \begin{itemize}\item{}Interfaces are pairs of objects in \(\mathsf {C}\) notated \({I \choose  O}\).
    \item{}A system with interface \({I \choose  O}\) is a span in \(\mathsf {proj}_\mathsf {C}\) of the form
      
\begin{center}
        \begin {tikzcd}
          {S \times  FS} & {S \times  I} & {O \times  I} \\
          S & S & O
          \arrow ["{\pi _S}"', from=1-1, to=2-1]
          \arrow ["{(\pi _S, u)}"', from=1-2, to=1-1]
          \arrow ["{e \times  I}", from=1-2, to=1-3]
          \arrow ["{\pi _S}"', from=1-2, to=2-2]
          \arrow ["\lrcorner "{anchor=center, pos=0.125}, draw=none, from=1-2, to=2-3]
          \arrow ["{\pi _O}", from=1-3, to=2-3]
          \arrow [equals, from=2-1, to=2-2]
          \arrow ["e"', from=2-2, to=2-3]
        \end {tikzcd}
      \end{center}

      which consists ultimately of maps \(e : S \to  O\) and \(u : S \times  I \to  FS\). This is precisely an \hyperref[djm-00HD]{open coalgebra}.
  \item{}The parallel product of open coalgebras \({u_1 \choose  e_1} : {FS_1 \choose  S_1} \mathrel {\mkern 3mu\vcenter {\hbox {$\shortmid $}}\mkern -10mu{\to }} {I_1 \choose  O_1}\) and \({u_2 \choose  e_2} : {FS_2 \choose  S_2} \mathrel {\mkern 3mu\vcenter {\hbox {$\shortmid $}}\mkern -10mu{\to }} {I_2 \choose  O_2}\) is given by \(e_1 \times  e_2 : S_1 \times  S_2 \to  O_1 \times  O_2\) together with
    \[
      (S_1 \times  S_1) \times  (I_1 \times  I_2) \cong  (S_1 \times  I_1) \times  (S_2 \times  I_2) \xrightarrow {u_1 \times  u_2} FS_1 \times  FS_2 \xrightarrow {\mu _F} F(S_1 \times  S_2)
    \]\end{itemize}\end{explication}
\par{}	In the remainder of this section we detail examples particular systems theories in \hyperref[djm-00HC]{the doctrine of open coalgebras}.
\begin{example}[{Deterministic Moore machines}]\label{ssl-005R}
\par{}Under \hyperref[djm-00HC]{the doctrine of open coalgebras}, the identity endomorphism on \(\mathsf {Set}\) is a systems theory that defines a module of determinstic Moore machines 
\begin{center}
				\begin {tikzcd}[column sep=8em]
					\bullet  & {\mathbb {L}\mathsf {ens}(\mathsf {proj}_{\mathsf {Set}})}
					\arrow ["{\mathbb {M}\mathsf {oore}(\mathsf {proj}_{\mathsf {Set}}, (-)^2)}", "\shortmid "{marking}, from=1-1, to=1-2]
				\end {tikzcd}
			\end{center}
 
		via the simple fibration \(\mathsf {proj}_{\mathsf {Set}}\) and the tangency \((-)^2 \colon  \mathsf {Set} \to  \mathsf {Set} \times  \mathsf {Set}\). This module of systems  is covered in  detail in \hyperref[ssl-0029]{Section \ref{ssl-0029}}.\end{example}

\begin{example}[{Non-deterministic Moore machines}]\label{ssl-005S}
\par{}Under \hyperref[djm-00HC]{the doctrine of open coalgebras}, the powerset monad on \(\mathsf {Set}\) is a systems theory that defines a module of non-determinstic Moore machines.\end{example}

\begin{example}[{Partially observable Markov decision processes (POMDPs)}]\label{ssl-005T}
\par{}Under \hyperref[djm-00HC]{the doctrine of open coalgebras}, the \emph{Giry monad} \cite{giry-1982-categorical} on measurable spaces --- \(f : \mathsf {Meas} \to  \mathsf {Meas}\) --- is a systems theory that defines a module of \emph{partially observable Markov decision process} (POMDPs).\par{}In this module, we interpret an interface \({I \choose  O}\) as consisting of our set of potential \emph{observations} \(O\) as well as our menu of \emph{actions} \(I\), and interpret a system \({u \choose  e} : {fS \choose  S} \mathrel {\mkern 3mu\vcenter {\hbox {$\shortmid $}}\mkern -10mu{\to }} {I \choose  O}\) as consisting of a \emph{deterministic} observation \(e : S \to  O\) as well as a \emph{stochastic} update \(u : S \times  I \to  f S\) that defines a conditional probability kernel \(P(- \mid  s, i) := u(s, i)\) of the next state on the current state.\par{}We can augment the Giry monad to include a reward term \(\mathbb {R}\). For example, we can consider the endofunctor on measureable spaces that is defined on objects by \(X \mapsto   \mathbb {R} \times  f(X)\). Or by \(X \mapsto  f(\mathbb {R} \times  X)\) if we want to preserve monadicity.\par{}There are of course many different probability monads, such as the monad of finitary distributions on \(\mathsf {Set}\) \cite{fritz-2009-convex}, the Radon monad on the category of compact Hausdorff spaces \cite{keimel-2008-monad}, the probability monad on quasi-Borel spaces \cite{heunen-2017-convenient}, and many more. All of these give specific definitions of POMDPs with particular notions of state space and probability distribution.\end{example}

\begin{example}[{Other non-deterministic coalgebras}]\label{ssl-005U}
\par{}There are many different sorts of non-deterministic coalgebras other than \hyperref[ssl-005S]{the one captured by powerset}. For example, Fritz, Perrone, and Rezagholi \cite{fritz-2021-probability} construct a Hoare hyperspace monad on topological spaces sending a topological space to a space of its closed subsets. They also construct a probability monad on topological spaces and show that the support of a distribution gives a commutative monad morphism to the Hoare hyperspace monad. This gives an example of a morphism \(F \Rightarrow  G\) of endofunctors which, by the pseudo-functoriality of the doctrine, gives us a compositionality theorem relating \hyperref[ssl-005T]{POMDPs} (with respect to their probability monad) and \hyperref[ssl-005S]{non-deterministic Moore machines} (with repsect to the Hoare hyperspace monad).\end{example}
\subsection{Doctrine of ODEs in tangent categories}\label{djm-00HF}\par{}	In this section, we present a \hyperref[djm-00A8]{doctrine} of ordinary differential equations in \emph{tangent categories} \cite{cockett-2013-differential}, expressed as parameterized sections of tangent bundles. This expands the approach taken in Section 3.5.2 of \cite{jaz-2021-book}.\par{}	In \cite{capucci-2024-fibrational}, Capucci, Crutwell, Ghani, and Zanasi consider a notion of \textbf{first order differentiable structures} (FODS) (Definition 58 of \cite{capucci-2024-fibrational}) in a 2-category \(\mathcal {K}\) and compare this notion to various notions of tangent category appearing prior in the literature. A FODS in \(\mathcal {K}\) consists of a fibration \(\pi  : E \to  B\) in \(\mathcal {K}\) with structured biproducts together with a section \(T : B \to  E\).
\begin{definition}[{First order differential structures (FODS)}]\label{ssl-005V}
\par{}	Let \(\mathcal {K}\) be a 2-category with finite 2-limits and a finite 2-limit preserving forgetful functor \(U : \mathcal {K} \to  \mathcal {C}\mathsf {at}\). Define the 2-category \(\mathcal {F}\mathsf {ODS}(\mathcal {K})\) of \textbf{first order differential structures} in \(\mathcal {K}\) (following Definition 58 of \cite{capucci-2024-fibrational}) as the following pullback of 2-categories:
	
\begin{center}
\begin {tikzcd}
	{\mathcal {F}\mathsf {ODS}(\mathcal {K})} && {\mathsf {2}\mathcal {C}\mathsf {at}^{\mathsf {colax}}(\Delta [2], \mathcal {K})} \\
	{\mathcal {F}\mathsf {ib}^+(\mathcal {K})} & {2\mathcal {C}\mathsf {at}(\Delta [1], \mathcal {C}\mathsf {at}) \times  \mathcal {C}\mathsf {at}} & {\mathsf {2}\mathcal {C}\mathsf {at}^{\mathsf {colax}}(\Delta [1], \mathcal {K})^2}
	\arrow [from=1-1, to=1-3]
	\arrow [from=1-1, to=2-1]
	\arrow ["\lrcorner "{anchor=center, pos=0.125}, draw=none, from=1-1, to=2-3]
	\arrow ["{(d_0^{\ast }, d_1^{\ast })}", from=1-3, to=2-3]
	\arrow ["{(\mathsf {forget}, d_0^{\ast })}", from=2-1, to=2-2]
	\arrow ["{(\mathsf {inc}, \mathrm {id})}", from=2-2, to=2-3]
\end {tikzcd}		
	\end{center}

	where \(\mathcal {F}\mathsf {ib}^+(\mathcal {K})\) is the 2-category of fibrations with biproducts in \(\mathcal {K}\).\end{definition}

\begin{remark}[{Remark on tangency morphisms and tangent categories}]\label{djm-00HH}
\par{}	While morphisms of FODS aren't considered in \cite{capucci-2024-fibrational}, Lanfranchi does define a 2-category of tangent categories and ``lax'' morphisms in \cite{lanfranchi-2025-formal} (see Definitions 2.5, 2.7 and 2.10 of \cite{lanfranchi-2025-formal}). What Lanfranchi calls a ``lax'' morphism of tangent categories (Definition 2.7 of \cite{lanfranchi-2025-formal}) corresponds to \hyperref[djm-00BK]{our (colax) morphisms of tangencies} when considered as a FODS by Theorem 56 of \cite{capucci-2024-fibrational}. For this reason, we feel justified in using our morphisms of tangencies as those of tangent categories or FODS more generally.\par{}	In \cite{leung-2017-classifying}, Leung shows (Theorem 14.1 of \cite{leung-2017-classifying}) that a tangent category structure on a category \(C\) is equivalent to a pseudo-monoidal functor \(\mathsf {Weil}_1 \to  \mathsf {End}(C)\) from the monoidal category of first order Weil algebras satisfying two limit preservation properties. Such a monoidal functor corresponds to a pseudofunctor \(\mathcal {B}\mathsf {Weil}_1 \to  \mathcal {C}\mathsf {at}\). The morphisms Lanfranchi consideres in Definition 2.7 of \cite{lanfranchi-2025-formal} should correspond to \emph{colax} morphisms of these pseudo-functors.\par{}	A core idea of \cite{lanfranchi-2025-formal} is that Leung's tangent structures as pseudo-functors \(\mathcal {B}\mathsf {Weil}_1 \to  \mathcal {C}\mathsf {at}\) can be thought of as generalized monads (themselves functors \(\mathcal {B}\Delta _+ \to  \mathcal {C}\mathsf {at}\)). Colax morphisms of monads are the kind for which the Kleisli construction is functorial; colax morphisms of tangent categories (what Lanfranchi unfortunately calls "lax") are the kind for which vector fields are functorial. This is ultimately why our use of \hyperref[djm-00BK]{colax morphisms of tangencies} is the correct choice for understanding the compositionality of systems of ODEs.\end{remark}
\par{}	Next we will see that if \(\mathcal {K}\) is a concrete 2-category with a finite 2-limit preserving forgetful 2-functor \(U : \mathcal {K} \to  \mathcal {C}\mathsf {at}\), then a FODS will induce a \hyperref[djm-00BK]{tangency}.
\begin{definition}[{Doctrine of ODEs in first order differential structures}]\label{djm-00HG}
\par{}	Let \(\mathcal {K}\) be a concrete 2-category with a finite 2-limit preserving forgetful 2-functor \(U : \mathcal {K} \to  \mathcal {C}\mathsf {at}\). The forgetful functor \(U : \mathcal {K} \to  \mathcal {C}\mathsf {at}\)  determines a square of cartesian 2-functors:
	
\begin{center}
		\begin {tikzcd}
	{\mathcal {F}\mathsf {ODS}(\mathcal {K})} & {\mathcal {T}\mathsf {an}} \\
	{\mathcal {F}\mathsf {ib}^+(\mathcal {K})} & \mathcal {F}\mathsf {ib}
	\arrow ["U", from=1-1, to=1-2]
	\arrow [from=1-1, to=2-1]
	\arrow [from=1-2, to=2-2]
	\arrow ["U"', from=2-1, to=2-2]
\end {tikzcd}
	\end{center}

	which when post composing with the \hyperref[djm-00HA]{doctrine of systems in a tangency} gives us a \hyperref[djm-00A8]{doctrine} of \textbf{ordinary differential equations} in a first order differential structure:
	
\begin{center}
		\begin {tikzcd}
	{\mathcal {F}\mathsf {ODS}(\mathcal {K})} & {\mathcal {T}\mathsf {an}} & {\ell }\mathcal {M}\mathsf {od}_{\mathsf {r}} \\
	{\mathcal {F}\mathsf {ib}^+(\mathcal {K})} & \mathcal {F}\mathsf {ib} & \mathcal {D}\mathsf {bl}
	\arrow ["U", from=1-1, to=1-2]
	\arrow [from=1-1, to=2-1]
	\arrow ["{\mathbb {M}\mathsf {oore}}", from=1-2, to=1-3]
	\arrow [from=1-2, to=2-2]
	\arrow [from=1-3, to=2-3]
	\arrow ["U"', from=2-1, to=2-2]
	\arrow ["{\mathbb {L}\mathsf {ens}}"', from=2-2, to=2-3]
\end {tikzcd}
	\end{center}\end{definition}
\par{}	There are many examples of FODS, tangent categories, and cartesian differential categories. Each of these is a systems theory that will produce a module of systems whose systems are generalized ordinary differential equations expressed as a parameterized section of a tangent bundle:
	
\begin{center}
		\begin {tikzcd}
			TS & {e^{\ast }I} & I \\
			S & S & O
			\arrow [dashed, maps to, from=1-1, to=2-1]
			\arrow ["u"', from=1-2, to=1-1]
			\arrow ["{\mathsf {lift}(e)}", from=1-2, to=1-3]
			\arrow [dashed, maps to, from=1-2, to=2-2]
			\arrow ["\lrcorner "{anchor=center, pos=0.125}, draw=none, from=1-2, to=2-3]
			\arrow [dashed, maps to, from=1-3, to=2-3]
			\arrow [equals, from=2-1, to=2-2]
			\arrow ["e"', from=2-2, to=2-3]
		\end {tikzcd}
	\end{center}
\begin{example}[{Euclidean ODEs}]\label{djm-00IJ}
\par{}	Let \(\mathsf {Euc}\) be the category whose objects are natural numbers \(n\) and whose maps are smooth functions \(\mathbb {R}^n \to  \mathbb {R}^m\). The category \(\mathsf {Euc}\) is a \emph{cartesian differential category}, which are an example of a first-order differential structure \cite{capucci-2024-fibrational}. Explicitly, the underlying fibration of the FODS induced by \(\mathsf {Euc}\) is the \hyperref[djm-00E4]{simple fibration} associated to \(\mathsf {Euc}\). The section of the FODS induced by \(\mathsf {Euc}\) on objects sends \(\mathbb {R}^n\) to \(\mathbb {R}^n \times  \mathbb {R}^n \to  \mathbb {R}^n\), and on morphisms sends \(f : \mathbb {R}^n \to  \mathbb {R}^m\) to its derivative
	\[
		\begin {aligned}
			\mathbb {R}^n \times  \mathbb {R}^n & \to  \mathbb {R}^m \times  \mathbb {R}^m\\
			(x, v) &\mapsto  (f(x), D_xf(v))
		\end {aligned}
	\]\par{}	The FODS induced by \(\mathsf {Euc}\) is a systems theory which produces a module of systems whose:
	\begin{itemize}\item{}Interfaces consist of numbers \(n_o\) and \(n_i\) of exposed variable and parameters, respectively.
		\item{Systems consist of:}\begin{itemize}\item{}					A number \(n_s\) of \emph{state variables}.
				\item{}					A map \(\vec {e} : \mathbb {R}^{n_s} \to  \mathbb {R}^{n_o}\) exposing \(n_o\) variables of state.
				\item{}					A map \(\vec {u} : \mathbb {R}^{n_s} \times  \mathbb {R}^{n_i} \to  \mathbb {R}^{n_s}\) giving a vector of state displacement vectors, representing the system of ODEs
					\[
						\frac {d\vec {x}}{dt} = \vec {u}(\vec {x}, \vec {i}).
					\]\end{itemize}\end{itemize}\end{example}
\section{Restricting doctines to free interactions}\label{djm-00G5}\par{}	In \hyperref[ssl-003W]{Section \ref{ssl-003W}}, \hyperref[ssl-0041]{Section \ref{ssl-0041}}, and \hyperref[djm-00G6]{Section \ref{djm-00G6}}, we constructed a number of examples of doctrines. However, the resulting modules of systems produced by these doctrines were over double categories of interactions which were qualitatively complex, often including the systems and system maps themselves as special cases. therefore, composing systems theories through these interactions could have arbitarily complex effects. For example:
	\begin{enumerate}\item{}			The  were loose right modules on the double category of \emph{all} (co)spans in \(\mathsf {C}\). These (co)spans could include arbitary component systems; they certainly do more than \emph{just} compose.
		\item{}			The \hyperref[djm-00G6]{systems theories of generalized Moore machines associated to tangencies} are acted on by all \hyperref[djm-00B6]{generalized lenses}; these do more than compose, since they can allow arbitrary intermediate transformations of the outputs and parameters of the Moore machines.\end{enumerate}
	In this section, we will show how to restrict these more general interactions to simpler \emph{coupling schemes} or \emph{wiring diagrams}. We base our approach on the following observation.
\begin{observation}[{Wiring diagrams are free interactions}]\label{djm-00G7}
\par{}	A wiring diagram for a given \hyperref[djm-00A8]{doctrine} is an interaction in a free theory of interactions in that doctrine.\end{observation}
\subsection{Wiring diagrams are free interactions}\label{ssl-0042}\par{}	Recall that given a \hyperref[djm-00A8]{doctrine} is a commuting diagram as follows.
	
  \begin{center}
		\begin {tikzcd}
			{\mathcal {D}\mathsf {}_{\mathsf {sys}}} & {{\ell }\mathcal {M}\mathsf {od}_{\mathsf {r}}} \\
			{\mathcal {D}\mathsf {}_{\mathsf {inter}}} & {\mathcal {D}\mathsf {bl}}
			\arrow ["{\mathbb {S}}", from=1-1, to=1-2]
			\arrow ["{\pi _{\mathcal {D}\mathsf {}}}"', from=1-1, to=2-1]
			\arrow ["{(-)_1}",from=1-2, to=2-2]
			\arrow ["{\mathbb {I}}"', from=2-1, to=2-2]
		\end {tikzcd}
\end{center}
	Given a systems theory \(T \in  \mathcal {D}\mathsf {}_{\mathsf {sys}}\), we call \(\pi _{\mathcal {D}\mathsf {}}(T)\) the \emph{interaction theory}. Under \(\mathbb {I}\), \(T\) produces a double category of interactions which acts on the systems in the modulee of systems defined by \(T\). The loose arrows of this double category of interactions are what we term \emph{interactions}. In the literature, many interactions are in the form of \emph{wiring diagrams} and that these are freely generated interactions (see  \cite{spivak-2013-operad}, \cite{vagner-2014-algebras}, \cite{libkind-2022-operadic}).\par{}	Let's see \hyperref[djm-00G7]{Observation \ref{djm-00G7}} borne out in a few special cases.
\begin{example}[{Undirected wiring diagrams are free interactions}]\label{djm-00G8}
\par{}			Systems in the  are well known to form \emph{hypergraph categories}, which by the work of Fong and Spivak \cite{fong-2018-hypergraph} are equivalently described as algebras for operads of cospans of typed finite sets (see \cite{fong-2018-hypergraph} for a review of hypergraph categories and how they arise as (co)spans and (co)relations).\par{}			If \(X\) is a set of "types", then the category of typed finite sets is defined to be the slice category \(\mathsf {Finset} \downarrow  X\) consisting of finite sets \(I\) equipped with a typing map \(\tau  : I \to  X\). Cospans of typed finite sets are interpreted as \emph{undirected wiring diagrams} (\cite{spivak-2013-operad}).

      The following undirected wiring diagram --- introduced in \hyperref[ssl-004E]{Example \ref{ssl-004E}} --- is an undirected wiring diagram over a single set of types.
      \begin{center}\includegraphics[scale=1]{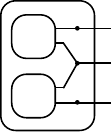}\end{center}

			Composition of such cospans gives nesting.\par{}Next, we will show how these undirected wiring diagrams arise as free interactions in the \hyperref[djm-00GQ]{port-plugging doctrine}. As the variable sharing nad port-plugging doctrines are highly related, we treat them together.\par{}	The \hyperref[djm-00GP]{variable sharing doctrine} and the \hyperref[djm-00GQ]{port-plugging doctrine}  have \(\mathcal {L}\mathsf {ex}\) and \(\mathcal {R}\mathsf {ex}\) as their \hyperref[djm-00A8]{interaction theory}, respectively. Both \(\mathcal {L}\mathsf {ex}\) and \(\mathcal {R}\mathsf {ex}\) are monadic over \(\mathcal {C}\mathsf {at}\) via adjunctions \(U : \mathcal {L}\mathsf {ex} \leftrightarrows  \mathcal {C}\mathsf {at} : F\) and \(U : \mathcal {R}\mathsf {ex} \leftrightarrows  \mathcal {C}\mathsf {at} : F\). If \(X\) is a set of "types", considered as a discrete category, then it generates a \emph{free interaction theory} \(FX\) in \(\mathcal {L}\mathsf {ex}\) or \(\mathcal {R}\mathsf {ex}\) (considering \(X\) as a discrete category).\par{}	Let us first consider the case of \(\mathcal {R}\mathsf {ex}\). The free finite colimit completion on a category is given by the finite colimits of representable presheaves. In the case that \(X\) is a discrete category, the category \(\mathsf {Psh}(X)\) of presheaves on \(X\) is equivalent to the slice category \(\mathsf {Set} \downarrow  X\). Therefore, the free rex category \(FX\) on \(X\) is the slice category \(\mathsf {Finset} \downarrow  X\) of \emph{\(X\)-typed finite sets}. So under the port-plugging doctrine, the interaction theory \(FX\) defines the double category of interactions the symmetric monoidal double category of cospans of \(X\)-typed finite sets, \(\mathbb {C}\mathsf {ospan}(\mathsf {Finset} \downarrow  X)\). We therefore see that a \emph{free interaction} in the port-plugging doctrine is precisely an undirected wiring diagram in then sense of \cite{spivak-2013-operad}.\par{}	Now consider the case of \(\mathcal {L}\mathsf {ex}\). The free finite limit category on a category \(C\) is the opposite of the free finite colimit completion of \(C^{\mathsf {op}}\). By the reasoning above, we therefore see that the free lex category on a discrete category \(X\) is \((\mathsf {Finset} \downarrow  X)^{\mathsf {op}}\). Therefore under the variable sharing doctrine, the interaction theory \(FX\) defines the double category of interactions \(\mathbb {S}\mathsf {pan}((\mathsf {Finset} \downarrow  X)^{\mathsf {op}})\); but this is equivalently \({\mathbb {C}\mathsf {ospan}(\mathsf {Finset} \downarrow  X)}^{\mathsf {op}}\). Again, we find that free interactions in the span doctrine are precisely undirected wiring diagrams.\end{example}

\begin{example}[{Directed wiring diagrams are free interactions}]\label{djm-00G9}
\par{}	In \cite{vagner-2014-algebras}, Vagner, Spivak, and Lerman consider algebras of dynamical systems (specifically, ODEs) on operads of \emph{directed wiring diagrams}. These diagrams are defined to be \emph{prisms} in the cocartesian category of typed finite sets, which are by definition lenses in the cartesian opposite of the category of typed finite sets.\par{}The following is a directed wiring diagram with a single type.
  \begin{center}\includegraphics[scale=1]{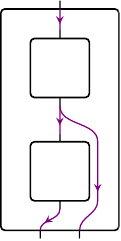}\end{center}\par{}	Now, the doctrine of simple Moore machines in a cartesian category has cartesian categories as its interaction theories. The 2-category \(\mathcal {C}\mathsf {art}\) of cartesian categories is monadic over \(\mathcal {C}\mathsf {at}\) via the adjunction \(U : \mathcal {C}\mathsf {art} \leftrightarrows  \mathcal {C}\mathsf {at} : F\), where for a category \(\mathsf {C}\) the free cartesian category \(F\mathsf {C}\) may be constructed as the opposite of a colax slice:
	\[
		FC := (\mathsf {Finset} \downarrow ^{\mathsf {colax}} \mathsf {C})^{\mathsf {op}}
	\]
	In particular, we find that free cartesian categories on a discrete category (set) \(T\) are the opposite categories of typed finite sets. Therefore, interactions in the theory of simple Moore machines --- which are lenses --- in these free cartesian categories are precisely directed wiring diagrams. This observation was made in Section 1.3.3 of \cite{jaz-2021-book}.\end{example}
\par{}	Furthermore, \hyperref[djm-00G7]{Observation \ref{djm-00G7}} suggests ways to define wiring diagrams which may have "function boxes" that transform values on the wires, or other such decorations: simply work in the free interaction theories generated by algebraic theories of the "functions" which you would like to place on the boxes (see Section 1.3.4 of \cite{jaz-2021-book}).
	For example, the free interaction theory consisting of  lenses in the Lawvere theory generated by a single unary operation \(d : X \to  X\) gives directed wiring diagrams \emph{with delays}.\subsection{Restricting a doctrine to free interactions}\label{ssl-005A}\par{}	We will now describe a general construction that lets us restrict a doctrine so that its interaction theories are freely generated. 

	In what follows, assume that we are given a \hyperref[djm-00A8]{doctrine}:
  \begin{center}
		\begin {tikzcd}
			{\mathcal {D}\mathsf {}_{\mathsf {sys}}} & {{\ell }\mathcal {M}\mathsf {od}_{\mathsf {r}}} \\
			{\mathcal {D}\mathsf {}_{\mathsf {inter}}} & {\mathcal {D}\mathsf {bl}}
			\arrow ["{\mathbb {S}}", from=1-1, to=1-2]
			\arrow ["{\pi _{\mathcal {D}\mathsf {}}}"', from=1-1, to=2-1]
			\arrow ["{(-)_1}",from=1-2, to=2-2]
			\arrow ["{\mathbb {I}}"', from=2-1, to=2-2]
		\end {tikzcd}
\end{center} and a pair of functors \(F : \mathcal {C} \leftrightarrows  \mathcal {D}\mathsf {}_{\mathsf {inter}} : U\) which form a 2-adjunction \(F \dashv  U\).
\begin{definition}[{The 2-category of restricted systems theories to free interactions}]\label{ssl-0057}
\par{}Define the cartesian 2-category \(\mathcal {C} \downarrow ^{\mathsf {adm}} U\pi _{\mathcal {D}\mathsf {}}\) to be the sub-2-category of the colax slice \(\mathcal {C} \downarrow ^{\mathsf {colax}} U\pi _{\mathcal {D}\mathsf {}}\) spanned by all objects and 1-cells 
  
\begin{center}
    \begin {tikzcd}
      {C_0} & {C_1} \\
      {U\pi _{\mathcal {D}\mathsf {}}(T_0)} & {U\pi _{\mathcal {D}\mathsf {}}(T_1)}
      \arrow ["{f_0}", from=1-1, to=1-2]
      \arrow ["{m_0}"', from=1-1, to=2-1]
      \arrow ["{m_1}", from=1-2, to=2-2]
      \arrow ["{\overline {f}}", shorten <=14pt, shorten >=14pt, Rightarrow, from=2-1, to=1-2]
      \arrow ["{U\pi _{\mathcal {D}\mathsf {}}(f_1)}"', from=2-1, to=2-2]
    \end {tikzcd}	
  \end{center}

  whose colaxator \(\overline {f}\) is \emph{admissible} meaning that \(\mathbb {I}\) sends its transpose \(\mathsf {tr}(\overline {f})\) in \(\mathcal {D}\mathsf {}_{\mathsf {inter}}\) along the \(F \dashv  U\) adjunction to a \hyperref[djm-00EX]{companion commuter transformation}.\par{}For notational clarity we define \[\mathcal {D}\mathsf {}_{\mathsf {sys}}|_F  \coloneqq  \mathcal {C} \downarrow ^{\mathsf {adm}} U\pi _{\mathcal {D}\mathsf {}}.\]\end{definition}

\begin{explication}[{The 2-category of restricted systems theories to free interactions}]\label{ssl-0058}
\par{}The cartesian 2-category \(\mathcal {D}\mathsf {}_{\mathsf {sys}}|_F \) is the 2-category of systems theories when we restrict the doctrine \((\pi _{\mathcal {D}\mathsf {}} \colon  \mathcal {D}\mathsf {}_{\mathsf {sys}} \to   \mathcal {D}\mathsf {}_{\mathsf {inter}}, \mathbb {S}, \mathbb {I})\) along the left adjoint \(F \colon  \mathcal {C} \to  \mathcal {D}\mathsf {}_{\mathsf {inter}}\).\par{}A theory in this restricted systems theory consists of an object \(c \in  \mathcal {C}\), a theory \(T \in  \mathcal {D}\mathsf {}_{\mathsf {sys}}\), and a \textbf{marking} \(m \colon  c \to  U \pi _{\mathcal {D}\mathsf {}}(T)\), which defines the generating interactions for the restricted systems theories.\par{}Note that the transpose of a marking is the morphism \(\mathsf {tr}(m) \colon  Fc \to  \pi _{\mathcal {D}\mathsf {}}(T)\) in \(\mathcal {D}\mathsf {}_{\mathsf {inter}}\).\end{explication}
\par{}	We give several examples of restricted doctrines in \hyperref[ssl-0043]{Section \ref{ssl-0043}}. In the meantime, we focus on a restriction of \hyperref[djm-00GQ]{the port-plugging doctrine} in order to make sense of these definitions.
\begin{example}[{Systems theories in the port-plugging doctrine restricted along \(\mathcal {C}\mathsf {at} \to  \mathcal {R}\mathsf {ex}\)}]\label{ssl-0059}
\par{}In the \hyperref[djm-00GQ]{port-plugging doctrine} the module of systems induced by \hyperref[ssl-003X]{the rex category \(\mathsf {Petri}\) of Petri nets} is over the double category of interactions \(\mathbb {C}\mathsf {ospan}(\mathsf {Petri})\). In other words, Petri nets interact by gluing along any sub-Petri net. This formulation of interaction is too unrestricted. We may want --- following the original work \cite{baez-2017-reaction}  --- to have Petri nets interact by gluing along \emph{species}.\par{}We can restrict the \hyperref[djm-00GQ]{port-plugging doctrine} along the free rex category construction \(F \colon  \mathcal {C}\mathsf {at} \to  \mathcal {R}\mathsf {ex}\).\par{}An example of a theory of this restricted systems theory consists of the terminal category \(\bullet \), the rex category \(\mathsf {Petri}\), and a marking \(m \colon  \bullet  \to  \mathsf {Petri}\) which sends \(\bullet \) to the Petri net with a single species.  We will continue this example in \hyperref[ssl-005D]{Example \ref{ssl-005D}}.\end{example}
\par{}Next we will define a functor from a restricted systems theory to \hyperref[ssl-0040]{right niches} and then restrict to get a module of systems.
\begin{lemma}[{Restricted systems theories to right niches}]\label{ssl-005B}
\par{}There is a cartesian 2-functor \(\mathcal {D}\mathsf {}_{\mathsf {sys}}|_F \to  \mathcal {N}\mathsf {iche}_r\) that maps a \hyperref[ssl-0058]{restricted theory} \((c \in  \mathcal {C}, T \in  \mathcal {D}\mathsf {}_{\mathsf {sys}}, m \colon  c \to  U \pi _{\mathcal {D}\mathsf {}}(T))\) to the right niche 
  
\begin{center}
    \begin {tikzcd}
      \bullet  & {\mathbb {I}(Fc)} & {} \\
      \bullet  & {\mathbb {S}(T)_1  = \mathbb {I} \pi _{\mathcal {D}\mathsf {}}}
      \arrow [equals, from=1-1, to=2-1]
      \arrow ["{\mathbb {I} (\mathsf {tr}(m))}", from=1-2, to=2-2]
      \arrow ["{\mathbb {S}(T)}"', "\shortmid "{marking}, from=2-1, to=2-2]
    \end {tikzcd}
  \end{center}\begin{proof}\par{}We begin by observing that the \(F \vdash  U\) adjunction defines an isomorphism \[C \downarrow ^{\mathsf {adm}} U\pi _{\mathcal {D}\mathsf {}} \cong  F \downarrow ^{\mathsf {adm}} \pi _{\mathcal {D}\mathsf {}},\] where \(F \downarrow ^{\mathsf {adm}} \pi _{\mathcal {D}\mathsf {}}\) is the sub-2-category of the colax slice \(F \downarrow ^{\mathsf {colax}} \pi _{\mathcal {D}\mathsf {}}\) spanned by all objects and 1-cells which are sent to a \hyperref[djm-00EX]{companion commuter transformation} under \(\mathbb {I}\).
    
    Explicitly, the 1-cells of \(F \downarrow ^{\mathsf {adm}} \pi _{\mathcal {D}\mathsf {}}\) are objects of the following pullback

\begin{center}
      \begin {tikzcd}
        {\mathsf {2}\mathcal {C}\mathsf {at}^{\mathsf {colax}}(\Delta [1], \mathcal {D}\mathsf {}_{\mathsf {inter}})_{\mathsf {comp}}} & {\mathsf {2}\mathcal {C}\mathsf {at}^{\mathsf {colax}}(\Delta [1], \mathcal {D}\mathsf {bl})_{\mathsf {comp}}} \\
        {\mathsf {2}\mathcal {C}\mathsf {at}^{\mathsf {colax}}(\Delta [1], \mathcal {D}\mathsf {}_{\mathsf {inter}})} & {\mathsf {2}\mathcal {C}\mathsf {at}^{\mathsf {colax}}(\Delta [1], \mathcal {D}\mathsf {bl})} \\
        {\mathcal {D}\mathsf {}_{\mathsf {inter}} \times  \mathcal {D}\mathsf {}_{\mathsf {inter}}} & {\mathcal {D}\mathsf {bl} \times  \mathcal {D}\mathsf {bl}}
        \arrow [from=1-1, to=1-2]
        \arrow [hook, from=1-1, to=2-1]
        \arrow [hook, from=1-2, to=2-2]
        \arrow [""{name=0, anchor=center, inner sep=0}, "{\mathbb {I} \circ   -}"', from=2-1, to=2-2]
        \arrow [""{name=0p, anchor=center, inner sep=0}, phantom, from=2-1, to=2-2, start anchor=center, end anchor=center]
        \arrow ["{(d_1^*, d_0^*)}"', from=2-1, to=3-1]
        \arrow ["{(d_1^*, d_0^*)}", from=2-2, to=3-2]
        \arrow ["{\mathbb {I} \times  \mathbb {I}}"', from=3-1, to=3-2]
        \arrow ["\lrcorner "{anchor=center, pos=0.125, rotate=0}, draw=none, from=1-1, to=0p]
      \end {tikzcd}
    \end{center}

    while the 2-category \(F \downarrow ^{\mathsf {adm}} \pi _{\mathcal {D}\mathsf {}}\) is itself the pullback 
    
\begin{center}
      \begin {tikzcd}
        {F \downarrow ^{\mathsf {adm}} \pi _{\mathcal {D}\mathsf {}}} & {\mathsf {2}\mathcal {C}\mathsf {at}^{\mathsf {colax}}(\Delta [1], \mathcal {D}\mathsf {}_{\mathsf {inter}})_{\mathsf {comp}}} && {} \\
        {\mathcal {C} \times  \mathcal {D}\mathsf {}_{\mathsf {sys}}} & {\mathcal {D}\mathsf {}_{\mathsf {inter}} \times  \mathcal {D}\mathsf {}_{\mathsf {inter}}} & {} & {}
        \arrow [from=1-1, to=1-2]
        \arrow [from=1-1, to=2-1]
        \arrow ["{(d_1^*, d_0^*)}", from=1-2, to=2-2]
        \arrow [""{name=0, anchor=center, inner sep=0}, "{F \times  \pi _{\mathcal {D}\mathsf {}}}"', from=2-1, to=2-2]
        \arrow ["\lrcorner "{anchor=center, pos=0.125}, draw=none, from=1-1, to=0]
      \end {tikzcd}
    \end{center}\par{}The desired map from  restricted systems theory to right niches is induced by the following commutative diagram.

  \begin{center}
      \begin {tikzcd}[column sep = tiny]
        {F \downarrow ^{\mathsf {adm}} \pi _{\mathcal {D}\mathsf {}}} & {\mathcal {D}\mathsf {}_{\mathsf {sys}}} & {\ell }\mathcal {M}\mathsf {od}_{\mathsf {r}} \\
        & {\mathcal {N}\mathsf {iche}_r} & \mathcal {N}\mathsf {iche} & {\ell }\mathcal {B}\mathsf {imod} \\
        & {\mathcal {D}\mathsf {bl}_{\bullet } \times  \mathsf {2}\mathcal {C}\mathsf {at}^{\mathsf {colax}}(\Delta [1], \mathcal {D}\mathsf {bl})_{\mathsf {comp}}} & {\mathsf {2}\mathcal {C}\mathsf {at}^{\mathsf {colax}}(\Delta [1], \mathcal {D}\mathsf {bl})_{\mathsf {conj}} \times  \mathsf {2}\mathcal {C}\mathsf {at}^{\mathsf {colax}}(\Delta [1], \mathcal {D}\mathsf {bl})_{\mathsf {comp}}} & {\mathcal {D}\mathsf {bl} \times  \mathcal {D}\mathsf {bl}}
        \arrow [from=1-1, to=1-2]
        \arrow [dashed, from=1-1, to=2-2]
        \arrow [from=1-1, to=3-2]
        \arrow ["{\mathbb {S}}", from=1-2, to=1-3]
        \arrow [hook, from=1-3, to=2-4]
        \arrow [from=2-2, to=2-3]
        \arrow [from=2-2, to=3-2]
        \arrow [from=2-3, to=2-4]
        \arrow [from=2-3, to=3-3]
        \arrow ["{(()_0, ()_1)}", from=2-4, to=3-4]
        \arrow [""{name=0, anchor=center, inner sep=0}, from=3-2, to=3-3]
        \arrow [""{name=1, anchor=center, inner sep=0}, "{d_0^* \times  d_0^*}"', from=3-3, to=3-4]
        \arrow [""{name=1p, anchor=center, inner sep=0}, phantom, from=3-3, to=3-4, start anchor=center, end anchor=center]
        \arrow ["\lrcorner "{anchor=center, pos=0.125}, draw=none, from=2-2, to=0]
        \arrow ["\lrcorner "{anchor=center, pos=0.125}, draw=none, from=2-3, to=1p]
      \end {tikzcd}
    \end{center}

    where the map \(F \downarrow ^{\mathsf {adm}} \pi _{\mathcal {D}\mathsf {}} \to  \mathcal {D}\mathsf {bl}_{\bullet } \times  \mathsf {2}\mathcal {C}\mathsf {at}^{\mathsf {colax}}(\Delta [1], \mathcal {D}\mathsf {bl})_{\mathsf {comp}}\) is constant at the terminal pointed double category on the left and the composite of the forgetful 2-functors \[F \downarrow ^{\mathsf {adm}} \pi _{\mathcal {D}\mathsf {}} \to  \mathsf {2}\mathcal {C}\mathsf {at}^{\mathsf {colax}}(\Delta [1], \mathcal {D}\mathsf {}_{\mathsf {inter}})_{\mathsf {comp}}\] and \[\mathsf {2}\mathcal {C}\mathsf {at}^{\mathsf {colax}}(\Delta [1], \mathcal {D}\mathsf {}_{\mathsf {inter}})_{\mathsf {comp}} \to  \mathsf {2}\mathcal {C}\mathsf {at}^{\mathsf {colax}}(\Delta [1], \mathcal {D}\mathsf {bl})_{\mathsf {comp}}\] defined by the pullbacks above.\end{proof}\end{lemma}

\begin{definition}[{Restriction of a doctrine to free interactions}]\label{djm-00GR}
\par{}The \textbf{restricted doctrine} of the doctrine \((\pi _{\mathcal {D}\mathsf {}} \colon  \mathcal {D}\mathsf {}_{\mathsf {sys}} \to   \mathcal {D}\mathsf {}_{\mathsf {inter}}, \mathbb {S}, \mathbb {I})\) along the left adjoint \(F \colon  \mathcal {C} \to  \mathcal {D}\mathsf {}_{\mathsf {inter}}\) is defined by the square
	
\begin{center}
		\begin {tikzcd}
			{ \mathcal {D}\mathsf {}_{\mathsf {sys}}|_F} & {{\ell }\mathcal {M}\mathsf {od}_{\mathsf {r}}} \\
			{\mathcal {C}} & {\mathcal {D}\mathsf {bl}}
			\arrow ["{\mathbb {S}|_F}", from=1-1, to=1-2]
			\arrow ["{\pi _{\mathcal {D}\mathsf {}}|_F}"', from=1-1, to=2-1]
			\arrow [from=1-2, to=2-2]
			\arrow ["{\mathbb {I} \circ  F}"', from=2-1, to=2-2]
		\end {tikzcd}
	\end{center}

	where \({\mathbb {S}|_F}\) is the composite of the 2-functor defined in \hyperref[ssl-005B]{Lemma \ref{ssl-005B}} and the \hyperref[ssl-004H]{restriction} \(\mathsf {Res}_r \colon  \mathcal {N}\mathsf {iche}_r \to  {\ell }\mathcal {M}\mathsf {od}_{\mathsf {r}}\).\end{definition}

\begin{explication}[{Restriction of a doctrine to free interactions}]\label{ssl-005C}
\par{}Recall that a \hyperref[ssl-0058]{restricted systems theory} is a triple \((c \in  \mathcal {C}, T \in  \mathcal {D}\mathsf {}_{\mathsf {sys}}, m \colon  c \to  U \pi _{\mathcal {D}\mathsf {}}(T))\). Under the \hyperref[djm-00GR]{restricted doctrine} this systems theory produces a module of systems that is the restriction of the right niche presented in \hyperref[ssl-005B]{Lemma \ref{ssl-005B}}. This module essentially restricts the double category of interactions to those freely produced by the data \(c \in  \mathcal {C}\) and contains all systems and system maps whose interfaces and interface maps live in this restricted double category of interactions. The action of interactions on systems is the same as in the original doctrine.\end{explication}

\begin{example}[{The port-plugging doctrine restricted along \(\mathcal {C}\mathsf {at} \to  \mathcal {R}\mathsf {ex}\)}]\label{ssl-005D}
\par{}In \hyperref[ssl-0059]{Example \ref{ssl-0059}}, we saw that an example of a systems theory in the port plugging doctrine restricted along the free rex construction \(\mathcal {C}\mathsf {at} \to  \mathcal {R}\mathsf {ex}\) consists of the triple \[(\bullet  \in  \mathcal {C}\mathsf {at}, \mathsf {Petri} \in  \mathcal {R}\mathsf {ex}, m \colon  \bullet  \to  \mathsf {Petri} )\] which sends \(\bullet \) to the Petri net with a single species.\par{}The free rex category of the terminal category is \(\mathsf {Finset}\). Therefore, in the restricted doctrine, this restricted theory produces a module of systems over the double category of interactions \(\mathbb {C}\mathsf {ospan}(\mathsf {Finset})\) which is the restriction of the right niche
  
\begin{center}
    \begin {tikzcd}[column sep=huge]
      \bullet  & {\mathbb {C}\mathsf {ospan}(\mathsf {Finset})} \\
      \bullet  & {\mathbb {C}\mathsf {ospan}(\mathsf {Petri})}
      \arrow [equals, from=1-1, to=2-1]
      \arrow [from=1-2, to=2-2]
      \arrow ["{\mathbb {C}\mathsf {ospan}(\mathsf {Petri})(\emptyset , -)}"', "\shortmid "{marking}, from=2-1, to=2-2]
    \end {tikzcd}
  \end{center}

  where \(\mathbb {C}\mathsf {ospan}(\mathsf {Finset}) \to  \mathbb {C}\mathsf {ospan}(\mathsf {Petri})\) is induced by interpretting finite sets as discrete Petri nets with no transitions.\par{}In this module of systems:
  \begin{itemize}\item{}An interface is a finite set.
    \item{}A system with interface \(M\) consists of a Petri net \(P\) and  a map from \(M\) to the species of \(P\). Note that this map is equivalent to a map of Petri nets from the Petri net with \(M\) species and no transitions to \(P\).
    \item{}An interaction is a cospan of finite sets.
    \item{}Interactions act on Petri nets by gluing along species identified in the interaction.\end{itemize}
  This is exactly the data of the hypergraph double category of Petri nets defined in \cite{baez-2020-open}. We illustrated this module of systems in \hyperref[ssl-001K]{Section \ref{ssl-001K}}.\end{example}
\subsection{Examples of restricting doctrines}\label{ssl-0043}\par{}In this section, we construct the systems theories defined in \hyperref[ssl-001F]{Section \ref{ssl-001F}} and others, using the machinery built up in \hyperref[ssl-003W]{Section \ref{ssl-003W}}, \hyperref[ssl-0041]{Section \ref{ssl-0041}},  \hyperref[djm-00G6]{Section \ref{djm-00G6}}, and \hyperref[ssl-0042]{Section \ref{ssl-0042}}.
\begin{example}[{Examples of the port-plugging doctrine restricted along \(\mathcal {C}\mathsf {at} \to  \mathcal {R}\mathsf {ex}\)}]\label{ssl-005X}
\par{}Let \(F \colon  \mathcal {C}\mathsf {at} \to  \mathcal {R}\mathsf {ex}\) be the left-adjoint that maps a category to its free rex category. Following \cite{fong-2018-hypergraph}, which gives an equivalence between \emph{hypergraph categories} (with symmetric monoidal structure free on objects) and algebras over operads of cospans, we may think of a module of systems over \(\mathbb {C}\mathsf {ospan}(F\mathsf {C})\) as a \emph{hypergraph double category}. Systems theories in the port-plugging doctrine restricted along the left adjoint \(\mathcal {C}\mathsf {at} \to  \mathcal {R}\mathsf {ex}\) produce  modules of systems of this form.\par{}In \hyperref[ssl-005D]{Example \ref{ssl-005D}} we gave an example of a systems theory in the restricted port-plugging doctrine that marked Petri nets by their places. We then saw how the module of systems produced by this marking consists open Petri nets that glue along their places. 
  Here we will give several other examples of systems theories in this restricted doctrine and describe the module of systems that they define. 
  \begin{itemize}\item{}Directed graphs with edge-labels in a set \(L\) (including, for example, circuit diagrams as in Fong's thesis \cite{fong-2016-algebra}) may be expressed as the slice category \(\mathsf {Graph} \downarrow  BL\) where \(BL\) is the graph with a single node and edge set \(L\); as a slice category of a presheaf category, this is a rex category. There is a systems theory in the restricted port-plugging doctrine that consists of the terminal category \(\bullet \), the rex category of directed labeled graphs \(\mathsf {Graph} \downarrow  BL\), and the marking \(\bullet  \to  \mathsf {Graph} \downarrow  BL\) that selects the graph with a single vertex with its unique trivial edge-labelling. This systems theory produces a  a module of labelled graphs over the double category \(\mathbb {C}\mathsf {ospan}(\mathsf {Finset})\) of cospans of finite sets.  This module of systems is equivalent to the hypergraph \emph{double} category of labelled graphs whose objects and loose morphisms form  the hypergraph category of labelled graphs, per \cite{fong-2018-hypergraph} and whose  tight morphisms account for label-preserving graph homomorphisms.
    \item{}Similar to marking Petri nets by  places and labeled graphs by  vertices, we can form hypeter graph double category of causal loop diagrams by marking  objects.
    \item{}In Section 3.5 of \cite{baez-2023-compositional}, the authors describe a schema for stock-flow diagrams and a restricted schema for their interfaces (see also \cite{baez-2022-categorical}). A full-fledged stock-flow diagram is a presheaf on the schema for stock-flow diagrams along with an auxiliary function. 
       There is a systems theory which consists of the free category on the graph 
       \[\bullet  \leftarrow  \bullet  \rightarrow  \bullet ,\]
      the rex category of full-fledged stock-flow diagrams, and a marking which selects the representable stock-flow diagrams defined by the interface schema.
      This systems theory defines a module of full-fledged stock-flow diagrams  over cospans of presheaves on this restricted interface that recovers the compositionality of full-fledged stock-flow diagrams in \cite{baez-2023-compositional}.\end{itemize}\end{example}

\begin{remark}[{The operadic approach to systems defined as structured cospans}]\label{djm-00FR}
\par{}	Many of the module of systems defined in  \hyperref[ssl-005X]{Example \ref{ssl-005X}} are a direct repackaging of symmetric monoidal double categories of structured cospans. In his thesis \cite{courser-2020-open}, Courser proves that given a category \(\mathsf {C}\) with pushouts and a functor \(L \colon  \mathsf {I} \to  \mathsf {C}\), there is a symmetric monoidal double category of structured cospans \(_L \mathbb {C}\mathsf {ospan}(\mathsf {C})\) whose objects are the objects in \(\mathsf {I}\) and where a loose morphism is a structured cospan \(Li \to  s \leftarrow  Lj\).\par{}	In many examples, these structured cospans represent open systems which interact only via the monoidal product (the trivial parallel product interaction) or via composition (series interaction) in \(_L \mathbb {C}\mathsf {ospan}(\mathsf {C})\). For example,  systems \(Li \to  s \leftarrow  Lj\) and \(Li' \to  s' \leftarrow  L{j'}\) can only interact non-trivially if \(j = i'\)  in which they may be composed in series to produce the system \(L_i  \to  s +_{Lj} s' \leftarrow  L_j'\). Or conversely if \(i = j'\) in which case they may be composed in series in the opposite order.\par{}	Since non-trivial interactions are limited to composition and composition depends on the interface of systems, this framework leads one to the practice of designing a system's interface with a particular interaction in mind. As an example, consider the diagram on page 2 of \cite{baez-2020-open} where the place \(E\) of a Petri net is \emph{double labeled} (that is, we have interface \(i = \{4, 5\}\) and map \(Li \to  s\) sending both \(4\) and \(5\) to the stock \(E \in  s\))  so that it can compose with two different places \(C\) and \(D\) of another Petri net.\par{}	By contrast, the operadic approach taken in modules of systems allows one to design the interface up-front for many potential interactions. This feature is critical for best practices in modular design.\par{}	We can repackage the data of a symmetric monoidal double category of structured cospans in to a module of systems, thereby applying an operadic approach to systems defined as structured cospans.  A rex category  \(\mathsf {C}\) and a functor \(L \colon   \mathsf {I} \to  \mathsf {C}\) is a systems theory in the \hyperref[ssl-005D]{restricted port-plugging doctrine}. It produces a module of systems in which 
	\begin{itemize}\item{}A system is a cospan \(\emptyset  \to  s \leftarrow  Li\) in \(\mathsf {C}\) which is equivalent to a morphism \(Li \to  s\) in \(\mathsf {C}\).
		\item{}An interaction is a cospan  in \(\mathsf {I}\).\end{itemize}
	Component systems \(Li \to  s\) and \(Li' \to  s'\) interacting via the cospan \(i + i' \to  m \leftarrow  j\) defines the composite system \(Lj \to  (Lm) +_{Li + Li'} (s + s')\).

\begin{center}
		\begin {tikzcd}
			& {(Lm) +_{L_i + L_i'} (s + s')} \\
			{s+ s'} && Lm \\
			& {Li + Li'} && Lj
			\arrow ["\lrcorner "{anchor=center, pos=0.125, rotate=-45}, draw=none, from=1-2, to=3-2]
			\arrow [from=2-1, to=1-2]
			\arrow [from=2-3, to=1-2]
			\arrow [from=3-2, to=2-1]
			\arrow [from=3-2, to=2-3]
			\arrow [from=3-4, to=2-3]
		\end {tikzcd}
	\end{center}\par{}	In this \emph{operadic} point of view, which we take from Fong and Spivak's \emph{Hypergraph Categories} \cite{fong-2018-hypergraph}, the systems have their boundaries set up front and they remain valid for all compositions down the line. Changes in system boundary can happen \emph{though} specific interactions (for example, an interaction whose left leg is the identity will re-label the boundary without any gluing), rather than being a meta-theoretic operation performed to prepare a system for a particular interaction.\par{}	We can also adjust the double category of interactions for such modules of systems.\par{}	For example, we can  assume that the right leg of each interaction is a \emph{monomorphism} --- disallowing double labelling entirely --- with no loss of expressability for composition, since all the gluing happens along the left leg of the interaction.\par{}	As a second example, we could expand the double category of interactions to include all structured cospans \(Li \to  t \leftarrow  Lk\).  We can interpret these structured cospans as \emph{processes} by which systems may interact. Such a process may arise by partial interaction. For example, consider how the interaction \(i + i' \to  m \leftarrow  j\) acts on the system \(Li \to  s\) and the trivial system \(Li' \to  Li'\). It produces a composite system \(Lj \to  t\) where \(t = (Lm +_{Li} s) + Li'\).

\begin{center}
		\begin {tikzcd}
			& {(Lm +_{Li} s) + Li'} \\
			{s+ Li'} && Lm \\
			& {Li + Li'} && Lj
			\arrow ["\lrcorner "{anchor=center, pos=0.125, rotate=-45}, draw=none, from=1-2, to=3-2]
			\arrow [from=2-1, to=1-2]
			\arrow [from=2-3, to=1-2]
			\arrow [from=3-2, to=2-1]
			\arrow [from=3-2, to=2-3]
			\arrow [from=3-4, to=2-3]
		\end {tikzcd}
	\end{center}

	We may interpret the cospan \( Li' \to  t  \leftarrow  Lj\) as an interaction that acts on a system with interface \(i'\) by the composing it with the system \(L_i \to  s\) along the interaction  \(i + i' \to  m \leftarrow  j\). This idea draws from the Katis-Sabadini-Walters line of thought, considering symmetric monoidal double (or bi-)categories to be \emph{process theories} \cite{katis-1997-bicategories}.\end{remark}

\begin{example}[{Graphical regular logic in the restricted variable sharing doctrine}]\label{ssl-005Y}
\par{}Similar to the examples in \hyperref[ssl-005X]{Example \ref{ssl-005X}}, we can restrict \hyperref[djm-00GP]{the variable sharing doctrine} along the left adjoint \(F \colon  \mathcal {C}\mathsf {at} \to  \mathcal {L}\mathsf {ex}\) that maps a category to its free lex category. Here we give an example a systems theory in this restricted theory that produces a module of system specifications that packages the data of a graphical regular logic  \cite{fong-2018-graphical}.\par{}Let \(\mathsf {C}\) be a lex category and let \(P \colon  \mathsf {C}^{\mathsf {op}} \to  \mathsf {Cat}\) be a regular hyperdoctrine. Its Grothendieck construction \(\int  P\) is a lex category. Let \(S\subseteq  \mathsf {ob}C\) be a set of \emph{sorts}. There is a systems theory in the restricted variable sharing doctrine that consists of the discrete category on \(S\), the lex category \(\int  P\), and the marking that maps the object \(c \in  S\) to the object \((c, \top ) \in  \int  P\) where \(\top \) is the limit of the empty set in \(Pc\). We denote this marking \(L_S \colon  S \hookrightarrow  \int  P\).\par{}This systems theory produces the module of systems \[\mathbb {S}\mathsf {pan}(\int  P)(\bullet , L_S(-)) \colon  \bullet  \mathrel {\mkern 3mu\vcenter {\hbox {$\shortmid $}}\mkern -10mu{\to }} \mathbb {C}\mathsf {ospan}(\mathsf {Finset} \downarrow  S)\] which is the module of \(P\)-predicates over the double category of interfaces \(\mathbb {C}\mathsf {ospan}(\mathsf {Finset} \downarrow  S)\).  Since the cospans which form the interactions of \(\mathbb {C}\mathsf {ospan}(\mathsf {Finset} \downarrow  S)\) are undirected wiring diagrams (as discussed in \hyperref[djm-00G8]{Example \ref{djm-00G8}}), we recover a form of \emph{graphical regular logic} \cite{fong-2018-graphical}.\end{example}

\begin{example}[{Restricting to double categories of directed wiring diagrams}]\label{ssl-005Z}
\par{}In \cite{spivak-2013-operad} Spivak constructs an operad of \emph{directed wiring diagrams} and shows that ODEs form an algebra over this operad. We exemplified in \hyperref[djm-00HF]{Section \ref{djm-00HF}} a module of ODEs over double categories of lenses, and in \hyperref[djm-00G9]{Example \ref{djm-00G9}} we observed that lenses in free cartesian categories are directed wiring diagrams. We can recover this observation by marking they systems theory \(\mathsf {Euc}\) of \hyperref[djm-00IJ]{Example \ref{djm-00IJ}} by \(\mathbb {R}^1 \colon  1 \to  \mathsf {Euc}\) with the real numbers; the resulting module of systems consists of systems of ordinary differential equations that interact via directed wiring diagrams. This observation extends to all the lens-based systems, giving us modules of Moore machines, POMDPs, and so on over double categories of directed wiring diagrams.\end{example}
\printbibliography
\end{document}